\documentclass[preprint,12pt]{elsarticle}

\usepackage{amssymb}
\usepackage{amsmath,graphicx, xcolor,amssymb}
\usepackage[font=small,labelfont=bf,tableposition=top]{caption}
\usepackage{graphicx}
\usepackage{epstopdf}
\usepackage{mathabx}
\usepackage{mathrsfs}
\usepackage{ragged2e}
\usepackage[mathscr]{euscript}
\usepackage{multirow}
\usepackage{subfigure}
\usepackage{enumerate}
\usepackage{mathrsfs}
\usepackage{color,soul}
\usepackage{todonotes}
\usepackage{verbatim}
\usepackage{float}
\usepackage{framed}
\usepackage{algorithm, algorithmic}
\usepackage{lineno}

\begin{document}

\begin{frontmatter}



\title{Interface spaces based on physics for multiscale mixed methods applied to flows in fractured-like porous media}


\author[mymainaddress]{Franciane F. Rocha}
\ead{fr.franciane@usp.br}

\author[mymainaddress]{Fabricio S. Sousa\corref{mycorrespondingauthor}}
\cortext[mycorrespondingauthor]{Corresponding author}
\ead{fsimeoni@icmc.usp.br}

\author[mymainaddress]{Roberto F. Ausas}
\ead{rfausas@icmc.usp.br}

\author[mysecondaryaddress]{Felipe Pereira}
\ead{felipepereira62@gmail.com}

\address[mymainaddress]{Instituto de Ci\^encias Matem\'aticas e de Computa\c c\~ao, Universidade de S\~ao Paulo,\\ Av. Trabalhador S\~ao-carlense, 400, 13566-590, S\~ao Carlos, SP, Brazil}
\address[mysecondaryaddress]{Department of Mathematical Sciences, The University of Texas at Dallas,\\ 800 W. Campbell Road, Richardson, TX 75080-3021, USA}

\author[mymainaddress]{Gustavo C. Buscaglia}
\ead{gustavo.buscaglia@icmc.usp.br}

\begin{abstract}
It is well known that domain-decomposition-based multiscale mixed methods rely on interface spaces, defined on the skeleton of the decomposition, to connect the solution among the non-overlapping subdomains. Usual spaces, such as polynomial-based ones, cannot properly represent high-contrast channelized features such as fractures (high permeability) and barriers (low permeability) for flows in heterogeneous porous media. We propose here new interface spaces, which are based on physics, to deal with permeability fields in the simultaneous presence of fractures and barriers, accommodated respectively, by the pressure and flux spaces. 
Existing multiscale methods based on mixed formulations can take advantage of the proposed interface spaces, however, in order to present and test our results, we use the newly developed Multiscale Robin Coupled Method (MRCM) [Guiraldello, et al., \textit{J. Comput. Phys.},  355 (2018) pp. 1-21], which generalizes most well-known multiscale mixed methods, and allows for the independent choice of the pressure and flux interface spaces. 
An adaptive version of the MRCM [Rocha, et al., \textit{J. Comput. Phys.}, 409 (2020), 109316] is considered that automatically selects the physics-based pressure space for fractured structures and the physics-based flux space for regions with barriers, resulting in a procedure with unprecedented accuracy. The features of the  proposed approach are investigated through several numerical simulations of single-phase and two-phase flows, in different heterogeneous porous media. The adaptive MRCM combined with the interface spaces based on physics provides promising results for challenging problems with the simultaneous presence of fractures and barriers.
\end{abstract}




\end{frontmatter}




\section{Introduction}
Simulations of petroleum reservoirs deal with highly heterogeneous permeability fields with multiple scales and high contrast. A representation of the solution that captures the heterogeneity requires several billion cells, making the numerical simulations extremely expensive \cite{chen2006computational}. The so-called {\em multiscale methods} have been introduced to exploit the multiscale structure of the problem and provide efficient approximations at a reduced computational cost \cite{kippe}. 
The multiscale literature is vast and contemplates numerical implementations using techniques such as the finite volume method \cite{jenny2003multi,jenny2005adaptive}, the primal finite element method \cite{HouMulti, aarnes2002multiscale, efendiev2013generalized}, and several mixed finite element methods \cite{chen_hou, aarnes, arbogast, valentin, pereira, chung2015mixed}. 

The multiscale approaches solve the problem on a coarse decomposition of the domain, incorporating the fine-grid information through local basis functions. The accuracy of the multiscale methods is related to the calculation of the basis functions. If the heterogeneities are not well represented by the basis, inaccurate solutions are obtained.
Spaces that are polynomial on the interfaces of the domain decomposition work well for smooth or Gaussian permeability fields, but their performance for high contrast, channelized ones is not satisfactory \cite{chung2014adaptive,chung2016enriched,cortinovis2014iterative,guiraldello2019interface}. 
To remedy this, informed spaces obtained from sets of snapshots by algebraic dimensionality reduction were considered in \cite{efendiev2013generalized,chung2016enriched,guiraldello2019interface}. Another approach is to define the interface space through eigensolutions of local partial-differential problems \cite{madureira2017hybrid}. These approaches can also be coupled to a-posteriori error estimators \cite{chung2014adaptive, chung2017online, chung2018online}. 

Other authors have looked for approaches more directly based on the geometry of the heterogeneities. In \cite{aarnes2006hierarchical, moyner2016multiscale, klemetsdal2019accelerating}, for example, the authors consider polynomial bases but adapt the grid to the geological properties, while in \cite{peszynska2002mortar,cortinovis2014iterative,cortinovis2017zonal} local enrichment functions were added on high-permeability regions. These geometrical strategies are well suited for permeability fields containing either channels or barriers, but not both. In this contribution, we propose a strategy to deal with the simultaneous presence of channels and barriers, since such situation is not unfrequent in reservoir simulation \cite{hosseinimehr2019dynamic}.
For this purpose, we take advantage of the Multiscale Robin Coupled Method (MRCM, \cite{guiraldello2018multiscale}), which allows for the independent choice of the pressure and flux interface spaces. In fact, the pressure space is designed so as to accommodate channels and the flux space to accommodate barriers, and the adaptive version of the MRCM \cite{bifasico} is used to automatically select the appropriate parameters at each location.

Our intention is to use high-definition volumetric grids that capture the large scale features of the fracture network, especially when they are relatively large as happens in fractured karst reservoirs \cite{baomin2009classification, popov2009multiphysics, huang2013efficient,lopes2019new}. Typically, fractures are handled with separate discrete models that represent the fractures as lower-dimensional objects so as to incorporate sub-grid resolution \cite{martin2005modeling, reichenberger2006mixed, formaggia2014reduced, schwenck2015dimensionally, berkowitz2002characterizing, flemisch2018benchmarks}. A popular approach is the Discrete Fracture Model (DFM) which uses unstructured grids to place fractures at the interface between matrix cells \cite{kim2000finite, karimi2003efficient, hoteit2008efficient}. Other approaches are the Embedded DFM \cite{ctene2016algebraic,chai2018efficient}, the multi-continuum model \cite{chung2017coupling, wang2020generalized} and the hierarchical fracture models  \cite{hajibeygi2011hierarchical, li2008efficient, efendiev2015hierarchical}.

Since it is possible to combine discrete models with multiscale methods \cite{bosma2017multiscale,  zhang2017multiscale, devloo2019multiscale, xia2018enriched}, the final goal of the improved MRCM proposed here is to allow a unified treatment of fractured karst reservoirs in which the modeling of the fractures is shared, depending on the fracture's size, between the volumetric grid and the discrete models. For this reason, we consider here permeability fields containing multiple narrow and relatively straight features (channels, barriers) that mimic the largest structures of a fractured porous medium, and refer to them as ``fractured-like'' fields. The accuracy and efficiency of the proposed method are investigated through several numerical simulations of fractured-like reservoirs both in single-phase and two-phase (oil-water) situations.

The paper is organized as follows. In Section \ref{Model problem} we introduce the model problem and numerical approximation. The MRCM is briefly recalled in Section \ref{The multiscale Robin coupled method} followed by Section \ref{physics-basead interface spaces}, where the strategies for building the physics-based interface spaces are presented along with some numerical experiments. In section \ref{aMRCM-PBS} we explain our proposed combination of the physics-based spaces with the adaptive MRCM.
Numerical experiments for single and two-phase flows are presented in Section \ref{results}. Finally, the paper is concluded in Section \ref{conclusions}.

\section{Model problem}\label{Model problem}

We consider a high-contrast heterogeneous oil reservoir through which an immiscible and incompressible single-phase or multiphase flow takes place. The heterogeneity of the medium is represented by a space-dependent permeability coefficient in the elliptic Darcy model. This coefficient exhibits variations of many orders of magnitude over short distances. For simplicity, the capillary pressure and gravity effects are not considered. 

\subsection{Single-phase model}

The unknowns considered in the single-phase model are the Darcy velocity $\mathbf{u}(\mathbf x, t)$ and the fluid pressure $p(\mathbf x, t)$.
The governing equations are the Darcy's law with a statement of conservation of mass, given by
 \begin{equation}\label{Darcy0}
\begin{array}{rll}
\mathbf{u}&=-K(\mathbf{x})\nabla p  &\mbox{in}\ \Omega \\
\nabla \cdot \mathbf{u}&=q  &\mbox{in}\ \Omega \\ 
p &= g &\mbox{on}\ \partial\Omega_{p}\\
\mathbf{u} \cdot \mathbf{n}&= z &\mbox{on}\ \partial\Omega_{u}
\end{array}
\end{equation}
where $\Omega\subset\mathbb{R}^2 $ is the computational domain; $K(\mathbf{x})$ is the symmetric, uniformly positive definite absolute permeability tensor; $q = q(\mathbf x, t)$ is a source term; $g = g(\mathbf x, t)$ is the pressure condition specified at the boundary $\partial\Omega_{p}$ and $z = z(\mathbf x, t)$ is the normal velocity condition ($\mathbf{n}$ is the outward unit normal) specified at the boundary $\partial\Omega_{u}$.

\subsection{Two-phase model}

In the context of two-phase flows of water and oil (denoted by $w$ and $o$, respectively), we consider that the reservoir contains injection wells, from which water is injected to displace the trapped oil towards production wells. The oil and water saturations are the new unknowns. These saturations are related considering a fully saturated medium (the sum of both oil and water saturation is equal to one). Therefore the model considers only the water saturation $s(\mathbf x, t)$ in the transport problem. Here the elliptic problem (\ref{Darcy0}) is written as
\begin{equation}\label{Darcy}
\begin{array}{rll}
\mathbf{u}&=-\lambda(s)K(\mathbf{x})\nabla p  &\mbox{in}\ \Omega \\
\nabla \cdot \mathbf{u}&=q  &\mbox{in}\ \Omega \\ 
p &= g &\mbox{on}\ \partial\Omega_{p}\\
\mathbf{u} \cdot \mathbf{n}&= z &\mbox{on}\ \partial\Omega_{u}
\end{array}
\end{equation}
and coupled with the hyperbolic conservation law for the transport of water saturation \cite{peaceman2000fundamentals,ewing1983mathematics}, which is given by
\begin{equation}\label{BL2D}
\begin{array}{rll}
\dfrac{\partial s}{\partial t} + \nabla \cdot \left(f(s)\mathbf{u}\right)& = 0  &\mbox{in}\ \Omega \\
s(\mathbf{x},t=0) &= s^0(\mathbf{x}) &\mbox{in}\ \Omega\\
s(\mathbf{x},t) &= \bar{s}(\mathbf{x},t) &\mbox{in}\ \partial\Omega^-
\end{array}
\end{equation}
where $s^0$ is the initial condition for the saturation and $\bar{s}$ is the saturation at the inlet boundaries $\partial\Omega^-=\{\mathbf{x}\in\partial\Omega,\ \mathbf{u} \cdot \mathbf{n}<0\}$.
The total mobility $\lambda(s)=\lambda_w(s)+\lambda_o(s)$ and the fractional flow of water $f(s)$ are respectively given by
\begin{equation}
 \lambda(s) =\dfrac{k_{rw}(s)}{\mu_w}+\dfrac{k_{ro}(s)}{\mu_o}\qquad  \mbox{and} \qquad f(s) = \dfrac{\lambda_w(s)}{\lambda(s)},
 \label{fluxos} 
\end{equation}
where $k_{rj}(s)$ and $\mu_j$, $j \in \{w,o\} $, are respectively the relative permeability function and viscosity of phase $j$. For simplicity, the porosity is considered constant both in space and time and consequently, it can be easily scaled out by changing the time variable. 

\subsection{Numerical approximation}

For computing single phase flows we solve directly the elliptic problem (\ref{Darcy0}) by the Multiscale Robin Coupled Method (MRCM). To simulate two-phase flows we solve the system (\ref{Darcy})-(\ref{BL2D}) using an operator splitting procedure, in which we solve (\ref{Darcy}) to compute $p(\mathbf x, t)$ and $\mathbf{u}(\mathbf x, t)$ and (\ref{BL2D}) for $s(\mathbf x, t)$, sequentially \cite{douglas1997numerical}. 
We use the MRCM to solve the elliptic problem (\ref{Darcy}), and the upwind method \cite{leveque2002finite} for the transport of saturation (\ref{BL2D}), in which a CFL-type condition is enforced. The splitting procedure considers larger time steps for the elliptic problem compared to those used for the hyperbolic equation (\cite{douglas1997numerical}) improving the computational efficiency. 

We consider the relation $\Delta t_p = C \Delta t_s$, where $\Delta t_s$ is the time step for the saturation equation (for simplicity assumed constant but in practice we allow for variable $\Delta t_s$), $\Delta t_p$ is the time step for the elliptic problem and $C$ is a positive integer.
The elliptic problem is solved at times $t_n = n\Delta t_p$, for $n=0,1,\dots$, where we compute $p^{n}(\mathbf x)$ and $\mathbf u^{n}(\mathbf x)$ using the water saturation $s^{n}(\mathbf x)$ at time $t=t_n$. The saturation equation is solved at intermediate times $t_{n,k} = t_n + k\Delta t_s$, for $k=1,2,\dots, C$, such that $t_n < t_{n,k} \leq t_{n+1}$. 
For each saturation time step $t = t_{n,k}$ we set the velocity $\mathbf u^{n,k}$ to the following linear extrapolated value
\begin{equation}\label{extrapolation}
\mathbf u^{n,k}(\mathbf x, t) =
\left\{
\begin{array}{ll}
\mathbf u^0(\mathbf x), &   \text{ if } 0 \leq t\leq t_1, \\
\dfrac{t-t_{n-1}}{\Delta t_p}\ \mathbf u^n(\mathbf x) -\dfrac{t-t_n}{\Delta t_p}\ \mathbf u^{n-1}(\mathbf x),
& \text{ if } t_n<t\leq t_{n+1},\\ 
\end{array} \right.
\end{equation}
where $\mathbf u^n (\mathbf x)$ is the velocity obtained from Equation (\ref{Darcy}) at time $t=t_n$, approximating $\mathbf u(\mathbf x, t_n)$ (see  \cite{douglas1997numerical}).

\section{The multiscale Robin coupled method}\label{The multiscale Robin coupled method}

In this paper, we use the MRCM to approximate the Equation (\ref{Darcy0}) for single-phase flows and the Equation (\ref{Darcy}) for two-phase flows. For simplicity, we drop the time dependency and obtain the following elliptic problem 
\begin{equation}\label{MRCM}
\begin{array}{rll}
\mathbf{u}&=-\kappa(\mathbf{x})\nabla p  &\mbox{in}\ \Omega \\
\nabla \cdot \mathbf{u}&=q  &\mbox{in}\ \Omega \\ 
p &= g &\mbox{on}\ \partial\Omega_{p}\\
\mathbf{u} \cdot \mathbf{n}&= z &\mbox{on}\ \partial\Omega_{u}
\end{array}
\end{equation}
where $\kappa=K(\mathbf x)$ for single-phase and $\kappa = \lambda(s(\mathbf x)) K(\mathbf x)$ for two-phase flows.

The MRCM considers a decomposition of the domain $\Omega$ into non-overlapping subdomains $\Omega_i, \ i=1,2,\cdots,N$. Continuity of the normal fluxes and pressure are weakly imposed on a coarse scale $H$, which is significantly larger than the fine-scale of the discretization $h$. Here $H$ is the characteristic size of the subdomains.
The weak continuities are enforced to the multiscale solution $(\mathbf u_h,p_h)$ at the skeleton $\Gamma$ of the decomposition (union of all interfaces $\Gamma_{ij} = \Omega_i \cap \Omega_j$) by the following compatibility conditions 
\begin{equation}
\int_\Gamma (\mathbf u_h^+ - \mathbf u_h^-)\cdot \check{\mathbf n} \ \psi \ d\Gamma=0 \quad \text{and}\quad 
\int_\Gamma (p_h^+ - p_h^-)\ \phi \ d\Gamma=0
\label{eq:compatibility}
\end{equation}
for all $(\phi,\psi) \in \mathcal{U}_H \times \mathcal{P}_H$, which are the interface spaces defined over the edges $\mathcal{E}_h$ of the skeleton $\Gamma$. 
Here $\check{\mathbf n}$ is a fixed normal vector to the skeleton $\Gamma$, pointing outwards from the sudomain with the smallest index. The solution on each side of the interface $\Gamma$ is represented by the $+$ and $-$ superscripts. The solution inside each subdomain
$(\mathbf u_h^i,p_h^i)$ is related to the normal flux and pressure unknowns at the interfaces $(U_H,P_H) \in \mathcal{U}_H \times \mathcal{P}_H$ by the Robin boundary conditions on the local problems
\begin{equation}
-\frac{\alpha(\mathbf x) H}{\kappa_i(\mathbf x)} \ \mathbf u_h^i \cdot \check{\mathbf n}^i + p_h^i = 
-\frac{\alpha(\mathbf x) H}{\kappa_i(\mathbf x)}\  U_H \ \check{\mathbf n} \cdot \check{\mathbf n}^i + P_H, \quad \mathbf x \in \Gamma_{i,j},
\label{eq:robin}
\end{equation}
where the $\check{\mathbf n}^i$ is the normal vector to $\Gamma$ pointing outwards of $\Omega_i$.
The parameter for the Robin boundary condition in Equation (\ref{eq:robin}) is given by
\begin{equation}
\beta_i(\mathbf {x}) = \frac{\alpha(\mathbf{x}) H}{\kappa_i(\mathbf x)},
\label{eq:beta}
\end{equation}
where $\alpha(\mathbf {x})$ is a dimensionless algorithmic function that takes different values according to the permeability field (in the context of the adaptive MRCM).
We remark that by changing $\alpha$ we recover the Multiscale Mortar Mixed Finite Element Method (MMMFEM, \cite{arbogast}) when $\alpha \rightarrow 0$ and the Multiscale Hybrid-Mixed Finite Element Method (MHM, \cite{valentin}) when $\alpha \rightarrow +\infty$. The solution given by the Multiscale Mixed Method (MuMM, \cite{pereira}) can also be recovered under the right choice of parameters. The interested reader is referred to \cite{guiraldello2018multiscale} for more details. 

The differential formulation of the MRCM consists in finding local solutions $(\mathbf u_h^i, p_h^i)$ into each subdomain $\Omega_i$, and global interface unknowns $(U_H,P_H)$ satisfying the local problems
\begin{equation}
\begin{array}{rcll}
\mathbf u_h^i &=& -\kappa(\mathbf{x)} \ \nabla p_h^i & \text{in }\Omega_i \\
\nabla \cdot \mathbf u_h^i &=& q & \text{in }\Omega_i \\
p_h^i & = & g & \text{on }\partial\Omega_i\cap \partial\Omega_p \\
\mathbf u_h^i \cdot \check{\mathbf n}^i &=& z & \text{on }\partial\Omega_i\cap \partial\Omega_u \\
-\beta_i  \mathbf u_h^i \cdot \check{\mathbf n}^i + p_h^i &=&  -\beta_i U_H \check{\mathbf n} \cdot \check{\mathbf n}^i + P_H & \text{on }\partial\Omega_i\cap \Gamma 
\end{array}
\label{eq:mrc1}
\end{equation}
and the compatibility conditions
\begin{equation}
\begin{array}{rcl}
\displaystyle \sum_{i=1}^N \int_{\partial\Omega_i\cap \Gamma} (\mathbf u_h^i \cdot \check{\mathbf n}^i) \ \psi \ d\Gamma &=& 0\\
\displaystyle \sum_{i=1}^N \int_{\partial\Omega_i\cap \Gamma} \beta_i (\mathbf u_h^i \cdot \check{\mathbf n}^i - U_H\ \check{\mathbf n} \cdot \check{\mathbf n}^i)\ \phi  \ (\check{\mathbf n} \cdot \check{\mathbf n}^i)\ d\Gamma &=& 0
\end{array}
\label{eq:mrc2}
\end{equation}
for all $(\phi,\psi)\in \mathcal{U}_H\times\mathcal{P}_H$.

The interface spaces $\mathcal{U}_H$ and $\mathcal{P}_H$ are defined over the skeleton $\Gamma$ as subspaces of
\begin{equation}
F_h(\mathcal{E}_h) = \left\{ f:\mathcal{E}_h\to \mathbb{R};~f|_e\,\in\,\mathbb{P}_0~,~\forall\,e\,\in\,\mathcal{E}_h \right\} ~,
\end{equation}
where $\mathcal{E}_h$ is the set of all edges of the skeleton $\Gamma$ \cite{guiraldello2018multiscale}.
These spaces are local and independently built on each interface $\Gamma_{i,j}$ and spanned by the multiscale basis functions $\{\phi_1, \phi_2, \cdots , \phi_{N_U}\}$ and $\{\psi_1, \psi_2, \cdots , \psi_{N_P}\}$, where $N_U$ and $N_P$ are the respective dimensions of the interface spaces (see \cite{guiraldello2018multiscale} for more details). The multiscale basis functions are obtained by solving the local problems in Equations (\ref{eq:mrc1}), that can be computed in parallel. These basis functions are further used to compute the global solution as a linear combination, whose coefficients are obtained by the solution of the global interface system generated by Equations (\ref{eq:mrc2}).

In the previous works, the authors consider the interface spaces $\mathcal{U}_H$ and $\mathcal{P}_H$ as low-dimensional polynomial spaces \cite{guiraldello2018multiscale} or informed spaces \cite{guiraldello2019interface} for the MRCM. In the next section, we propose alternative physics-based interface spaces to improve the accuracy of the multiscale solution of the elliptic problem. These interface spaces are built to capture variations of the permeability field, especially those characterized by fractures and barriers.

\section{Physics-based interface spaces} \label{physics-basead interface spaces}
Low-dimensional polynomial spaces are not robust to represent variations of high-contrast permeability fields containing channelized structures as fractures and barriers. Even informed spaces with a low number of degrees of freedom are not enough. To better represent the variations of fractured-like permeability fields we present two physics-based interface spaces. The idea is to build the multiscale basis functions based on the pressure and flux solutions at each channelized structure (fracture/barrier). We present one space for pressure and another for flux. The new pressure space is suited for high permeability channels, whereas the new flux space for barriers or low permeability barriers, that cross the interface between subdomains.

The support of the new basis functions are the interfaces that contain channelized structures. In the remaining interfaces the spaces can be freely chosen (any low-dimensional polynomial or informed spaces). 
We denote $(\mathcal{U}_{H,k}, \mathcal{P}_{H,k})$ the choice of the interface spaces made up of the elementwise constant fine grid representation of polynomials over the interface elements, where $k$ is the degree of the polynomial. To introduce the concept we consider, for simplicity, the linear spaces $(\mathcal{U}_{H,1}, \mathcal{P}_{H,1})$ for both pressure and flux. This choice of interface spaces is used further in the numerical section. The basis of these spaces contain $N_P=2$ and $N_U=2$ functions on each interface between two subdomains.
If the permeability field contains fractures passing through the interface we propose to replace the space $\mathcal{P}_{H,1}$ by the physics-based pressure space. On the other hand, if the permeability filed contains barriers crossing the interface we replace the space $ \mathcal{U}_{H,1}$
by the physics-based flux space. In the following subsections, we describe how the physics-based interface spaces are built for capturing the high-contrast structures.


\subsection{A physics-based interface space for the pressure} \label{A physics-basead interface space for the pressure}

\begin{figure}[htbp]
\center 
\subfigure[ref1][Permeability]{\includegraphics[scale=0.178]{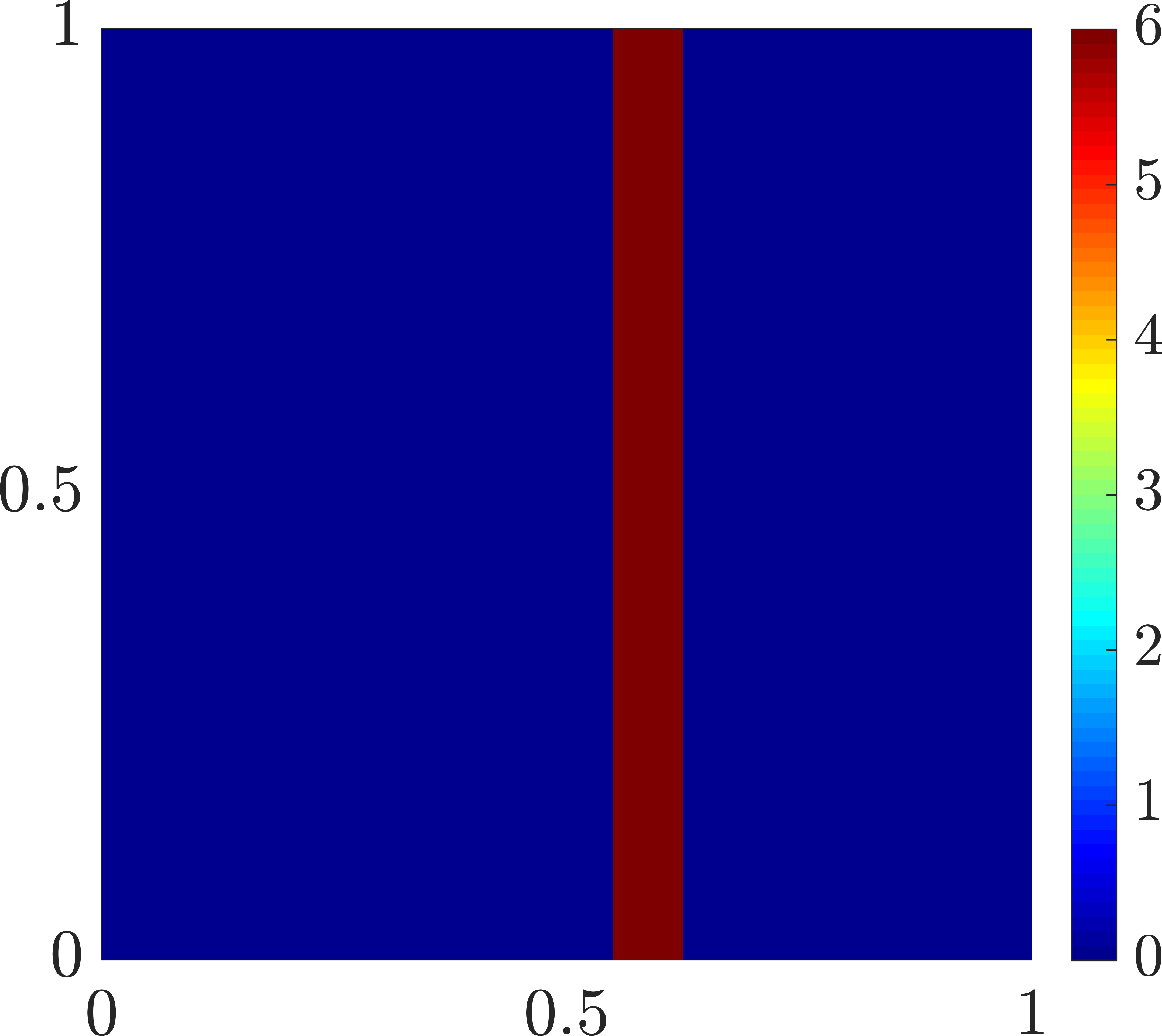}} \label{fig:frac1_a}
\subfigure[ref2][Pressure]{\includegraphics[scale=0.178]{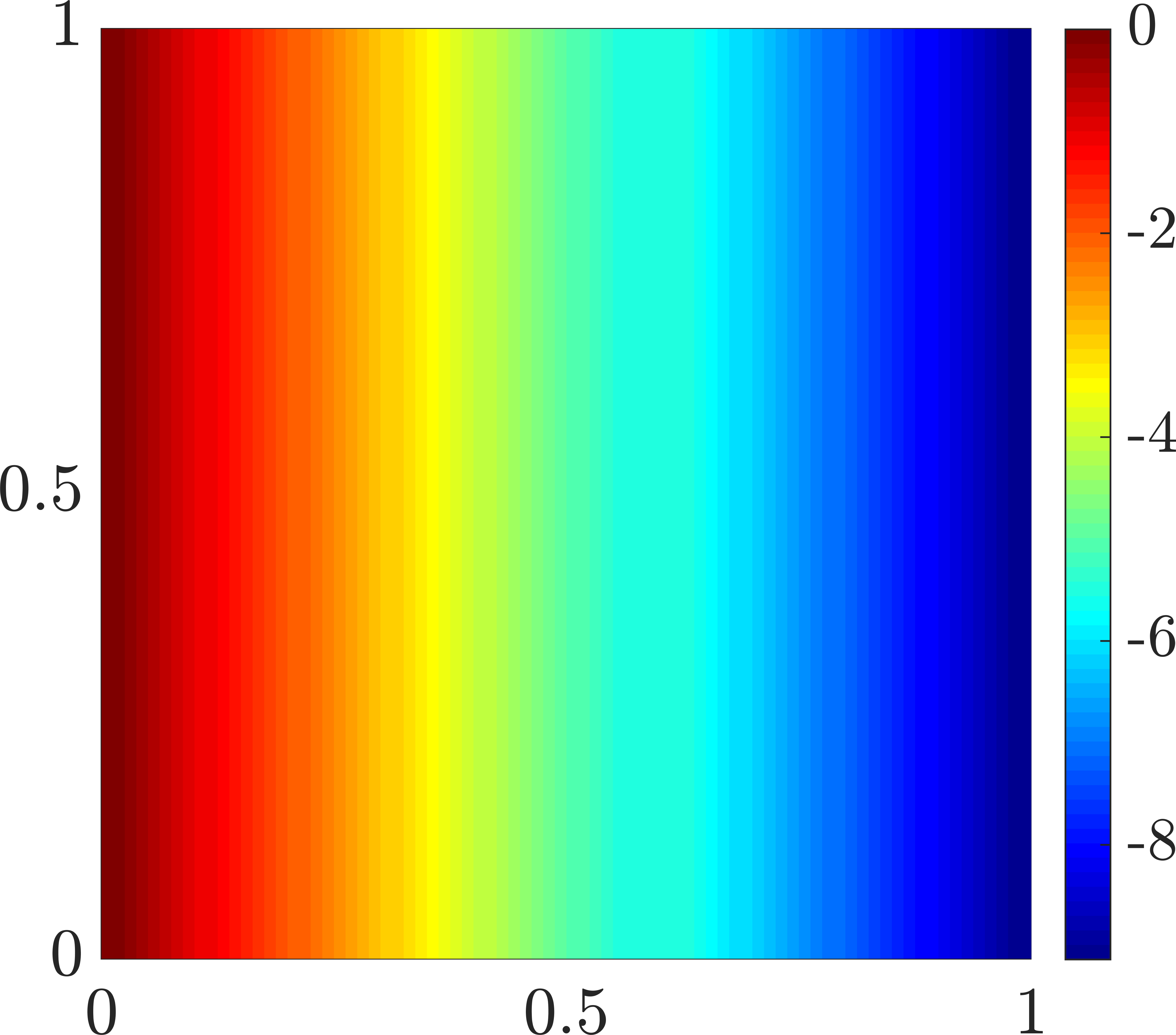}} \label{fig:frac1_b}
\subfigure[ref2][Flux]{\includegraphics[scale=0.178]{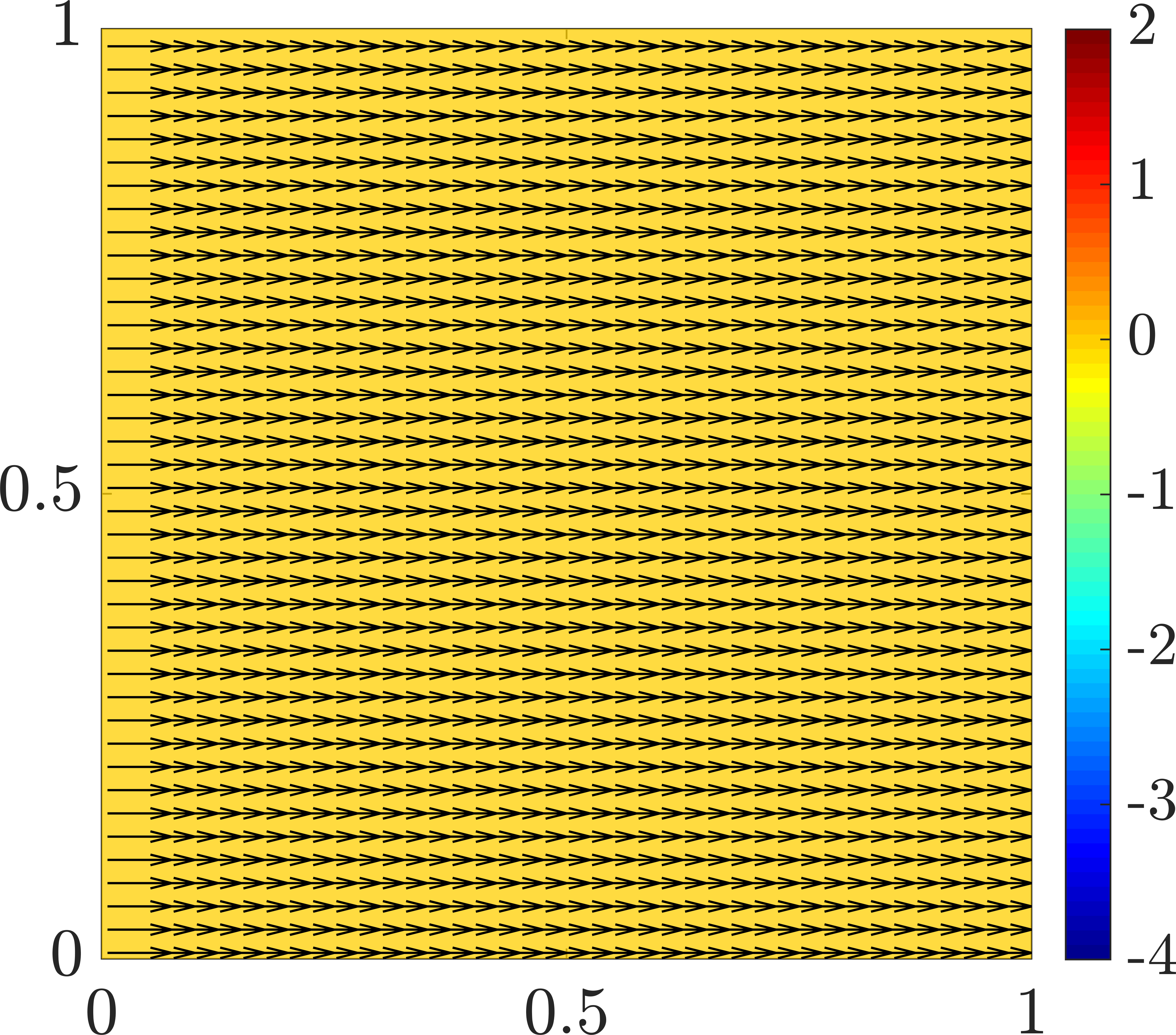}} \label{fig:frac1_c}
\subfigure[ref2][Pressure profile]{\includegraphics[scale=0.178]{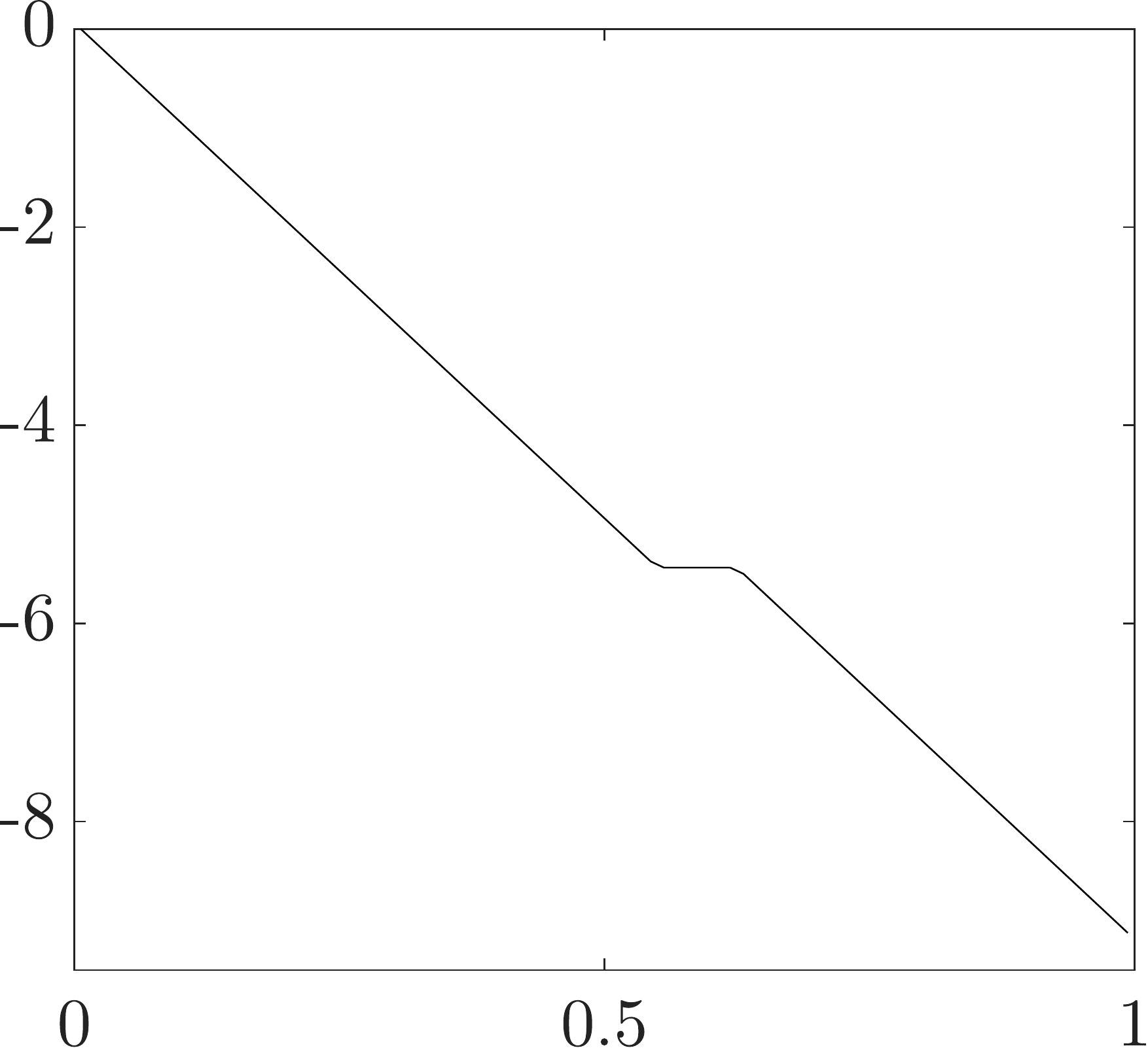} } \label{fig:frac1_d}
\caption{Vertical fracture problem. (a) Permeability field (log-scaled) containing a vertical fracture. (b) Pressure field. (c) Flux. The colors in the flux plot refer to the log-scale flux magnitude. (d) Pressure along a horizontal line.
}
 \label{fig:frac1}
\end{figure}

We define a simplified problem test to show the behavior of the solution to define the interface spaces.
In Figure \ref{fig:frac1} we consider a high-contrast permeability field containing a vertical fracture Figure \ref{fig:frac1}(a). We show the pressure Figure \ref{fig:frac1}(b) and flux Figure \ref{fig:frac1}(c) fine grid solutions. The flow is established by imposing a flux boundary condition from left to right and no-flow at top and bottom. This geometry induces a one-dimensional pressure solution which is plotted along a horizontal line in Figure \ref{fig:frac1}(d). Notice, the pressure is essentially constant over the fracture region.
Any domain decomposition with more than one subdomain in the $y$-direction contains horizontal interfaces through which the fracture passes. Let $\Gamma^{\text{frac}}$ be the set of all interfaces that contain at least one fine cell in which the absolute permeability is larger than a cutoff value $\zeta_{\max}$. 
We intend to replace the pressure linear space $\mathcal{P}_{H,1}=\text{span}\{\psi_1, \psi_2\}$ at the interfaces $\Gamma_{i,j}\subset\Gamma^{\text{frac}}$ by a physics-based pressure space, denoted by $\mathcal{P}_{H}^*$.

\begin{figure}[htbp]
    \centering 
    \includegraphics[scale=0.35]{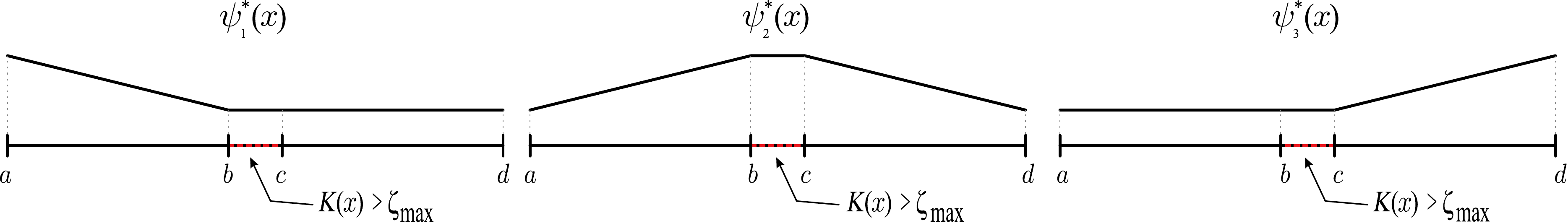} 
    \caption{Physics-based basis functions for pressure at the interfaces that contain fractures. 
Note that the set of functions is able to capture the pressure profile across the fracture.
}
    \label{fig:new_press_space}
\end{figure}

The new pressure space mimics the behavior of the pressure across the fractures. Let $\Gamma_{i,j}\subset\Gamma^{\text{frac}}$ be an interface with support in $[a,d]$ through which a fracture passes in $[b,c]\subset[a,d]$, as sketched in Figure \ref{fig:new_press_space}. The basis functions are defined as:
\begin{equation}\label{new_basis1}
\psi_1^*(x)=\left\{
\begin{array}{cl}
\dfrac{b-x}{b-a}  &\mbox{if}\ x\in(a,b) \\
0\   &\mbox{otherwise}
\end{array} \right. 
\end{equation}
\begin{equation}\label{new_basis2}
\psi_2^*(x)=\left\{
\begin{array}{cl}
\dfrac{x-a}{b-a}  &\mbox{if}\ x\in(a,b) \\
1  &\mbox{if}\ x\in(b,c) \\
\dfrac{d-x}{d-c}  &\mbox{if}\ x\in(c,d) 
\end{array} \right.
\end{equation}
\begin{equation}\label{new_basis3}
\psi_3^*(x)=\left\{
\begin{array}{cl}
\dfrac{x-c}{d-c}  &\mbox{if}\ x\in(c,d) \\
0\   &\mbox{otherwise}
\end{array} \right.
\end{equation}
The new interface pressure space is then defined as $\mathcal{P}_{H}^*=\text{span}\{\psi_1^*, \psi_2^*, \psi_3^*\}$ at $\Gamma_{i,j}\subset\Gamma^{\text{frac}}$. Since the definition of the basis depends only on the fine-grid discretization of the permeability at the interface, any fracture crossing the interface can be represented. Therefore, these basis functions are not restricted to fractures orthogonal to the interface. If the interface contains more than one fracture we need to define a new basis function with similar behavior to the $\psi_2^*$ for each fracture. The total number of basis functions per interface is thus $2+N_{\text{frac}}$, where $N_{\text{frac}}$ is the number of fractures at the interface.
 
\begin{figure}[htbp]

\center
\subfigure[ref1][Permeability]{\includegraphics[scale=0.172]{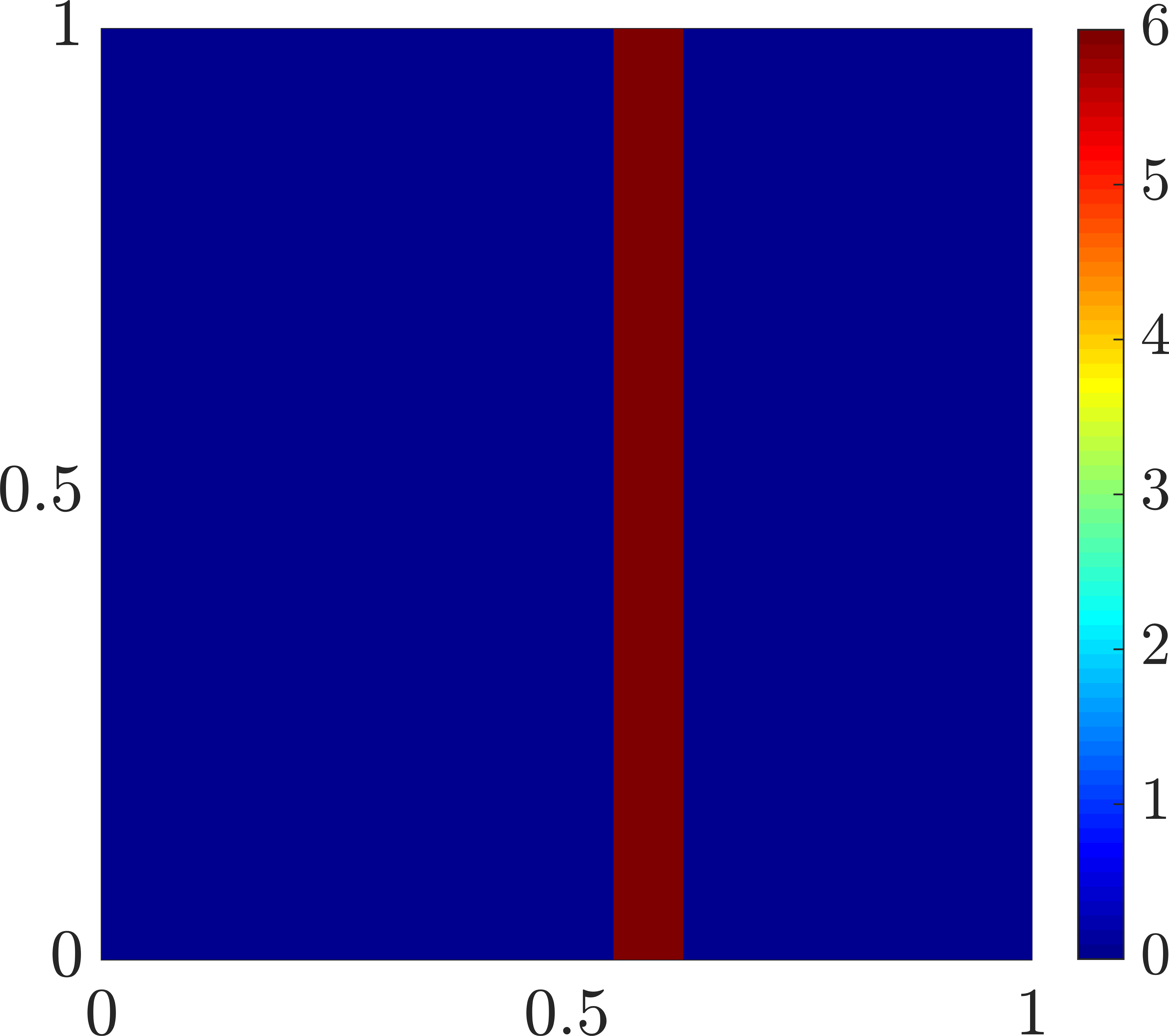} } \label{fig:frac2_a}
\subfigure[ref2][Pressure]{\includegraphics[scale=0.172]{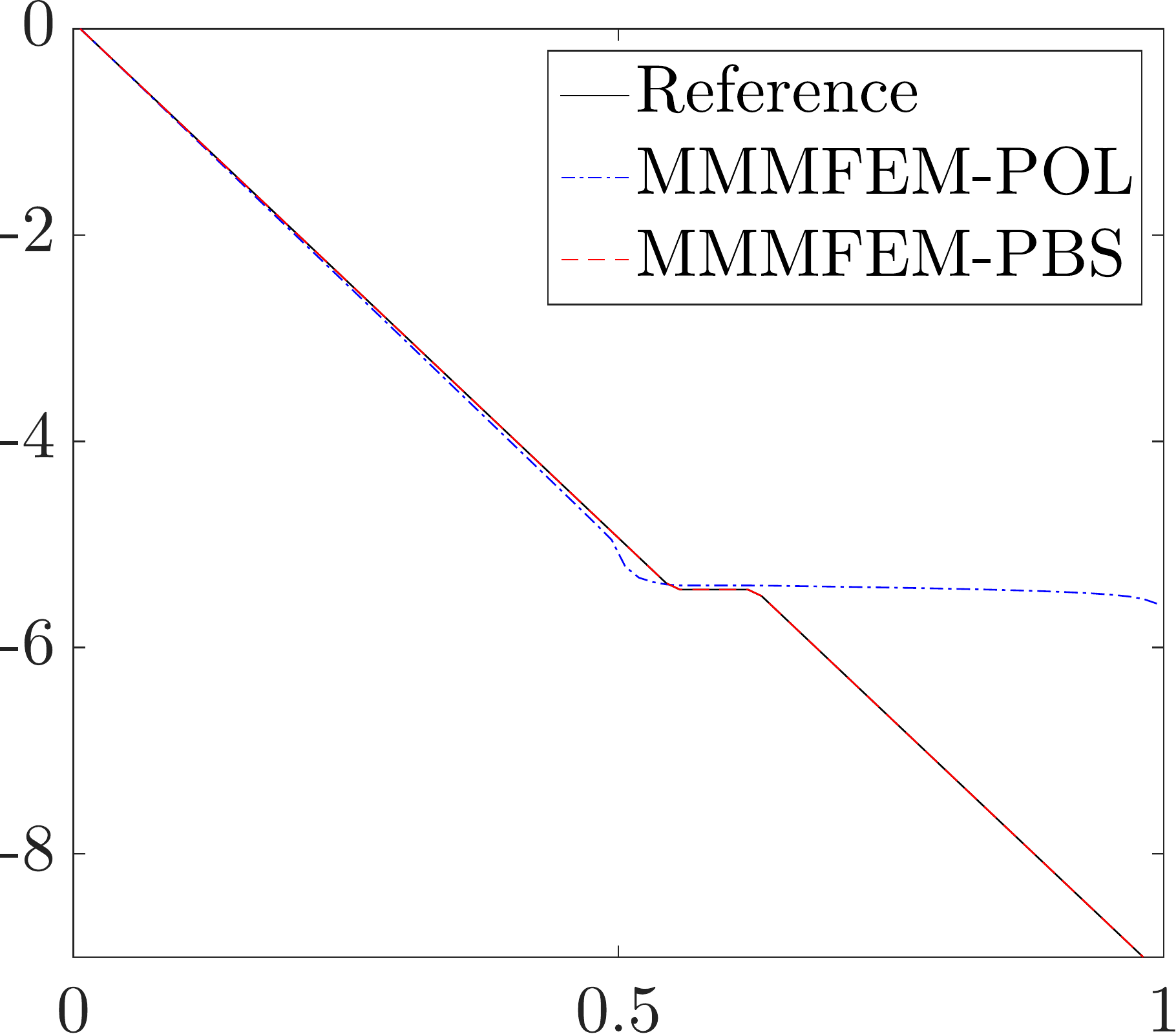} } \label{fig:frac2_b}
\subfigure[ref2][Permeability]{\includegraphics[scale=0.172]{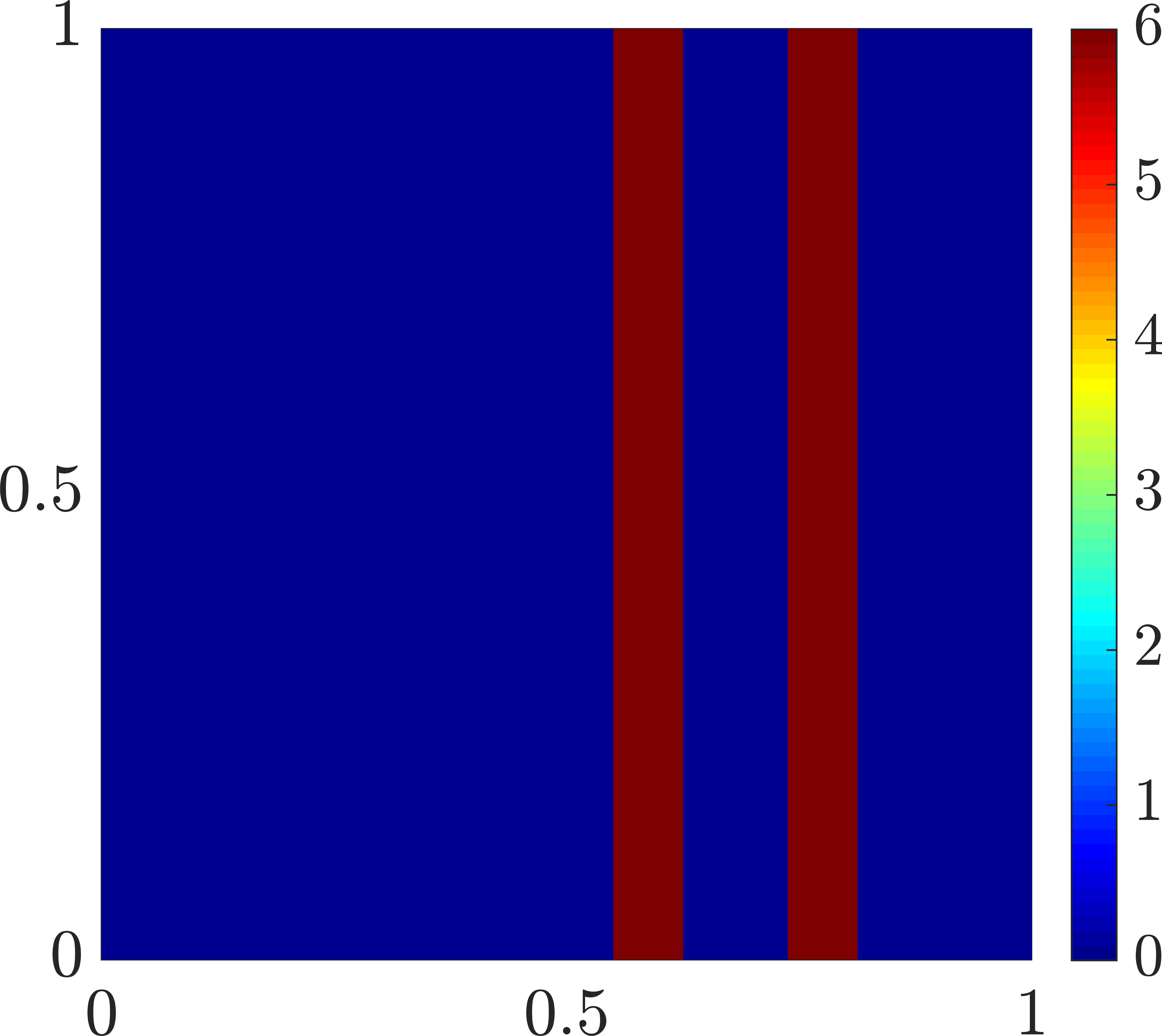} } \label{fig:frac2_c}
\subfigure[ref2][Pressure]{\includegraphics[scale=0.172]{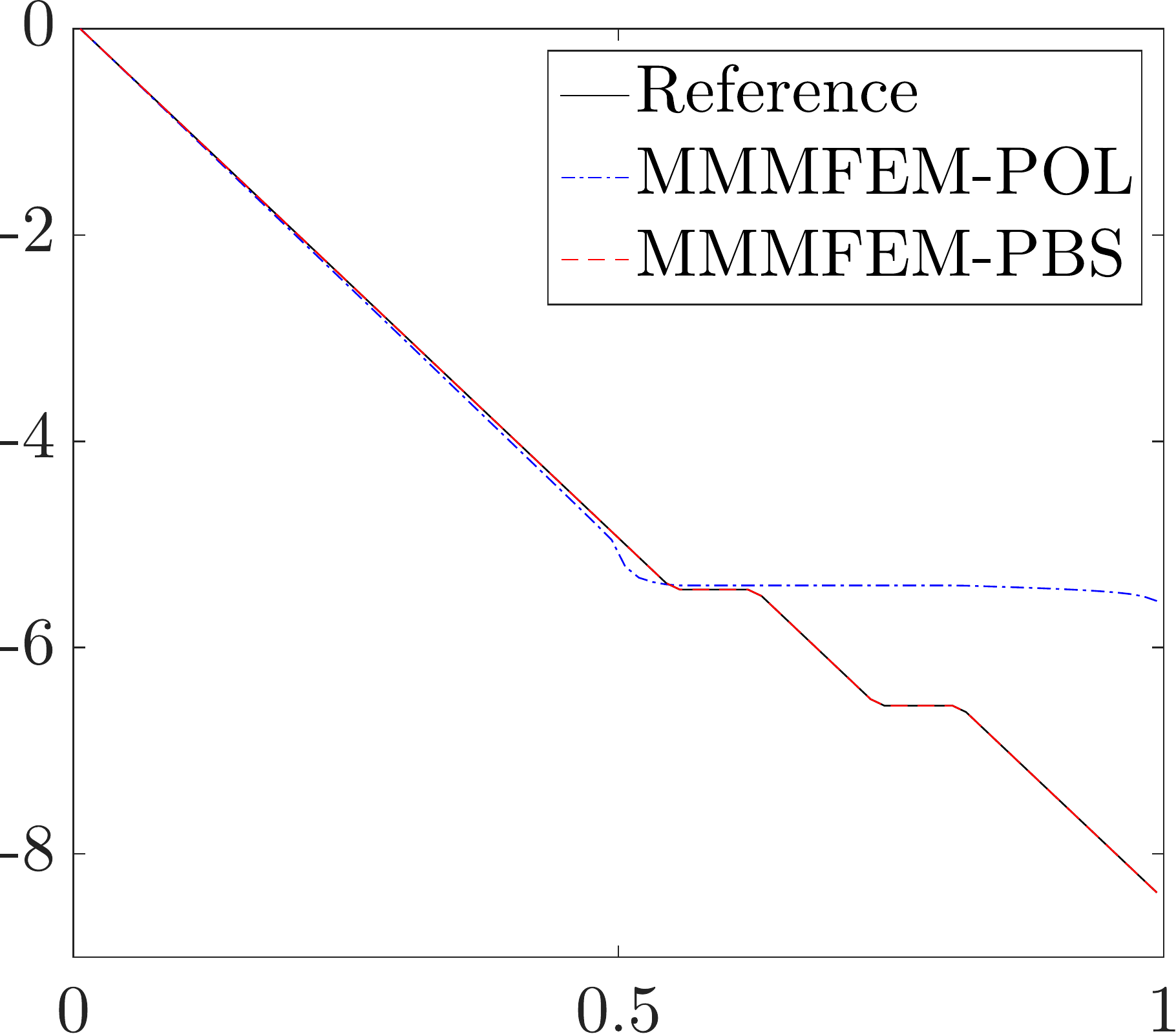} } \label{fig:frac2_d}
\caption{High-contrast permeability field (log-scaled) with one (a) and two fractures (c). Pressure solutions for one (b) and two fractures (d) computed by the fine-grid solver and the multiscale ones MMMFEM-POL and MMMFEM-PBS.  We note that for both one and two fractures the correct pressure solution is only captured by the MMMFEM-PBS.}
 \label{fig:frac2}
\end{figure}

In Figure \ref{fig:frac2}, we show the pressure solution at the horizontal line $y=0.49$ delivered by the MMMFEM (by setting the $\alpha$ parameter of the MRCM to the value $10^{-6}$) in a domain decomposition with $2\times2$ coarse cells. We denote by MMMFEM-PBS the multiscale method combined with the physics-based spaces and by MMMFEM-POL the multiscale method combined with polynomial spaces (linear in this paper). We note that for both single (Figure \ref{fig:frac2}(a) and Figure \ref{fig:frac2}(c)) and multiple (Figure \ref{fig:frac2}(b) and Figure \ref{fig:frac2}(d)) fractures the correct pressure solution is only captured by the MMMFEM combined with physics-based interface spaces.

\subsection{A physics-based interface space for the flux} \label{A physics-basead interface space for the flux}

\begin{figure}[htbp]
\center
\subfigure[ref1][Permeability]{\includegraphics[scale=0.173]{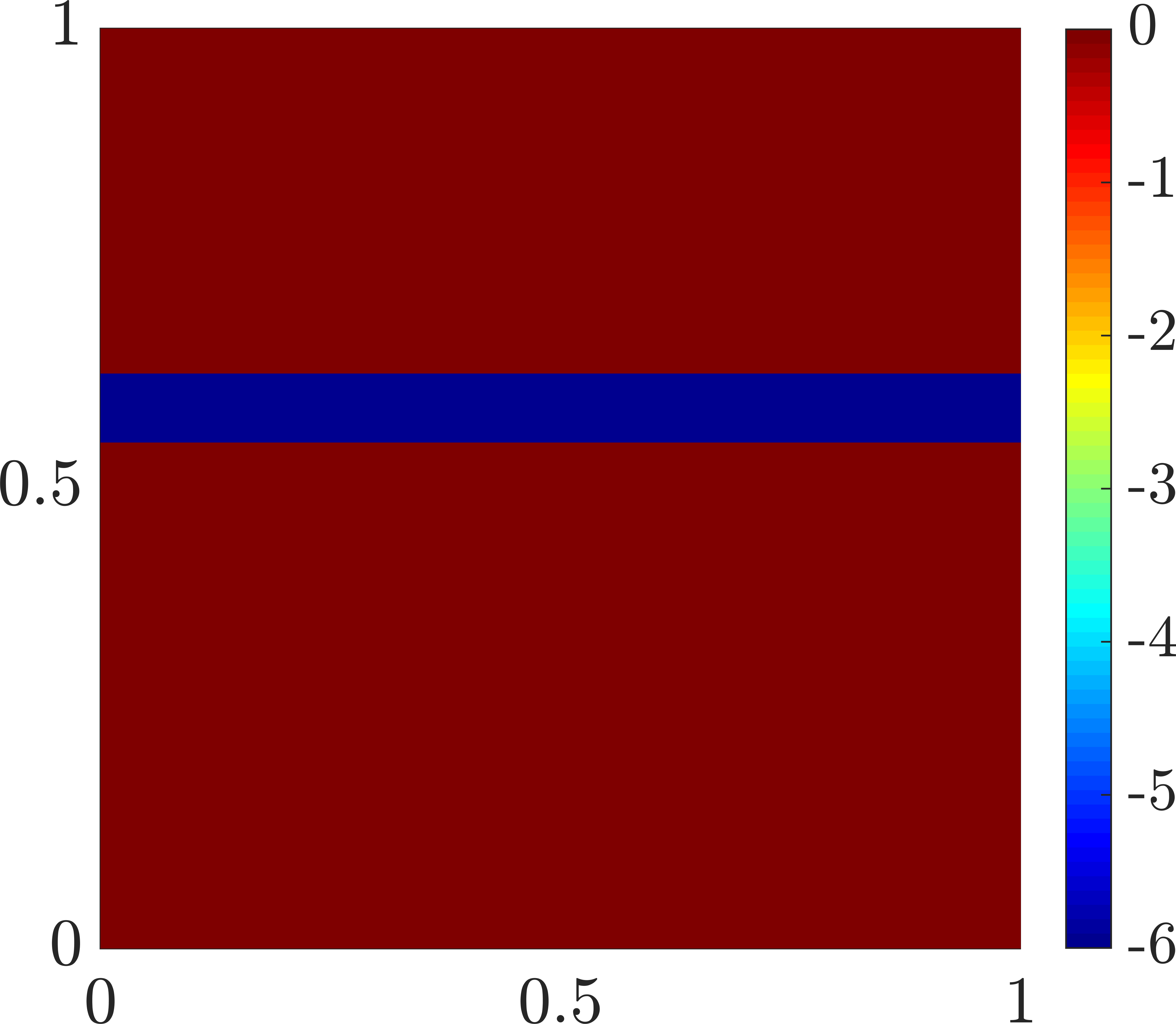}} \label{fig:barrier1_a}
\subfigure[ref2][Pressure]{\includegraphics[scale=0.173]{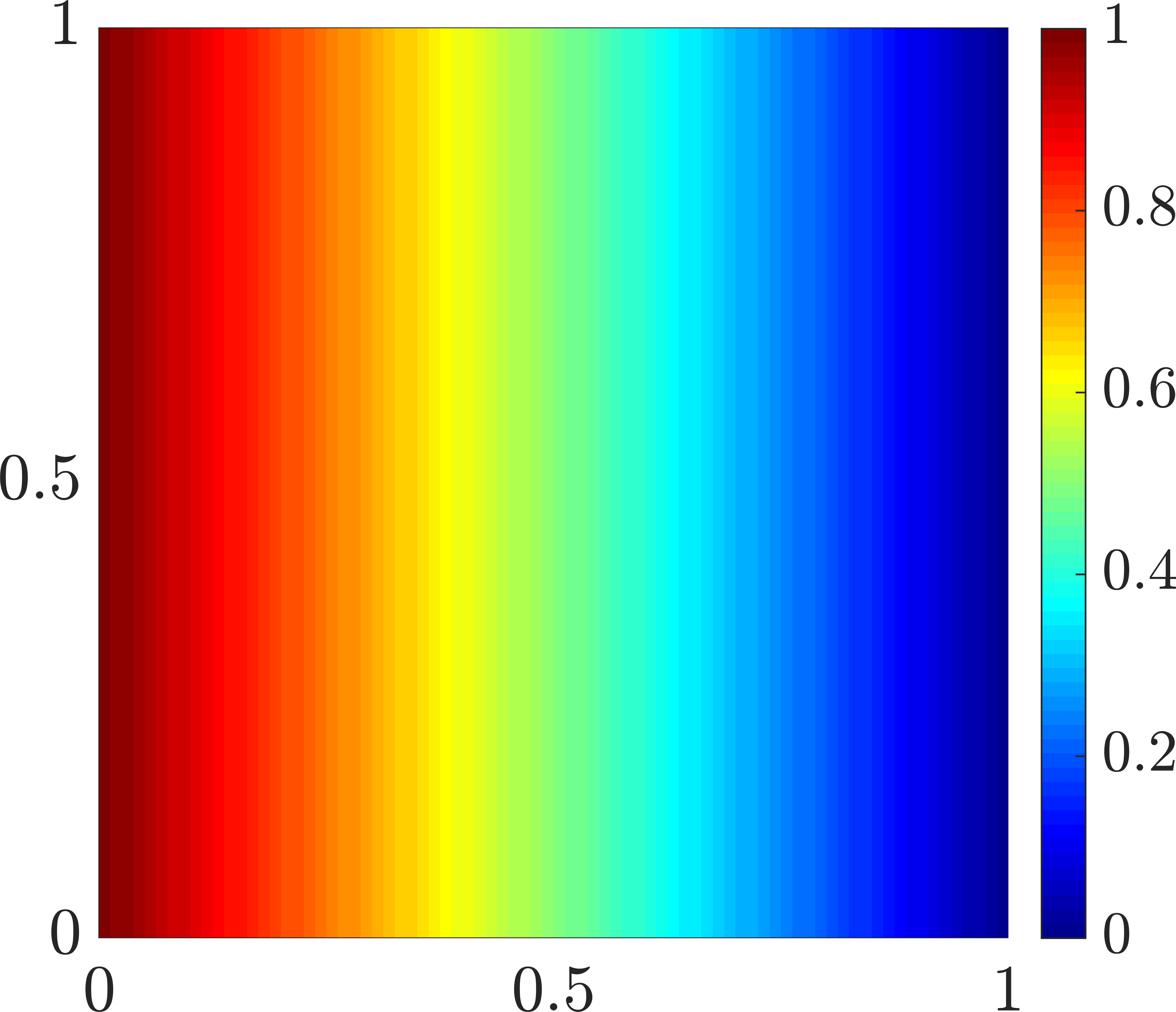}} \label{fig:barrier1_b}
\subfigure[ref2][Flux]{\includegraphics[scale=0.173]{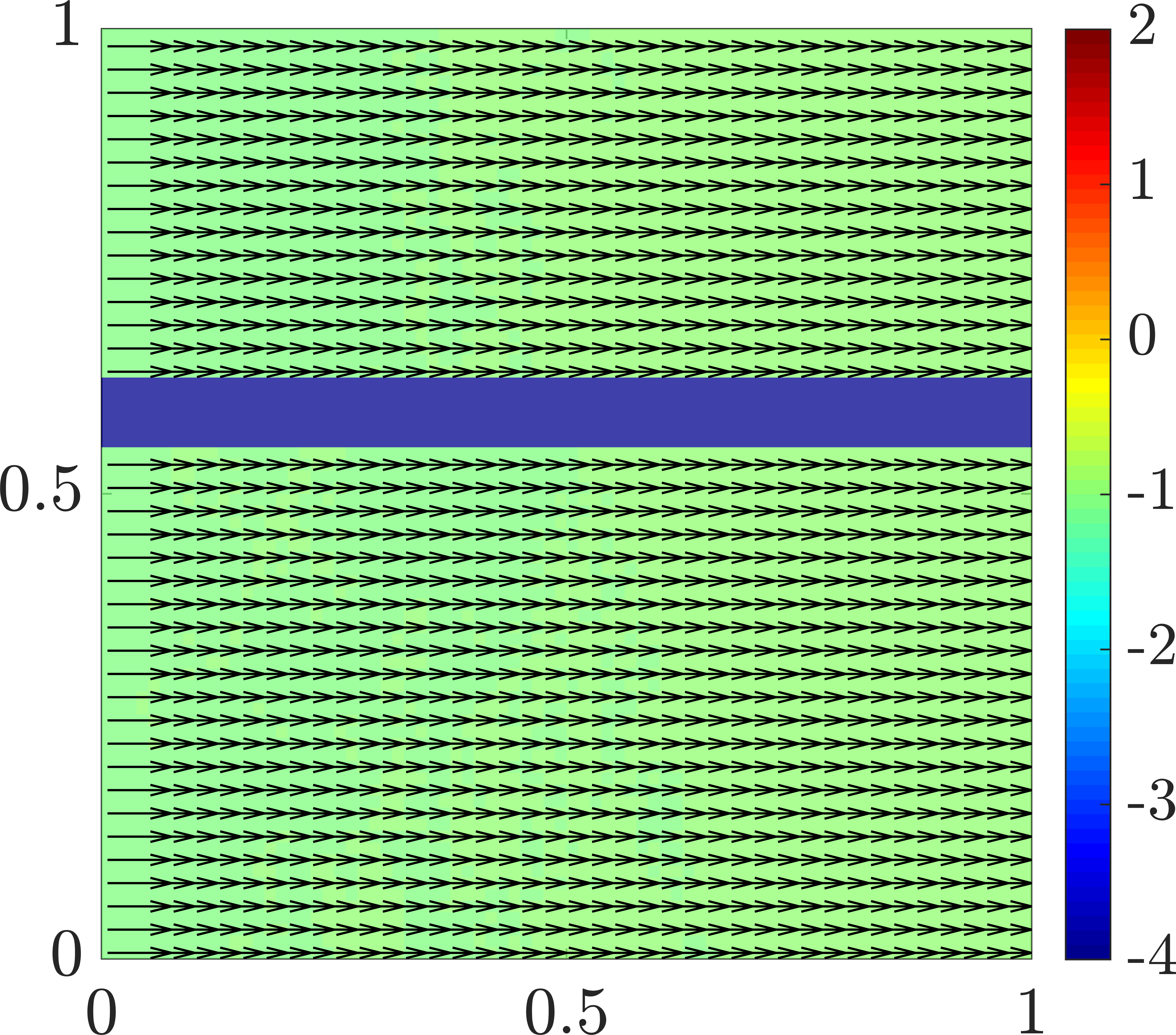}} \label{fig:barrier1_c}
\subfigure[ref2][Flux profile]{\includegraphics[scale=0.173]{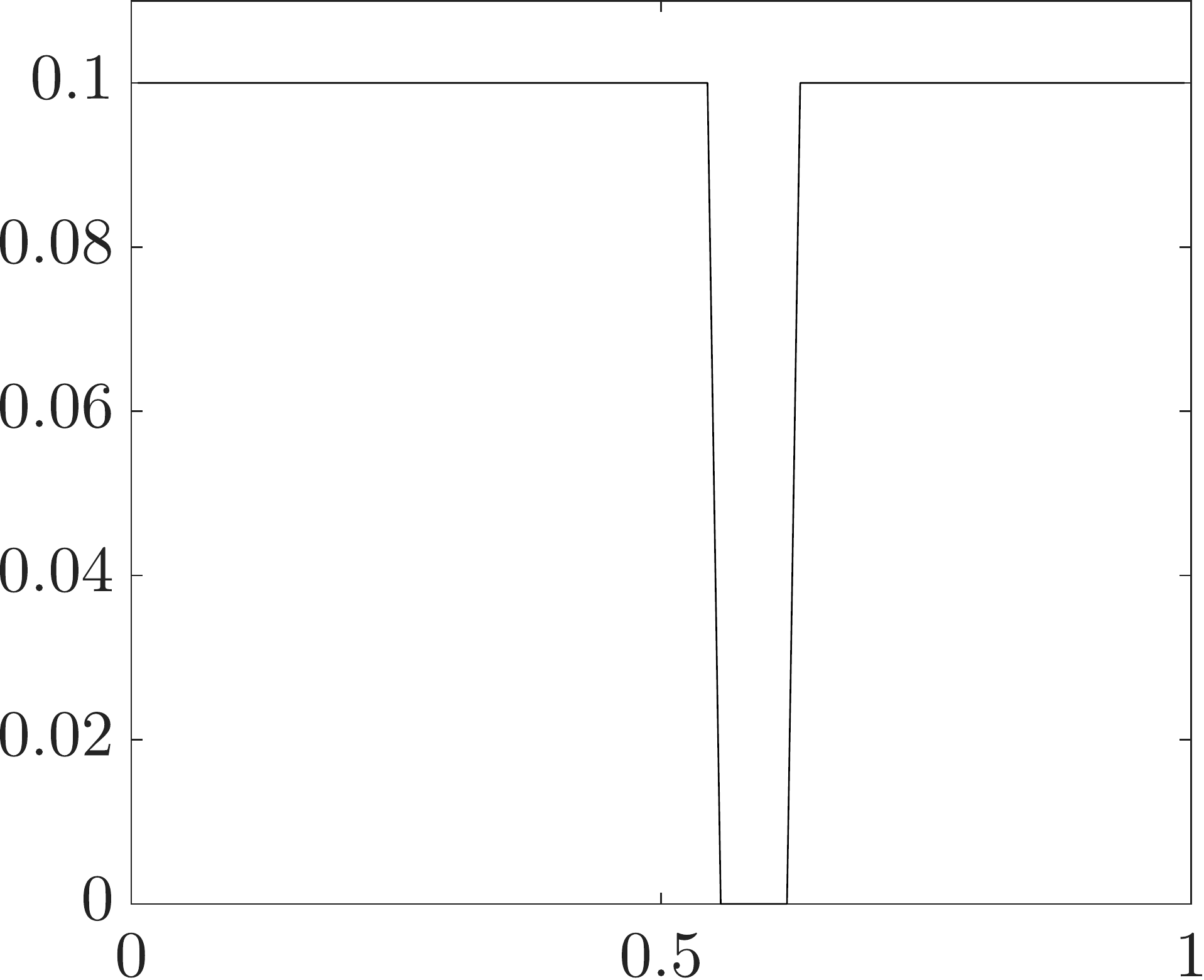} } \label{fig:barrier1_d}
\caption{Horizontal barrier problem. (a) Permeability field (log-scaled) containing a horizontal barrier. (b) Pressure field. (c) Flux. The colors in the flux plot refer to the log-scale flux magnitude.
Additionally, the $x$-component solution of the flux along a vertical line is illustrated in (d), where we note that the flux is discontinuous at the transitions to barrier regions.}
 \label{fig:barrier1}
\end{figure}

Now we focus on fields containing barriers. Once again we define a simplified problem to motivate the interface spaces.
In Figure \ref{fig:barrier1}(a) we consider a high-contrast permeability field containing a horizontal barrier. We show the fine grid solutions for pressure Figure \ref{fig:barrier1}(b) and for flux Figure \ref{fig:barrier1}(c). Here the flow is established by imposing a pressure gradient from left to right and no-flow at top and bottom. The $x$-component of the flux along a vertical line is illustrated in Figure \ref{fig:barrier1}(d), showing the discontinuities at the locations of transitions to barrier regions. Any domain decomposition with more than one subdomain in $x$-direction contains vertical interfaces through which the barrier passes. Let $\Gamma^{\text{barrier}}$ be the set of all the interfaces that contain at least one fine cell in which the absolute permeability is lower than a cutoff value $\zeta_{\min}$. We introduce the new flux spaces $\mathcal{U}_{H}^*$ according to the behavior of the flux solution by replacing the flux linear spaces $\mathcal{U}_{H,1}=\text{span}\{\phi_1, \phi_2\}$ at the interfaces $\Gamma_{i,j}\subset\Gamma^{\text{barrier}}$.

\begin{figure}[htbp]
    \centering
    \includegraphics[scale=0.35]{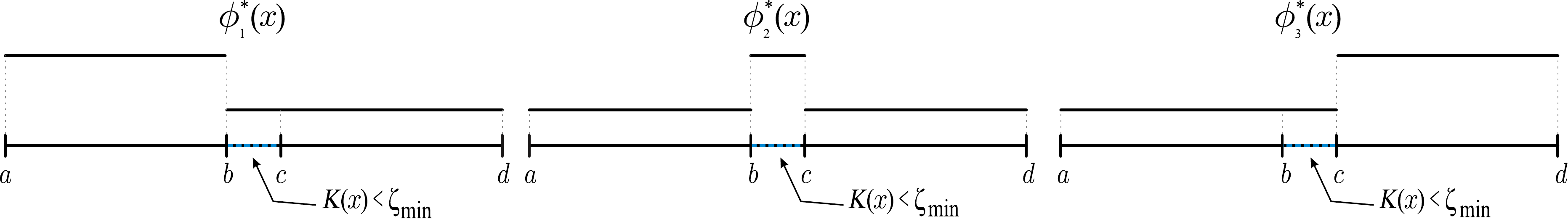} 
    \caption{Physics-based basis functions for flux at the interfaces that contain barriers. }
    \label{fig:new_flux_space}
\end{figure}

Let $\Gamma_{i,j}\subset\Gamma^{\text{barrier}}$ be an interface with support in $[a,d]$ through which a barrier passes in $[b,c]\subset[a,d]$, as sketched in Figure \ref{fig:new_flux_space}. We define the following basis functions:
\begin{equation}\label{new_basis4}
\phi_1^*(x)=\left\{
\begin{array}{cl}
1  &\mbox{if}\ x\in(a,b) \\
0\   &\mbox{otherwise}
\end{array} \right.
\end{equation}
\begin{equation}\label{new_basis5}
\phi_2^*(x)=\left\{
\begin{array}{cl}
1  &\mbox{if}\ x\in(b,c) \\
0\   &\mbox{otherwise} 
\end{array} \right.
\end{equation}
\begin{equation}\label{new_basis6}
\phi_3^*(x)=\left\{
\begin{array}{cl}
1  &\mbox{if}\ x\in(c,d) \\
0\   &\mbox{otherwise}
\end{array} \right.
\end{equation}
The new interface space is defined as $\mathcal{U}_{H}^*=\text{span}\{\phi_1^*, \phi_2^*, \phi_3^*\}$ at $\Gamma_{i,j}^{\text{barrier}}$. These basis functions are not restricted to barriers orthogonal to the interface, similar to the pressure basis. If the interface contains more than one barrier we need to define a new basis function with behavior similar to that of the $\phi_2^*$ for each barrier, plus a constant function for each region between two barriers. The total number of basis functions per interface is thus $1+2N_{\text{barrier}}$, where $N_{\text{barrier}}$ is the number of barriers.

\begin{figure}[htbp]

\center
\subfigure[ref1][Permeability]{\includegraphics[scale=0.166]{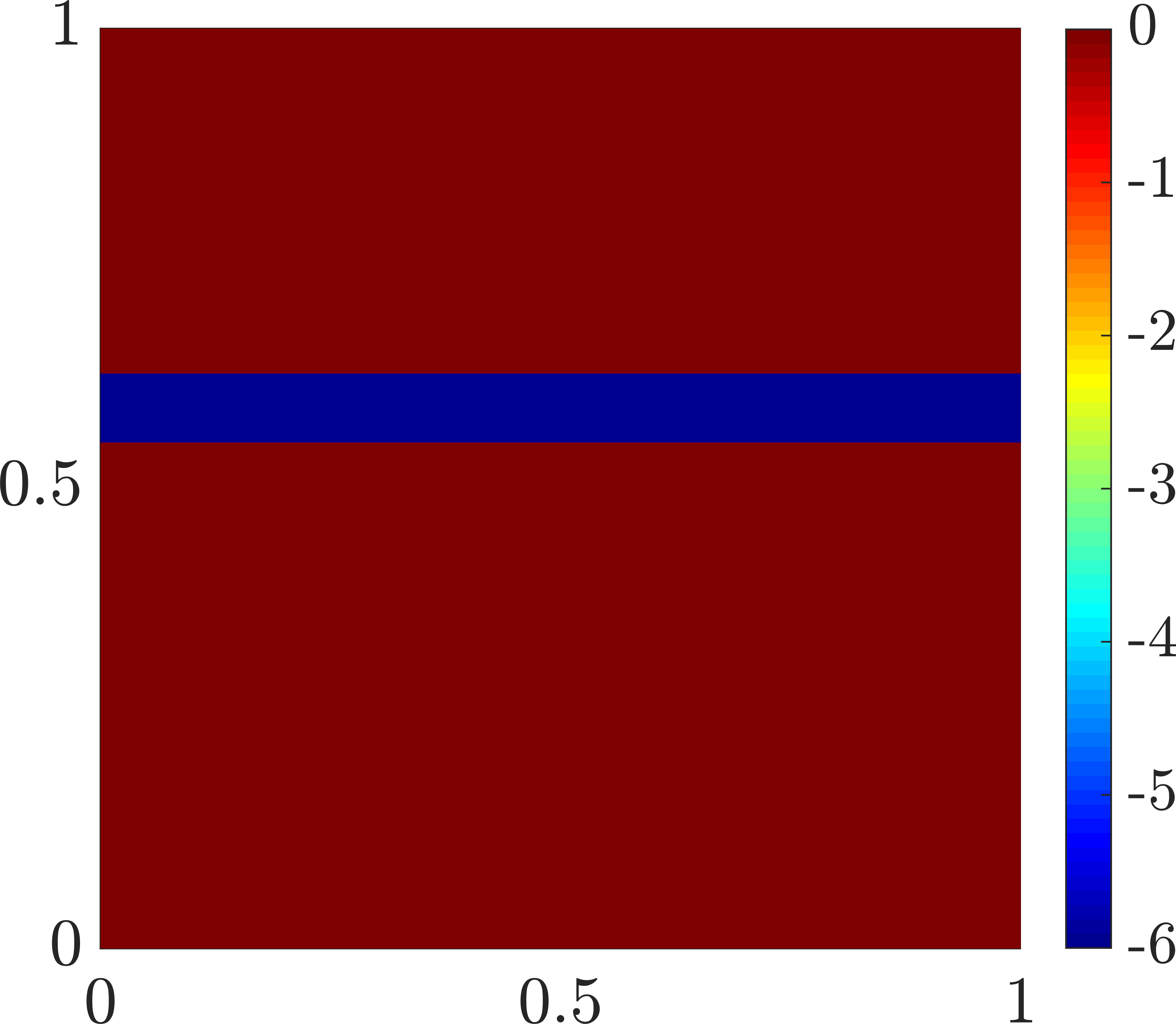} } \label{fig:barrier2_a}
\subfigure[ref2][Flux]{\includegraphics[scale=0.166]{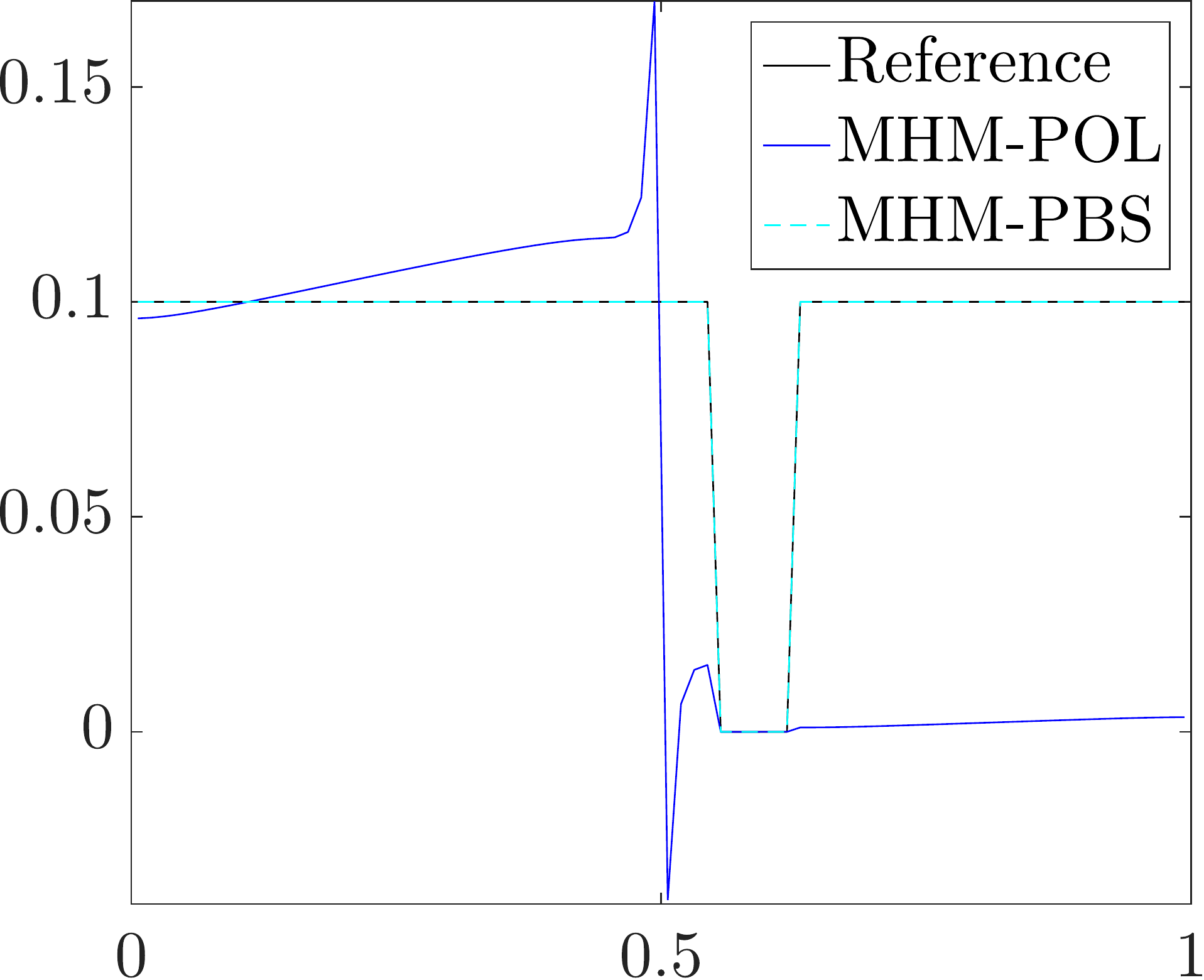} } \label{fig:barrier2_b}
\subfigure[ref2][Permeability]{\includegraphics[scale=0.166]{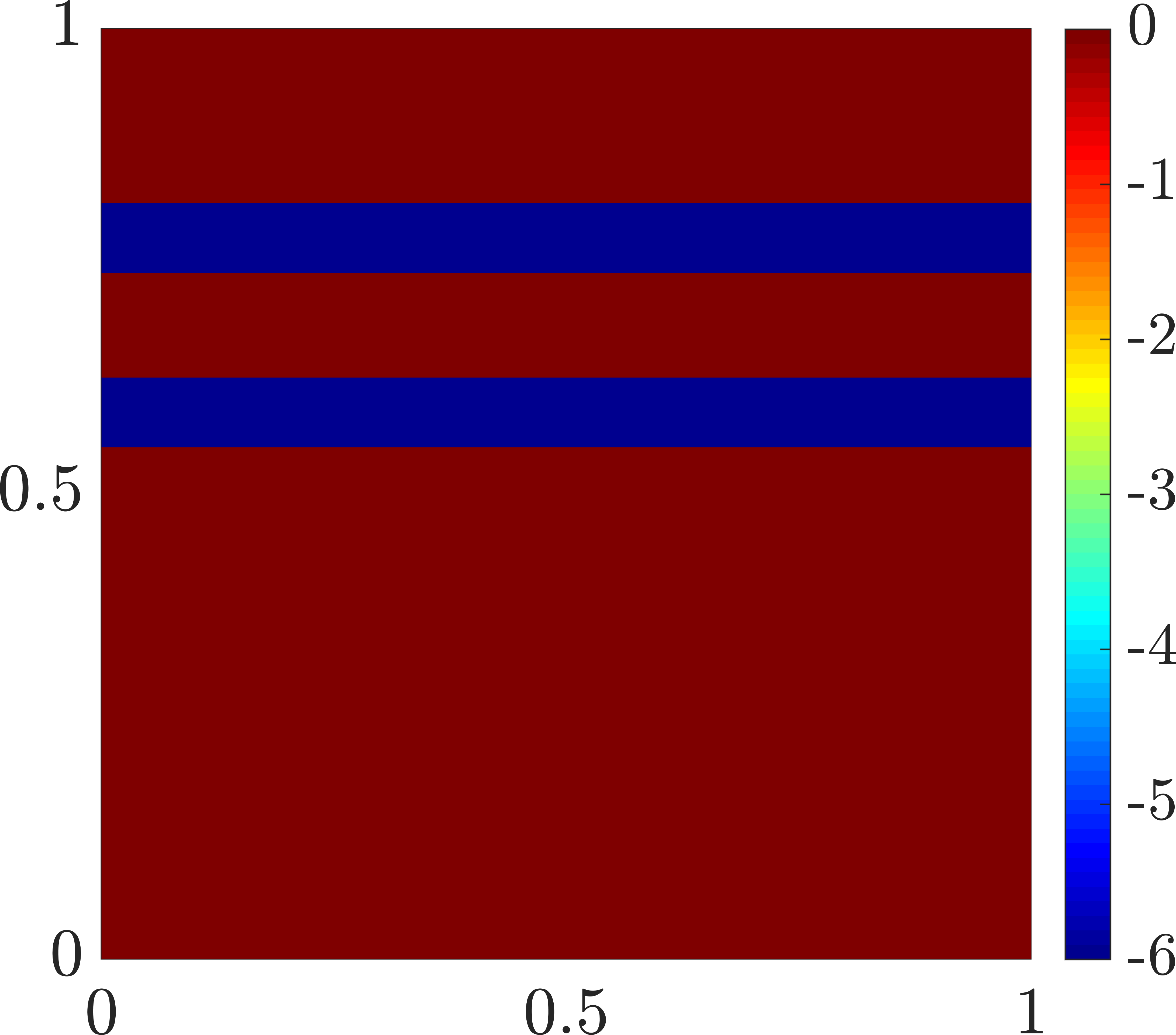} } \label{fig:barrier2_c}
\subfigure[ref2][Flux]{\includegraphics[scale=0.166]{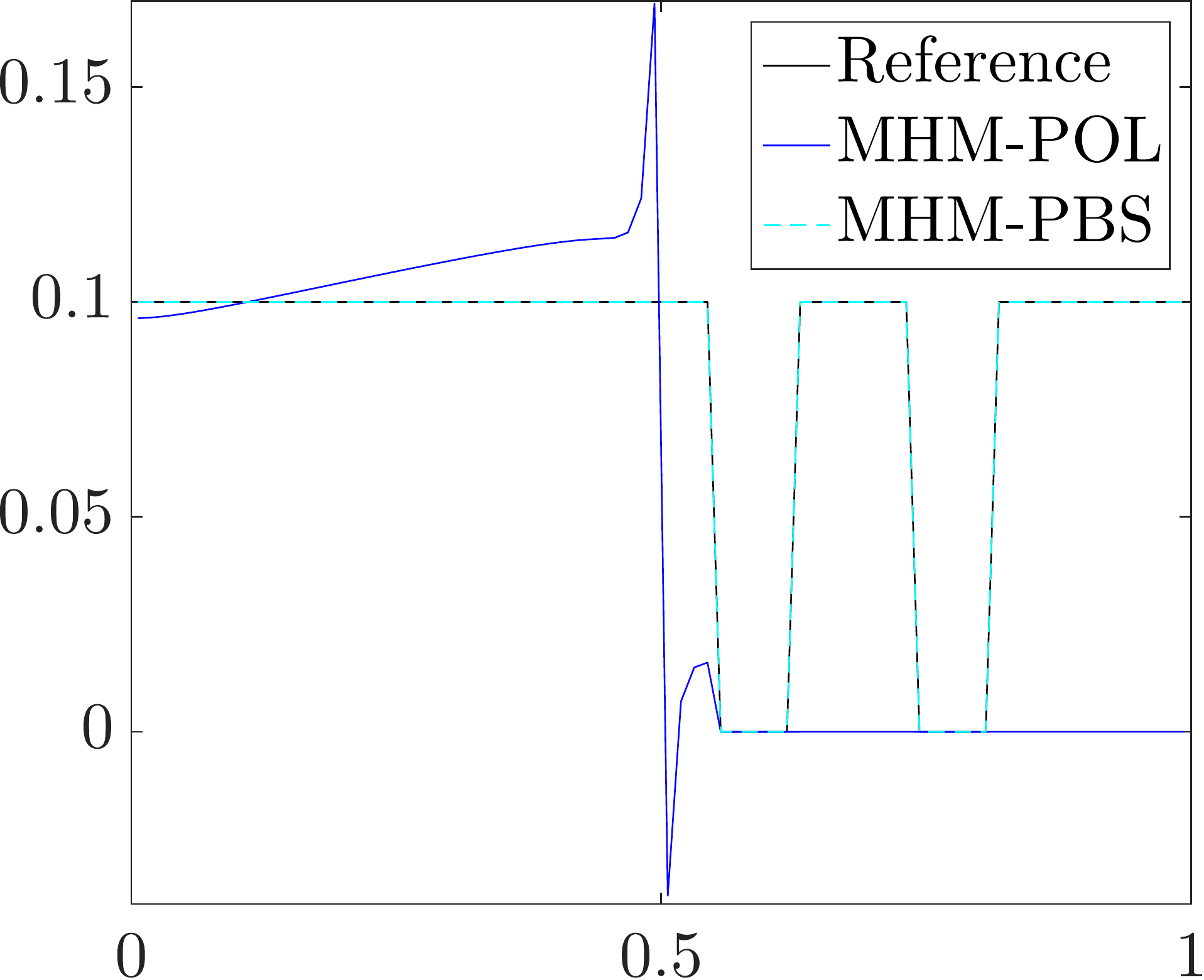} } \label{fig:barrier2_d}
\caption{High-contrast permeability field (log-scaled) with one (a) and two barriers (c). Fine-grid reference, MHM-POL and MHM-PBS solutions for the $x$-component of the flux considering one (c) and two barriers (d). Notice that for both one and two barriers the correct $x$-flux solution is only captured by the MHM-PBS.}
 \label{fig:barrier2}
\end{figure}

In Figure \ref{fig:barrier2} we show the $x$-component of the flux along $x=0.5$ provided by the MHM (by setting the $\alpha$ parameter of the MRCM to the value $10^{6}$) in a domain decomposition with $2\times2$ coarse cells. We compare the approximations provided by the MHM-POL (MHM combined with the linear spaces) and MHM-PBS (MHM combined with the physics-based spaces) with the fine-grid solution.
We note that for both single (Figure \ref{fig:barrier2}(a) and Figure \ref{fig:barrier2}(c)) and multiple (Figure \ref{fig:barrier2}(b) and Figure \ref{fig:barrier2}(d)) barriers the correct $x$-component of the flux is only captured by MHM-PBS.

\subsection{Experiments with the physics-based interface spaces} 
The examples in Figure \ref{fig:frac2} and Figure \ref{fig:barrier2} illustrate that the usual linear interface spaces fail to approximate the solution in the presence of fractures and barriers even for simple problems. The proposed physics-based spaces have been able to represent the behavior of the solution in those problems. 
In this subsection, we show that the physics-based spaces work well also in slightly more complex permeability fields.  
Initially, we consider a permeability field with fractures and use the physics-based pressure space combined with the MMMFEM. Then we consider a permeability field with barriers and use the physics-based flux space combined with the MHM.

 \subsubsection{The MMMFEM for permeability fields with fractures}  

\begin{figure}[htbp]
\center
\subfigure[ref1][Permeability]{\includegraphics[scale=0.37]{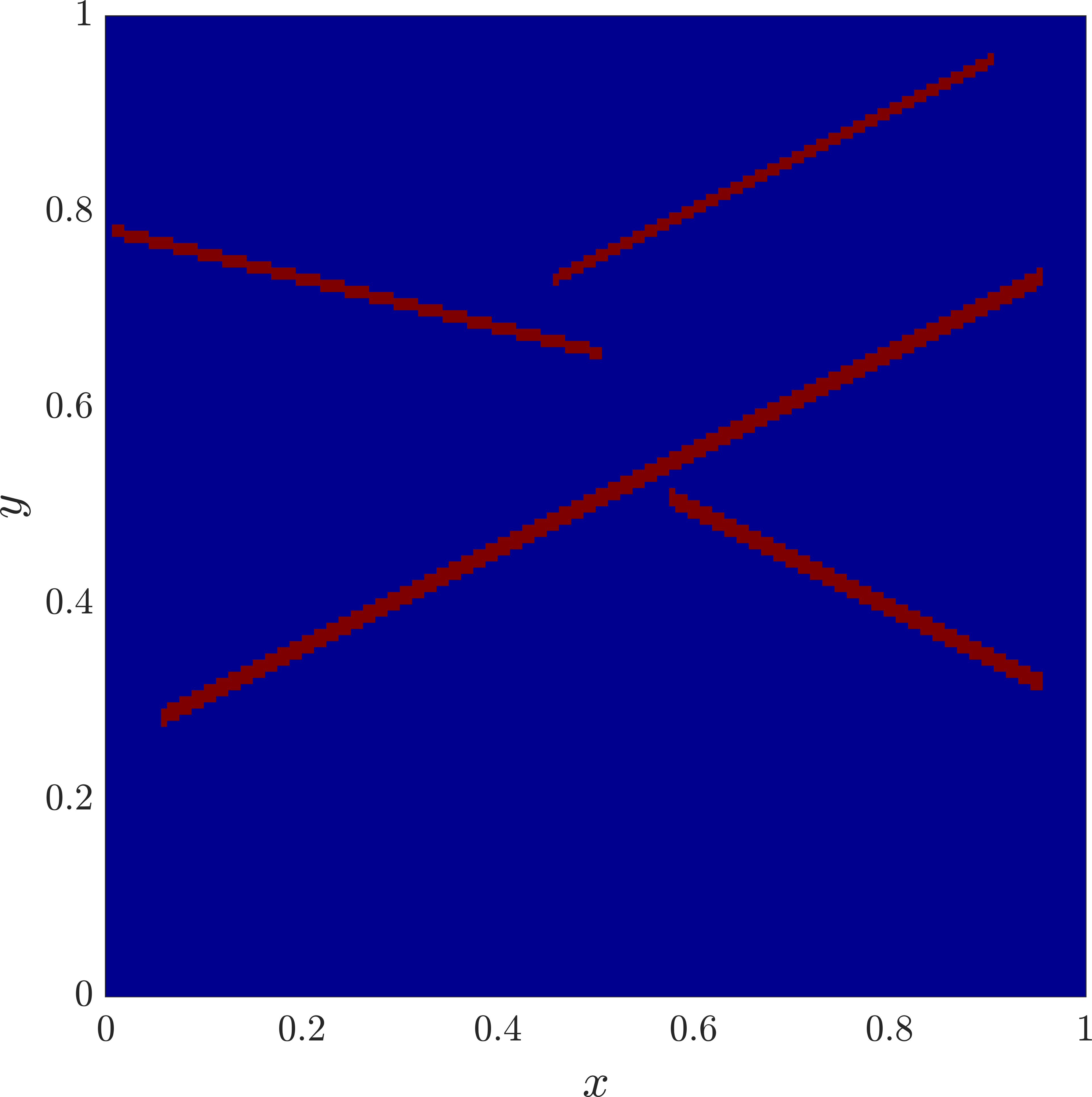} }   \label{fig:perm_fracs}
\subfigure[ref2][Pressure errors]{\includegraphics[scale=0.37]{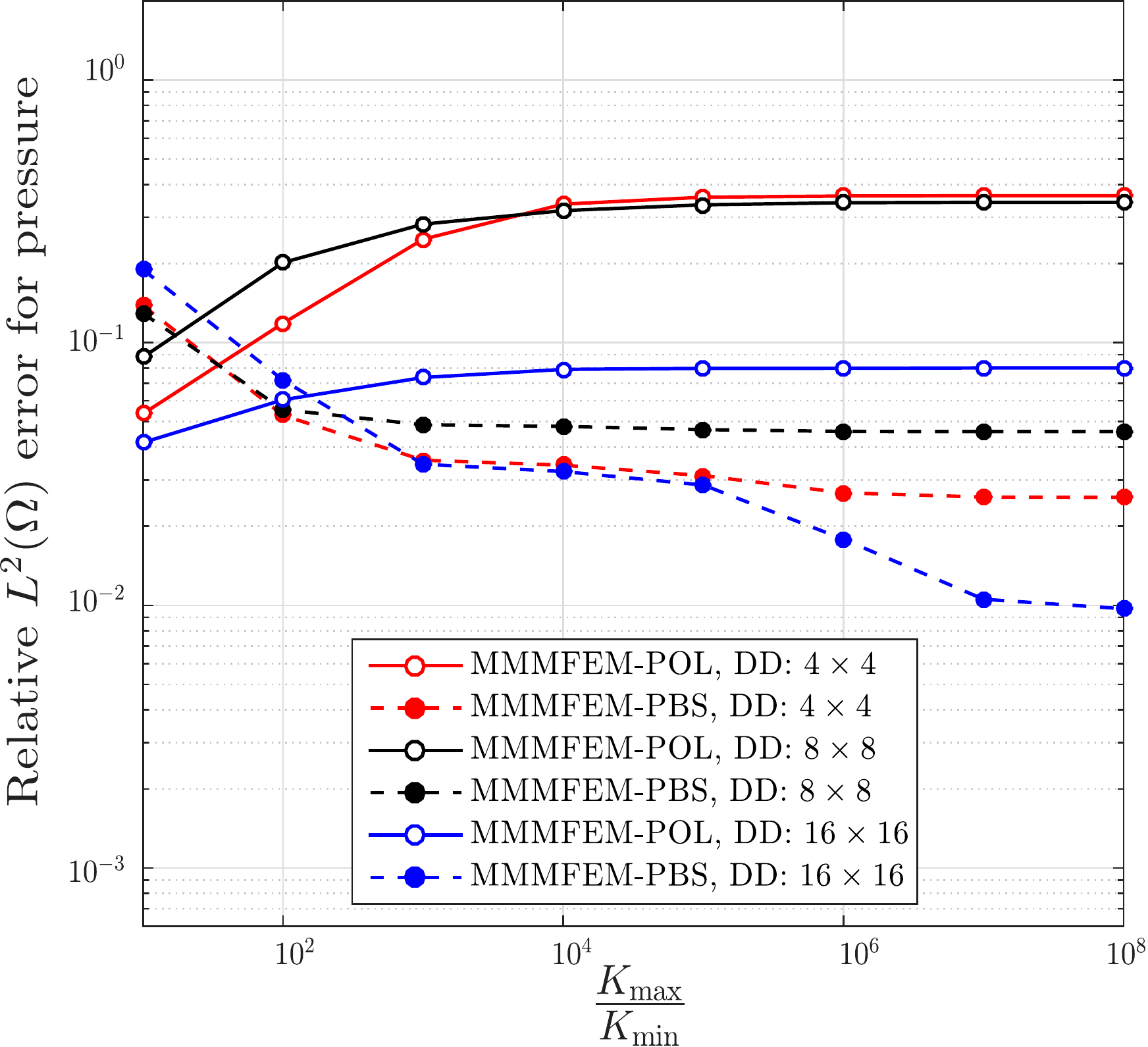} }   \label{fig:error_fracs}
\caption{High-contrast permeability field with fractures (a). We consider the permeability in the fracture, ${K_{\max}}$ varying from $10$ to $10^8$ whereas the background is homogeneous with $K=1$.
Relative $L^2(\Omega)$ pressure errors as function of the contrast are shown for the MMMFEM-POL and MMMFEM-PBS (b). Three domain decompositions are considered: with $4\times4$, $8\times8$ and $16\times16$ coarse cells. We note a significant improvement for the MMMFEM-PBS in all the meshes, especially for high-contrast.}
 \label{fig:press_for_fracs}
\end{figure}

\begin{figure}[htbp]
    \centering
        \includegraphics[scale=0.25]{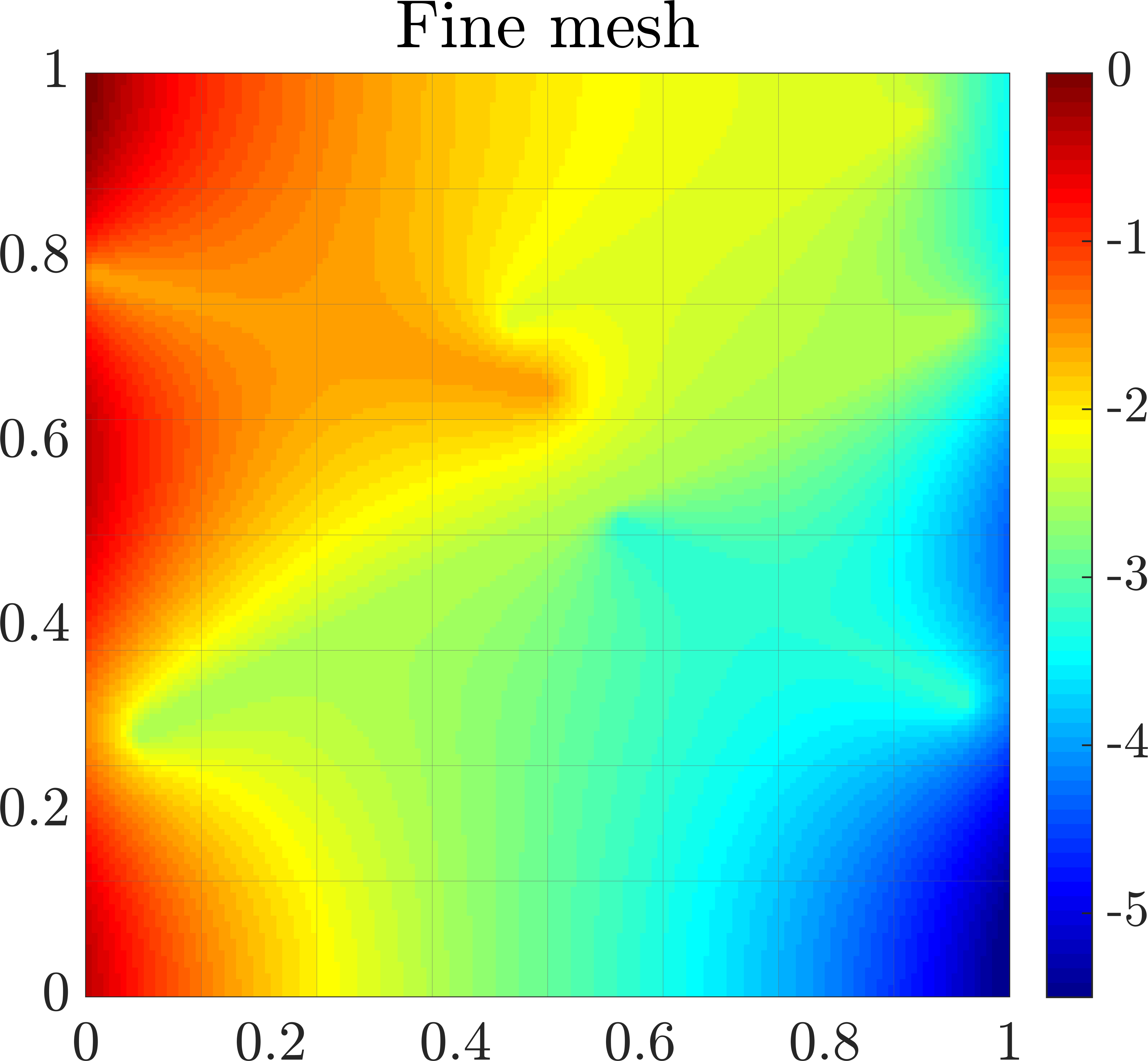}
    \includegraphics[scale=0.25]{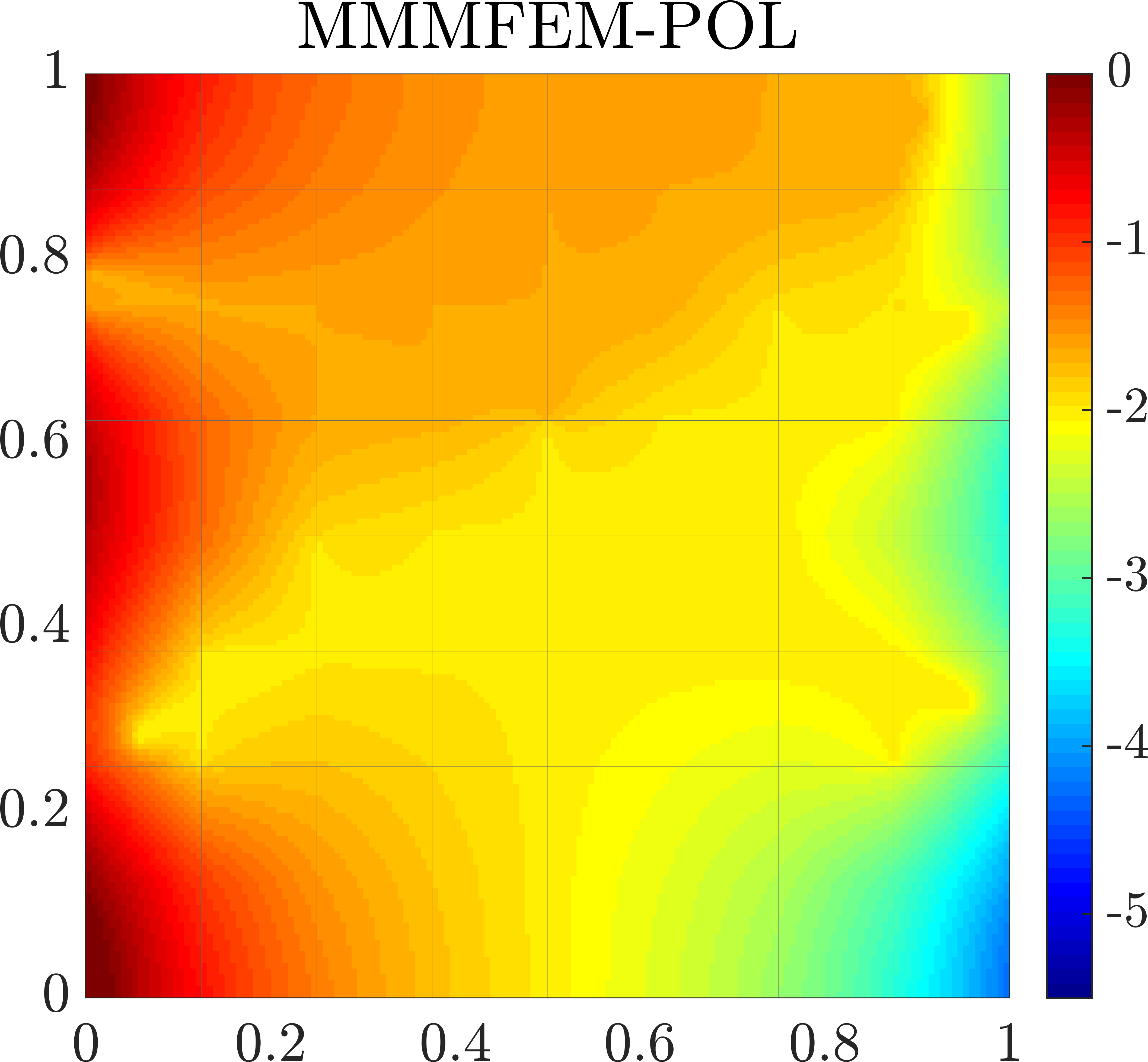}
    \includegraphics[scale=0.25]{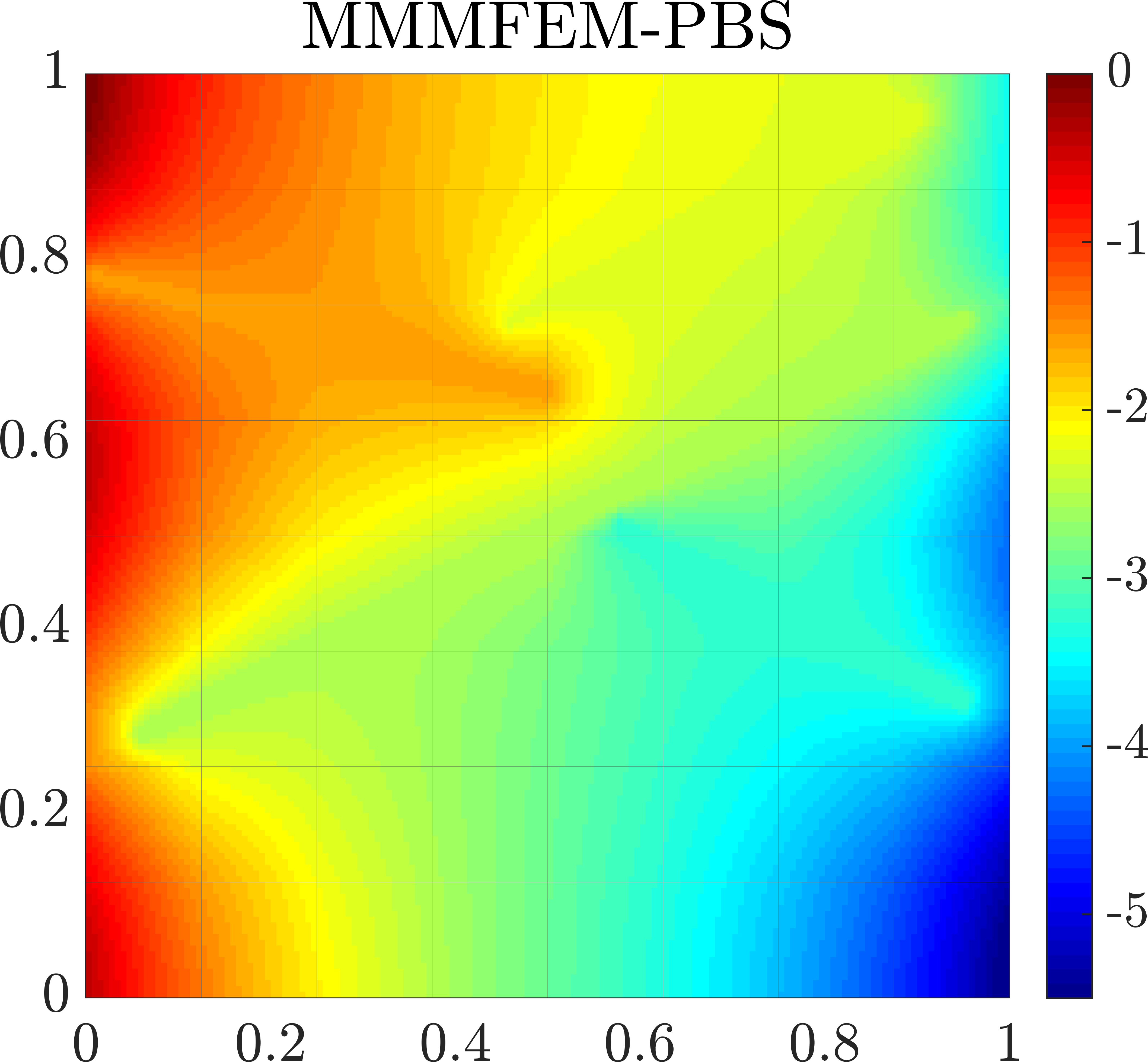}
      \caption{Pressure approximations considering the contrast of ${K_{\max}}/{K_{\min}}=10^8$. Left to right: fine mesh, MMMFEM-POL, and MMMFEM-PBS solutions. The domain decomposition considered contains $8\times8$ coarse cells and is illustrated by the lines in the plot. We note that the MMMFEM-PBS solution is more accurate than the MMMFEM-POL.}
       \label{fig:press_fracs}
\end{figure}

The first study considers a permeability field containing fractures. In this case, the physics-based pressure space is applied to the interfaces that contain cells in which the absolute permeability is larger than a cutoff value $\zeta_{\max}$. We combine the new interface space with the MMMFEM.
We consider the permeability in the fracture, ${K_{\max}}$ varying from $10$ to $10^8$ whereas the background is homogeneous with $K=1$, see Figure \ref{fig:press_for_fracs} (a). The cutoff value set to capture the fractures is $\zeta_{\max}=1$ in all cases. The domain considered is $\Omega=[0,1]\times[0,1]$ containing $160\times160$ cells. In Figure \ref{fig:press_for_fracs} (b) we show the relative $L^2(\Omega)$ errors for pressure as function of the permeability contrast. Three domain decompositions are considered: $4\times4$, $8\times8$ and $16\times16$ coarse cells, each one containing, respectively $40\times40$, $20\times20$ and $10\times10$ fine cells. The boundary conditions in the simulations of this subsection are no-flow at the top and bottom boundaries along with an imposed flux on the left and right boundaries. No source terms are considered. 
We compare the multiscale solution considering the usual linear polynomial (MMMFEM-POL) and the physics-based (MMMFEM-PBS) pressure interface spaces. We note that for  permeability contrasts larger than 100 the improvement with the physics-based spaces is significant for all the domain decompositions considered. Domain decompositions with more subdomains present smaller errors.
Figure \ref{fig:press_fracs} shows the pressure approximations for the decomposition of $8\times8$ coarse cells in the highest permeability contrast ${K_{\max}}/{K_{\min}}=10^8$. It is clear that the imprecisions of the MMMFEM solution with the linear interface spaces are corrected with the use of the physics-based ones. 

 \subsubsection{The MHM for permeability fields with barriers}  
 
\begin{figure}[htbp]
\center
\subfigure[ref1][Permeability]{\includegraphics[scale=0.37]{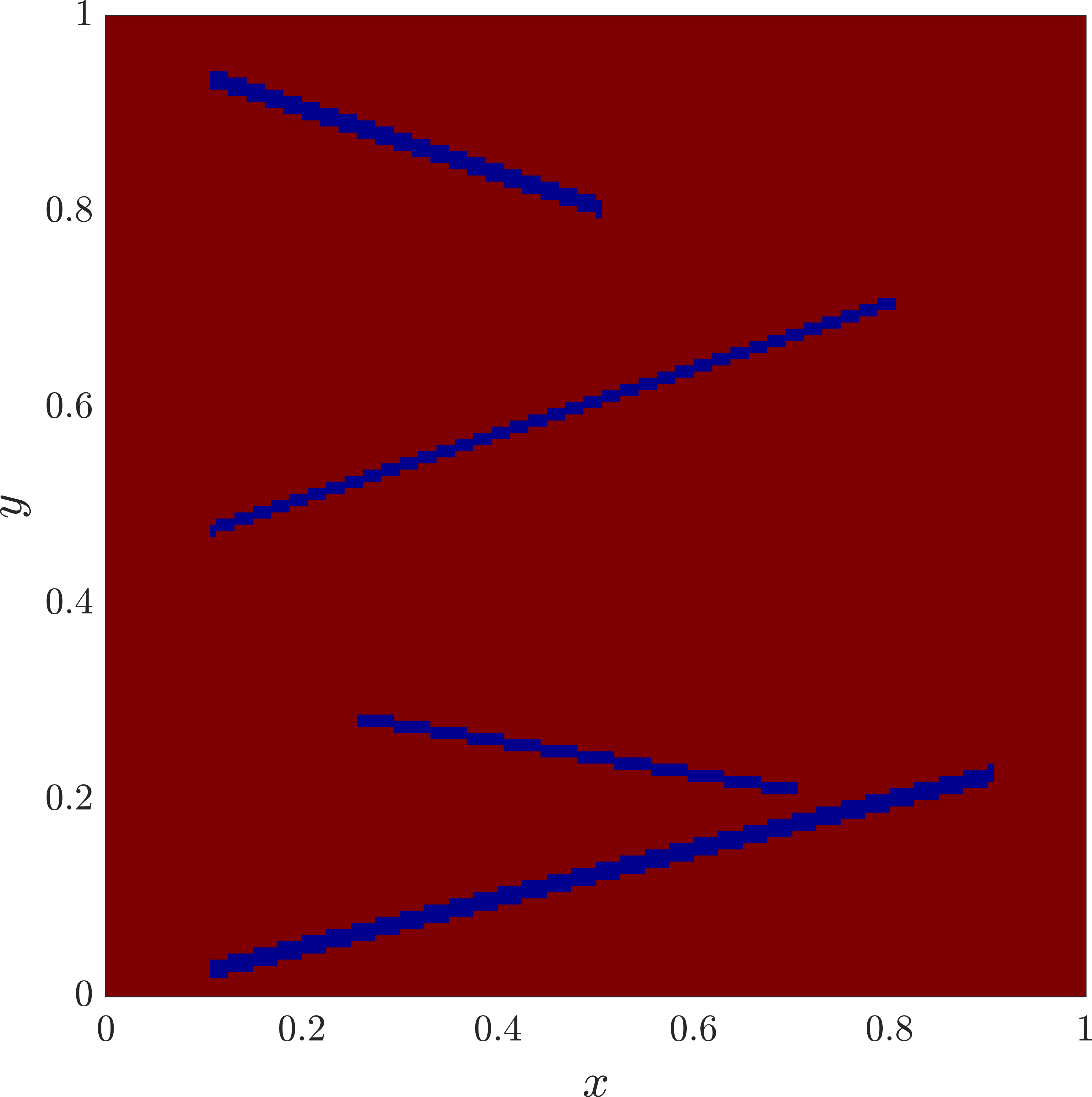} }   \label{fig:perm_barriers}
\subfigure[ref2][Flux errors]{\includegraphics[scale=0.37]{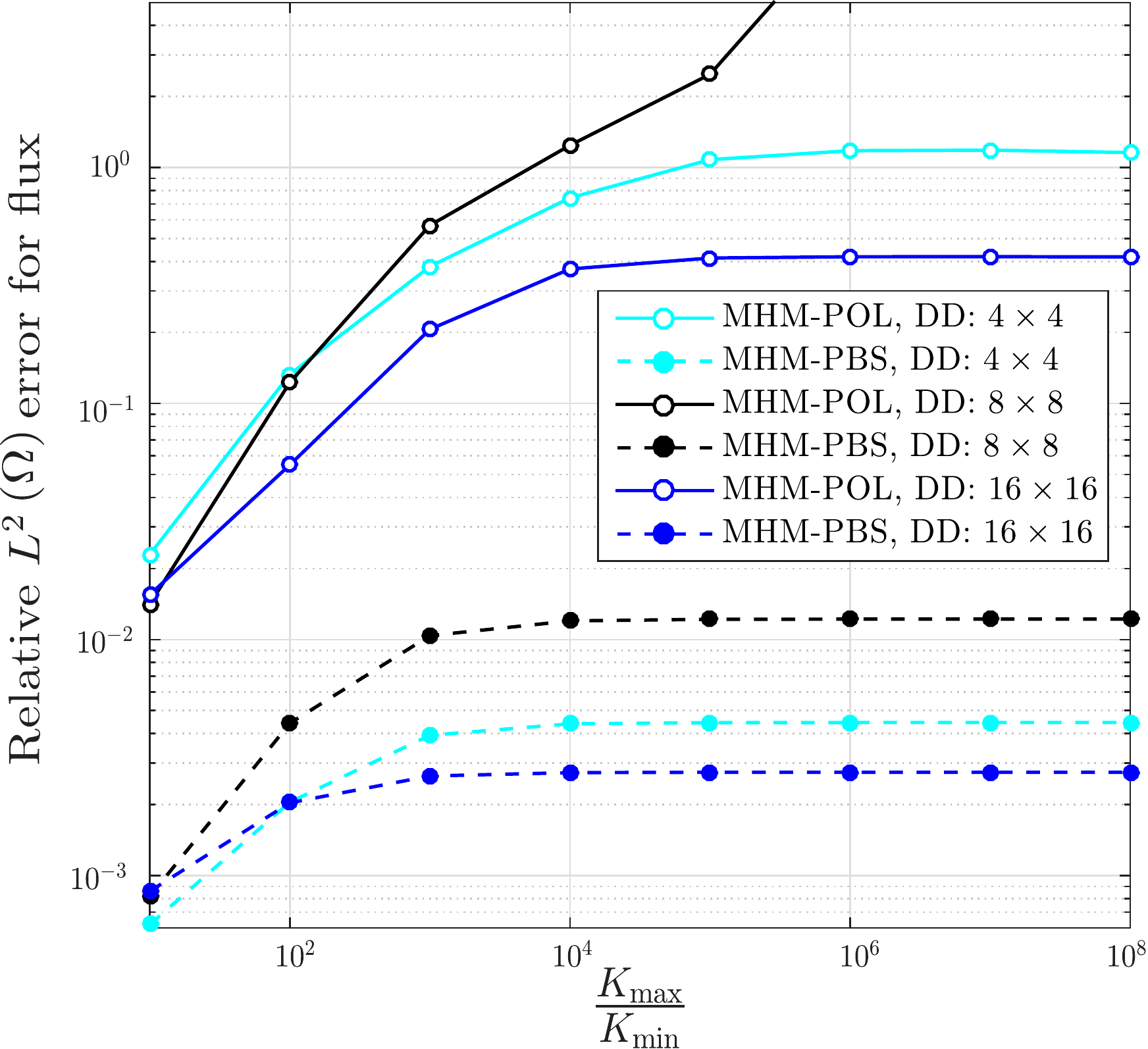} }   \label{fig:error_barriers}
\caption{High-contrast permeability field with barriers (a). We consider the permeability in the barrier, ${K_{\min}}$ varying from $10^{-8}$ to $10^{-1}$ whereas the background is homogeneous with $K=1$.
Relative $L^2(\Omega)$ flux errors as function of the contrast are shown for the MHM-POL and MHM-PBS (b). Three domain decompositions are considered: with $4\times4$, $8\times8$ and $16\times16$ coarse cells. We note that for high-contrast the only accurate approximations are produced by the MHM-PBS (in all the meshes).}
 \label{fig:flux_for_barriers}
\end{figure}

\begin{figure}[htbp]
    \centering
        \includegraphics[scale=0.25]{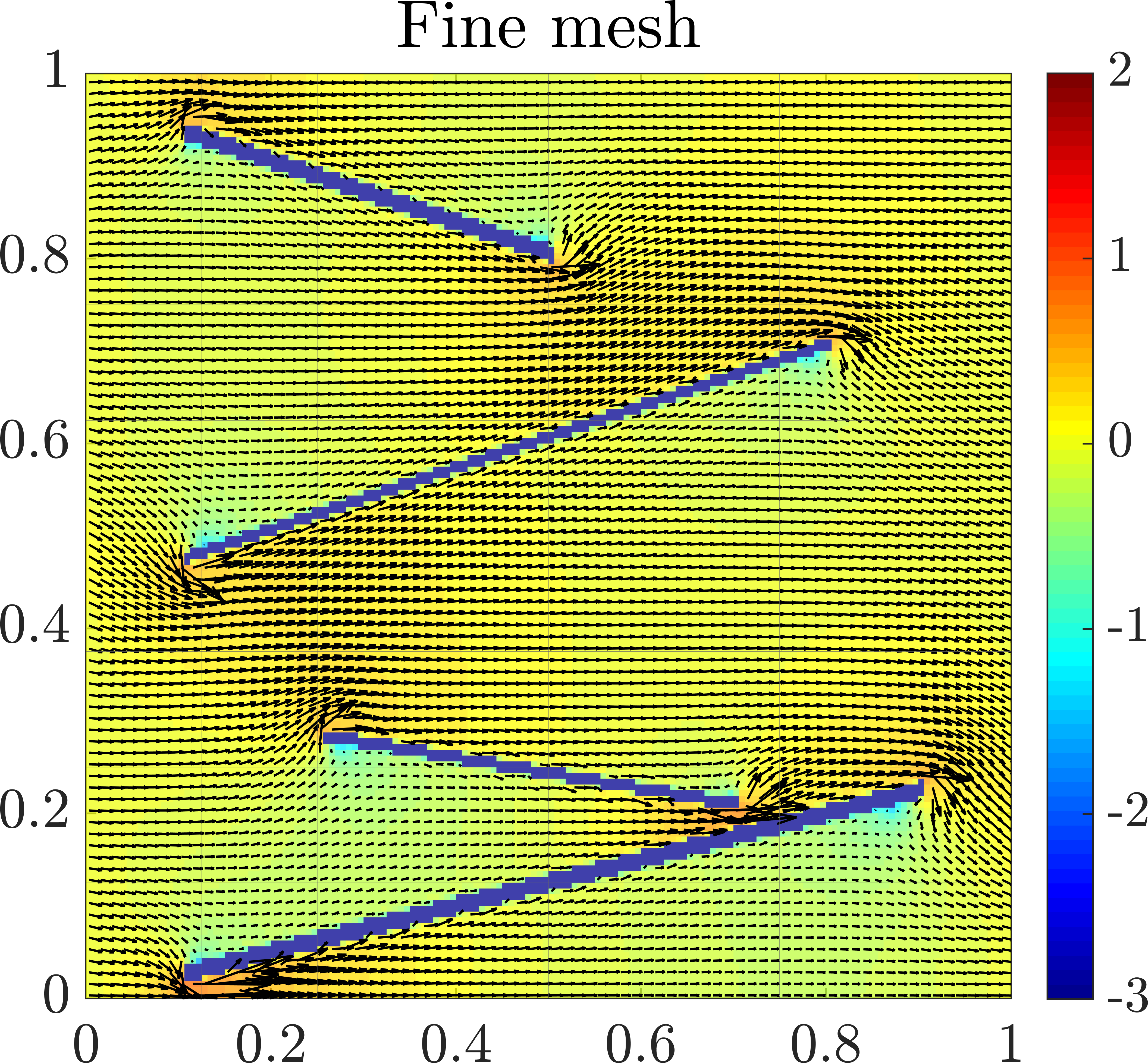}
    \includegraphics[scale=0.25]{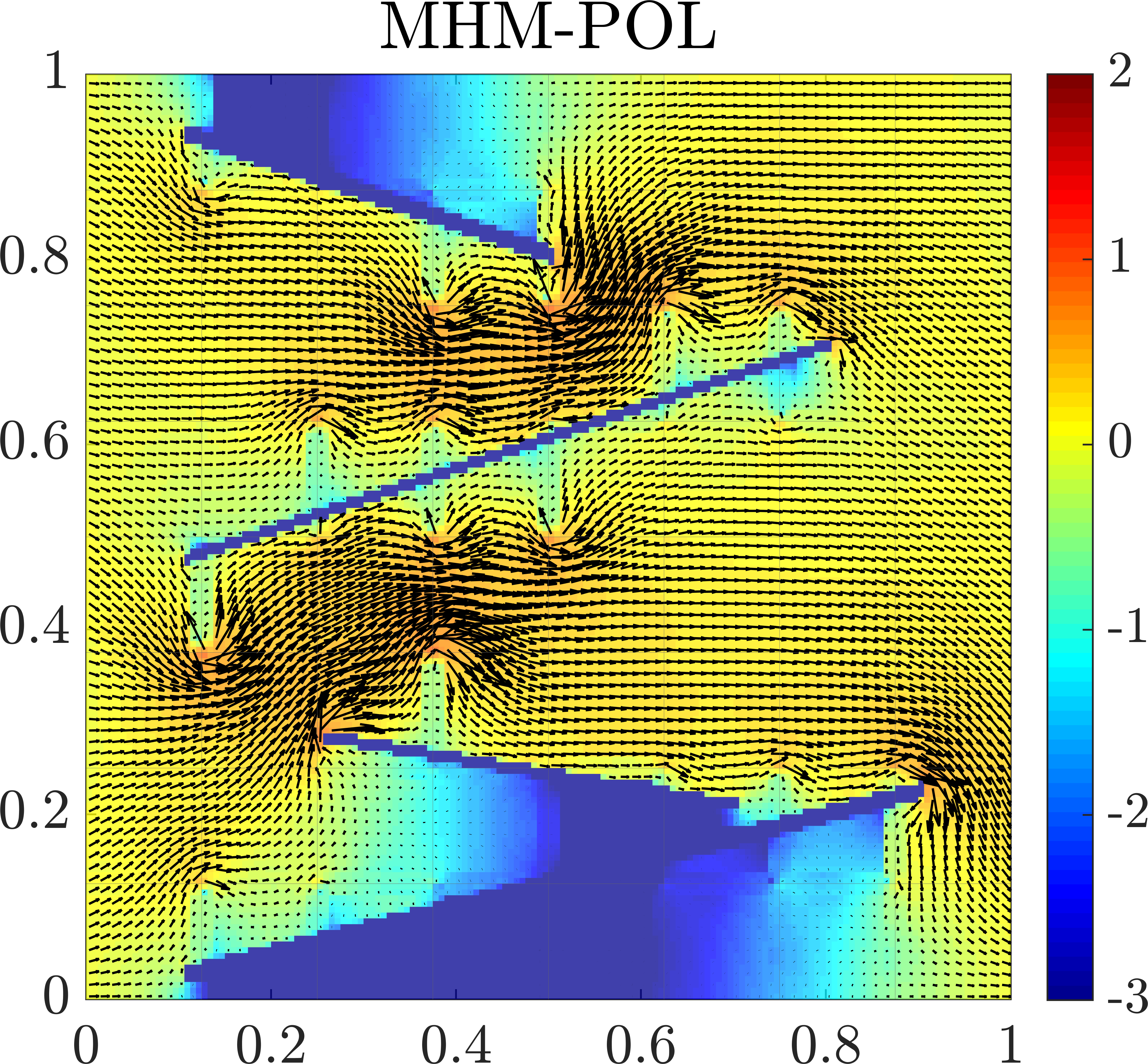}
    \includegraphics[scale=0.25]{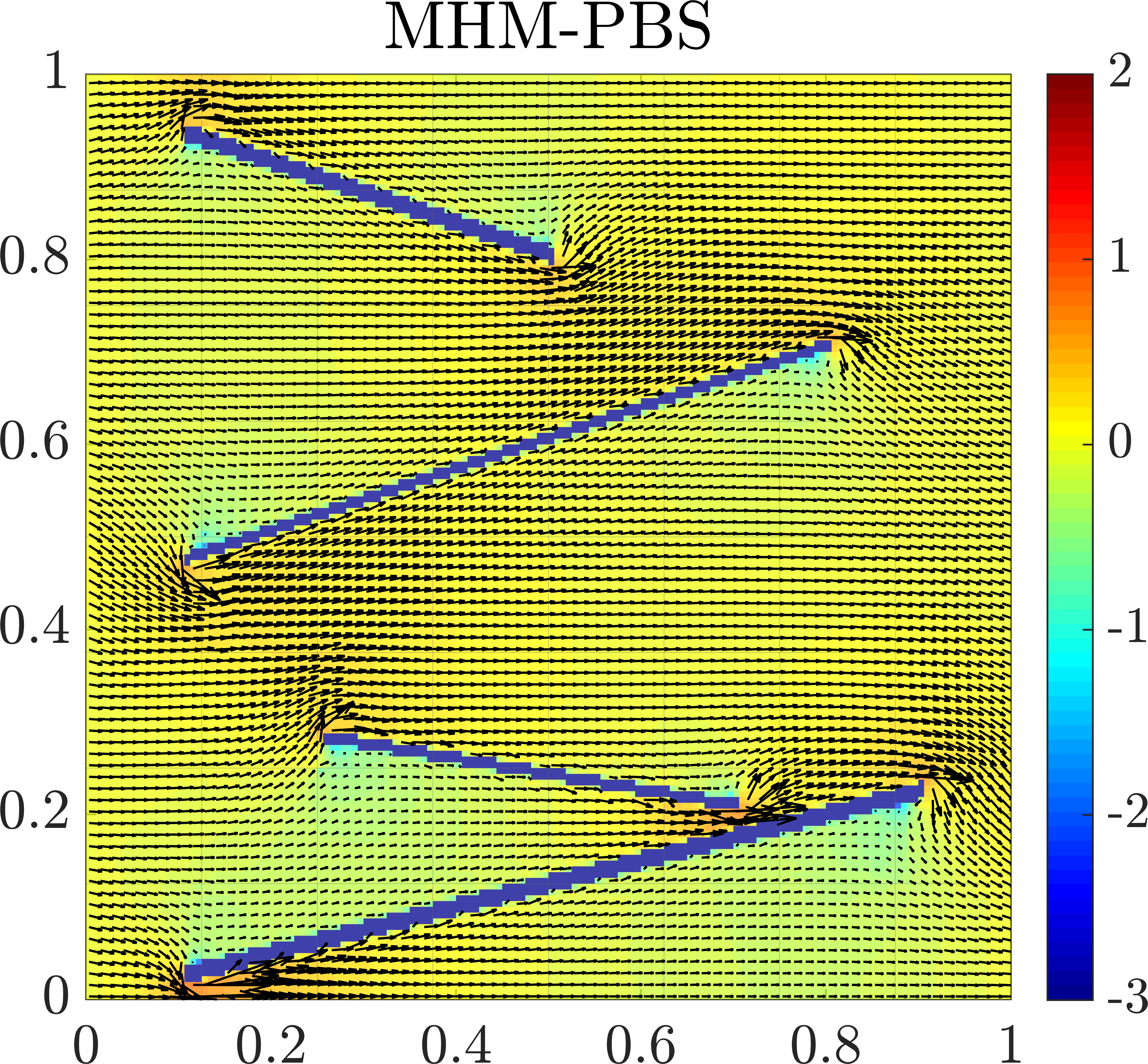}
      \caption{Flux approximations considering the contrast of ${K_{\max}}/{K_{\min}}=10^8$. Left to right: fine mesh, MHM-POL, and MHM-PBS solutions. The colors in the flux plot refer to the log-scale flux magnitude. The domain decomposition considered contains $8\times8$ coarse cells and is illustrated by the lines in the plot. We note that the MHM-PBS approximation is closely related to the reference solution and the MHM-POL is inaccurate.}
      \label{fig:flux_barriers}
\end{figure}

In the next experiment, the same problem of the previous subsection is considered, except that the fractures are replaced by barriers of low permeability. The physics-based flux space is used at interfaces that contain cells in which the absolute permeability is lower than a cutoff value $\zeta_{\min}$. We combine the flux space with the MHM to approximate the solution. Similar to the previous example, we compare different permeability contrasts. 
Here, we consider the permeability in the barrier, ${K_{\min}}$ varying from $10^{-8}$ to $10^{-1}$ whereas the background is homogeneous with $K=1$, see Figure \ref{fig:flux_for_barriers} (a). The cutoff value set to capture the barriers is $\zeta_{\min}=1$ in all cases.
The relative $L^2(\Omega)$ errors for flux as function of the contrast are displayed in Figure \ref{fig:flux_for_barriers} (b). Similar to the previous case, the solution is more accurate by using the physics-based space instead of the linear space. The flux errors provided by the linear spaces increase quickly with the contrast whereas the errors from the physics-based spaces are controlled. The results are consistent for the three domain decompositions, where the smaller errors are provided by the decompositions with more subdomains. 
In Figure \ref{fig:flux_barriers} we show the flux approximations for the decomposition of $8\times8$ coarse cells in the highest permeability contrast ${K_{\max}}/{K_{\min}}=10^8$. We note that the MHM-POL solution is inaccurate and the MHM-PBS approximation captures the correct behavior of the reference solution.

\section{Adaptive MRCM with the physics-based interface spaces}\label{aMRCM-PBS}

We can conclude from the numerical studies reported above that the physics-based pressure (respectively, flux) space is fundamental to produce an accurate pressure (resp., flux) solution in presence of fractures (resp., barriers). Therefore if fractures and barriers appear in a single interface one can achieve a better approximation of pressure and flux using both pressure and flux physics-based interface spaces. To attain this goal one feature necessary for the multiscale method used is the ability to include the interface spaces independently. For this purpose, we combine the physics-based interface spaces with the MRCM, whose formulation enables to include the spaces separately. In the MRCM framework, we consider the adaptive version (denoting by $a$MRCM) that consists of setting the function $\alpha(\mathbf {x})$ depending on the variation of the permeability field at the interfaces. 
In a previous work \cite{bifasico}, we found that the best choice is to set a small value $\alpha_{\text{small}}$ (pressure is favored) for regions where the permeability is larger than a cutoff value and a large value $\alpha_{\text{large}}$ (flux is favored) for the remaining areas. In this sense, the $a$MRCM controls the relative importance of each interface space at each location.

To solve the Equations (\ref{Darcy0}) or (\ref{Darcy}) by the MRCM we need to set the interface parameters (i.e., the interface spaces and the adaptive function $\alpha(\mathbf {x})$). In 
Algorithm \ref{alg1} we detail the preprocessing operations to set the interface spaces and the adaptive coefficient $\alpha(\mathbf {x})$ from the permeability field. At each interface $\Gamma_{k},\ k=1,\cdots, M$ (where $M$ is the total number of interfaces, horizontal or vertical), and at each interface cell $e_l\in \Gamma_k,\ l=1,\cdots, m_k$ (where $m_k$ is the number of fine cells in $\Gamma_k$), it is defined whether $e_l$ corresponds to a fracture, a barrier, or to background, and then $\alpha({\bf x}_l)$ (at the cell's center) is set accordingly.
 Then, we use the proposed physics-based interface spaces taking into account the observations presented in the previous section. We set the interface pressure (resp., flux) spaces $\mathcal{P}_{H,1}^*$ (respectively, $\mathcal{U}_{H,1}^*$) that consider the physics-based space for the interfaces containing fractures (resp., barriers) and linear spaces for the remaining interfaces.
Finally we are able to combine the $a$MRCM with the physics-based spaces $(\mathcal{U}_{H,1}^*, \mathcal{P}_{H,1}^*)$. We refer to this combination as $a$MRCM-PBS and consider for comparisons, the $a$MRCM-POL, which represents the $a$MRCM combined with linear polynomial spaces independently on the permeability field.

\begin{algorithm}
\caption{Setting interface parameters for solving Equations (\ref{Darcy0}) or (\ref{Darcy}) by the $a$MRCM-PBS}
\label{alg1}
\begin{algorithmic}[1]
\STATE Given $K(\mathbf x)$, $\zeta_{\min}$, $\zeta_{max}$, $\alpha_{\text{small}}$ and $\alpha_{\text{large}}$ 
\FOR{$k=1$ \TO $M$ }
	\FOR{$l=1$ \TO $m_k$}
    	\STATE Evaluate the permeability in both sides ($K^-({\bf x}_l)$ and $K^+({\bf x}_l)$) of the interface cell $e_l\in\Gamma_k $ 
    	\IF {($\max\{K^-({\bf x}_l),K^+({\bf x}_l)\}>\zeta_{\max}$)}
    		\STATE $\alpha({\bf x}_l)=\alpha_{\text{small}}$
    		\STATE Add $\Gamma_{k}$ to $\Gamma^{\text{frac}}$ 
    	\ELSE
    		\STATE $\alpha({\bf x}_l)=\alpha_{\text{large}}$
    	\ENDIF
    	\IF {($\min\{K^-({\bf x}_l),K^+({\bf x}_l)\}<\zeta_{\min} $)}
    		\STATE Add $\Gamma_{k}$ to $\Gamma^{\text{barrier}}$ 
		\ENDIF
	
	\ENDFOR  
  
	\IF {($\Gamma_{k}\subset\Gamma^{\text{frac}}$)}
		\STATE Compute the physics-based functions for pressure from Equations (\ref{new_basis1})-(\ref{new_basis3}) 
		\STATE Set $\mathcal{P}_{H,1}^*=\mathcal{P}_{H}^*$ at $\Gamma_k$	
	\ELSE
		\STATE Compute the linear polynomials functions for pressure	 
		\STATE  Set $\mathcal{P}_{H,1}^*=\mathcal{P}_{H,1}$ at $\Gamma_k$	
	\ENDIF

	\IF {($\Gamma_{k}\subset\Gamma^{\text{barrier}}$)}
		\STATE Compute the physics-based functions for flux from Equations (\ref{new_basis4})-(\ref{new_basis6}) 
		\STATE Set $\mathcal{U}_{H,1}^*=\mathcal{U}_{H}^*$ at $\Gamma_k$	
	\ELSE
		\STATE Compute the linear polynomials functions for flux	 
		\STATE  Set $\mathcal{U}_{H,1}^*=\mathcal{U}_{H,1}$ at $\Gamma_k$	
	\ENDIF 
\ENDFOR
\STATE Given $(\mathcal{U}_{H,1}^*, \mathcal{P}_{H,1}^*)$ and $\alpha(\mathbf {x})$, solve Equations (\ref{Darcy0}) or (\ref{Darcy}) to obtain $p(\mathbf x)$ and $\mathbf u(\mathbf x)$ by the $a$MRCM-PBS
\end{algorithmic}
\end{algorithm}

In the next section, we explore the proposed $a$MRCM-PBS in challenging high-contrast fractured-like porous media through numerical experiments.

\section{Numerical experiments }\label{results}

We present representative numerical experiments to investigate the accuracy of the introduced physics-based interface spaces for the approximation of flows in high-contrast fractured-like porous media.
We start with the investigation of the physics-based interface spaces combined with the $a$MRCM for the elliptic problem. Then we study the influence of the physics-based interface spaces in the approximation of two-phase problems.

In all simulations the fine grid solution is used as the reference solution for evaluating the multiscale approximations. We consider the $a$MRCM (by setting $\alpha_{\text{small}}=10^{-2}$ and $\alpha_{\text{large}}=10^{2}$), MMMFEM (by setting $\alpha=10^{-6}$) and MHM (by setting $\alpha=10^{6}$). 
To recover a continuous flux at the fine-scale we consider a velocity post-processing (downscaling)
\cite{guiraldello2019downscaling}. We choose the Stitch method, which has been indicated as the procedure that provides the best compromise between computational cost and precision.
The patch thicknesses of the oversampling regions for the Stitch method is fixed in two elements (that represents $10\%$ of the size of most subdomains considered in the examples).

\subsection{Single-phase flows in permeability fields with fractures and barriers} 

\begin{figure}[htbp]
    \centering
    \includegraphics[scale=0.36]{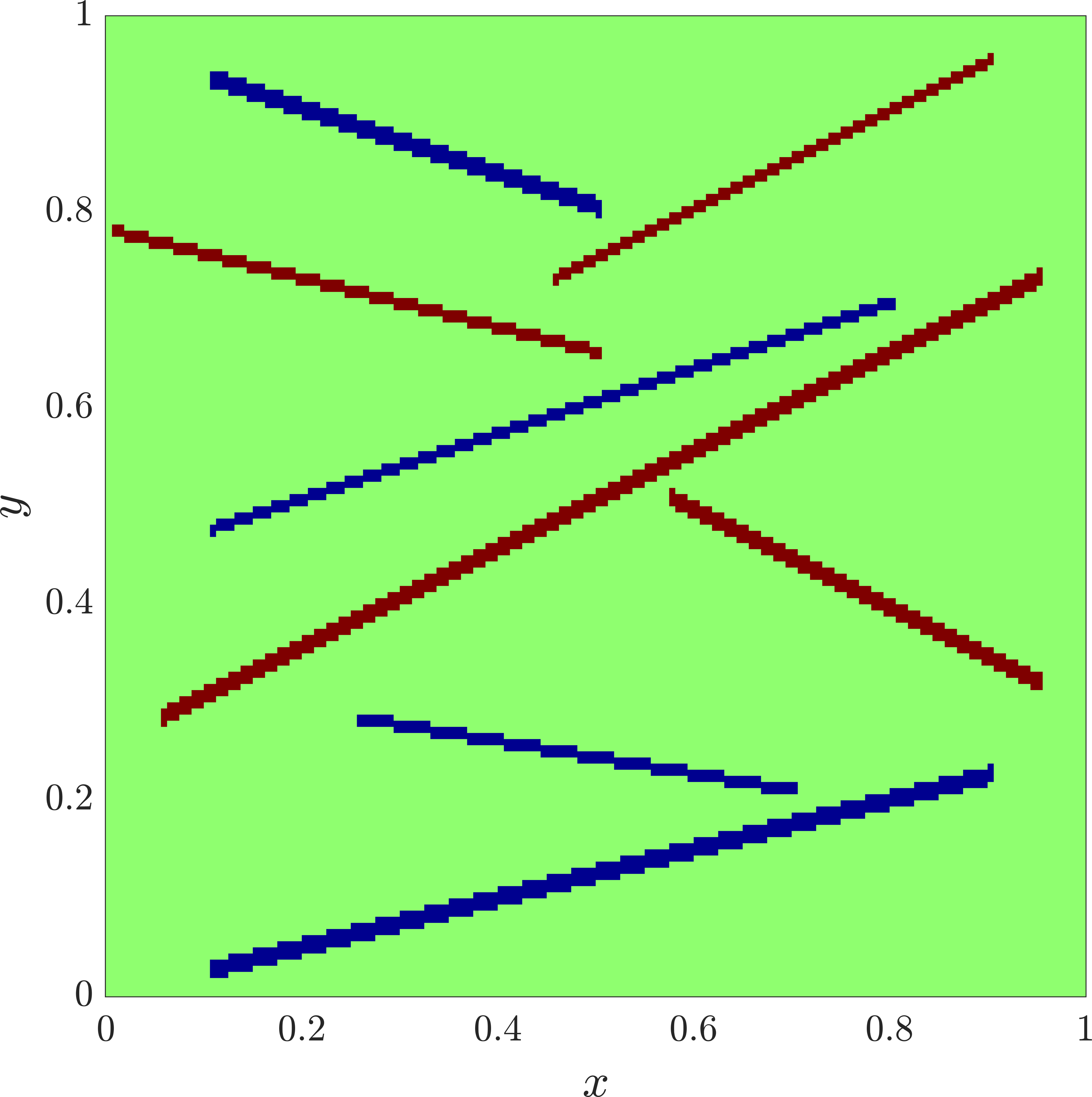}
      \caption{Permeability field containing both fractures (red) and barriers (blue). The contrast considered vary from $10$ to $10^8$.}
      \label{fig:perm_frac_and_barriers}
\end{figure}

We consider a permeability field containing both fractures and barriers, see Figure \ref{fig:perm_frac_and_barriers}. We fix the homogeneous background permeability with $K=1$ and consider the permeability in the barriers, ${K_{\min}}$ varying from $10^{-4}$ to $0.5$ whereas the permeability in the fractures, ${K_{\max}}$ varying from $10$ to $10^{4}$. The cutoff values considered are $\zeta_{\max}=\zeta_{\min}=1$.

\begin{figure}[htbp]
    \centering
        \includegraphics[scale=0.37]{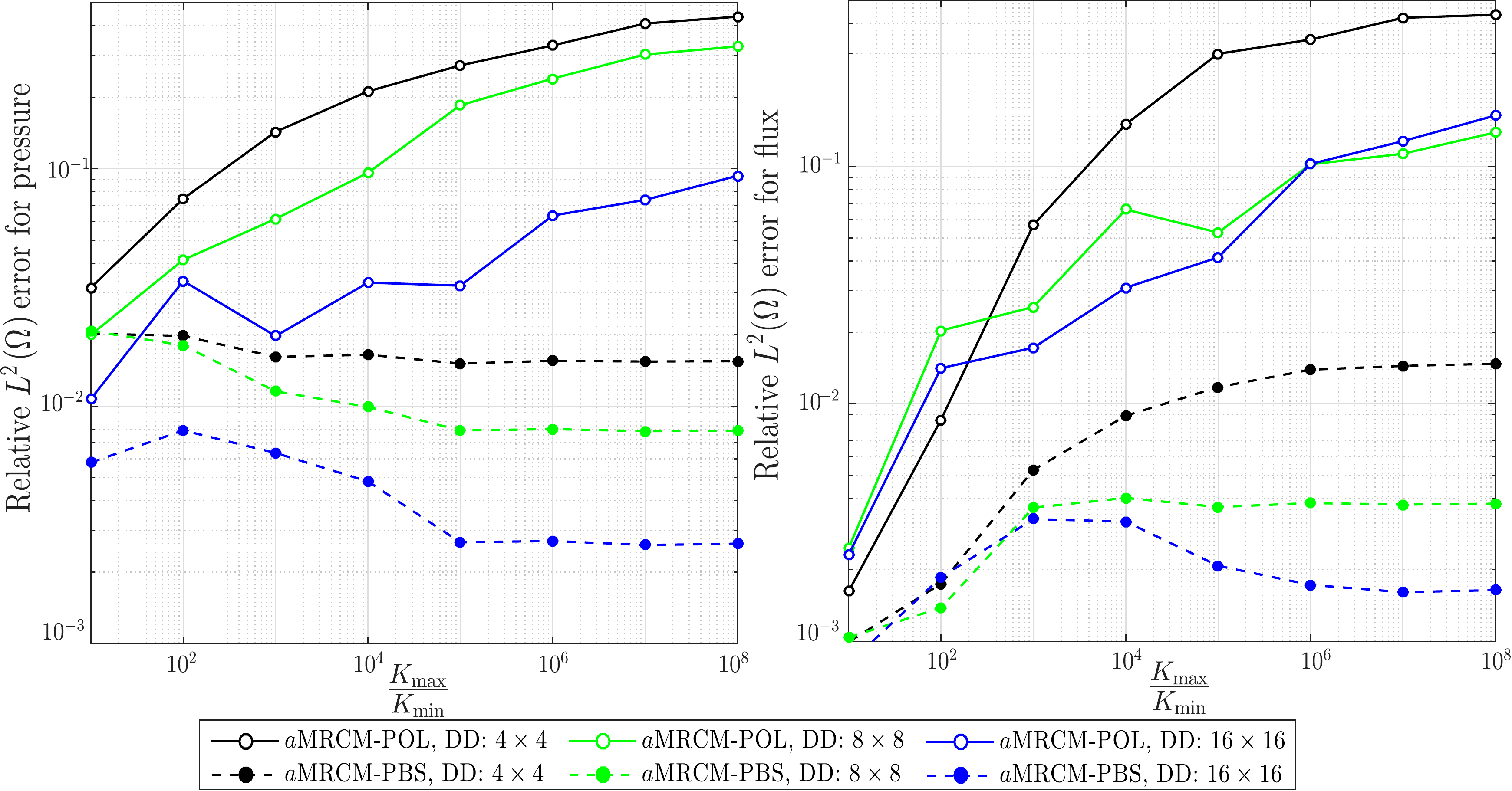}
       \caption{Relative $L^2(\Omega)$ pressure (left) and flux (right) errors as function of the contrast are shown for the $a$MRCM-POL and $a$MRCM-PBS. Three domain decompositions are considered: with $4\times4$, $8\times8$ and $16\times16$ coarse cells. We note that the $a$MRCM-PBS is more accurate than the $a$MRCM-POL in all the contrast and meshes considered.}
       \label{fig:erro_frac_and_barriers}
\end{figure}

\begin{figure}[htbp]
    \centering
        \includegraphics[scale=0.38]{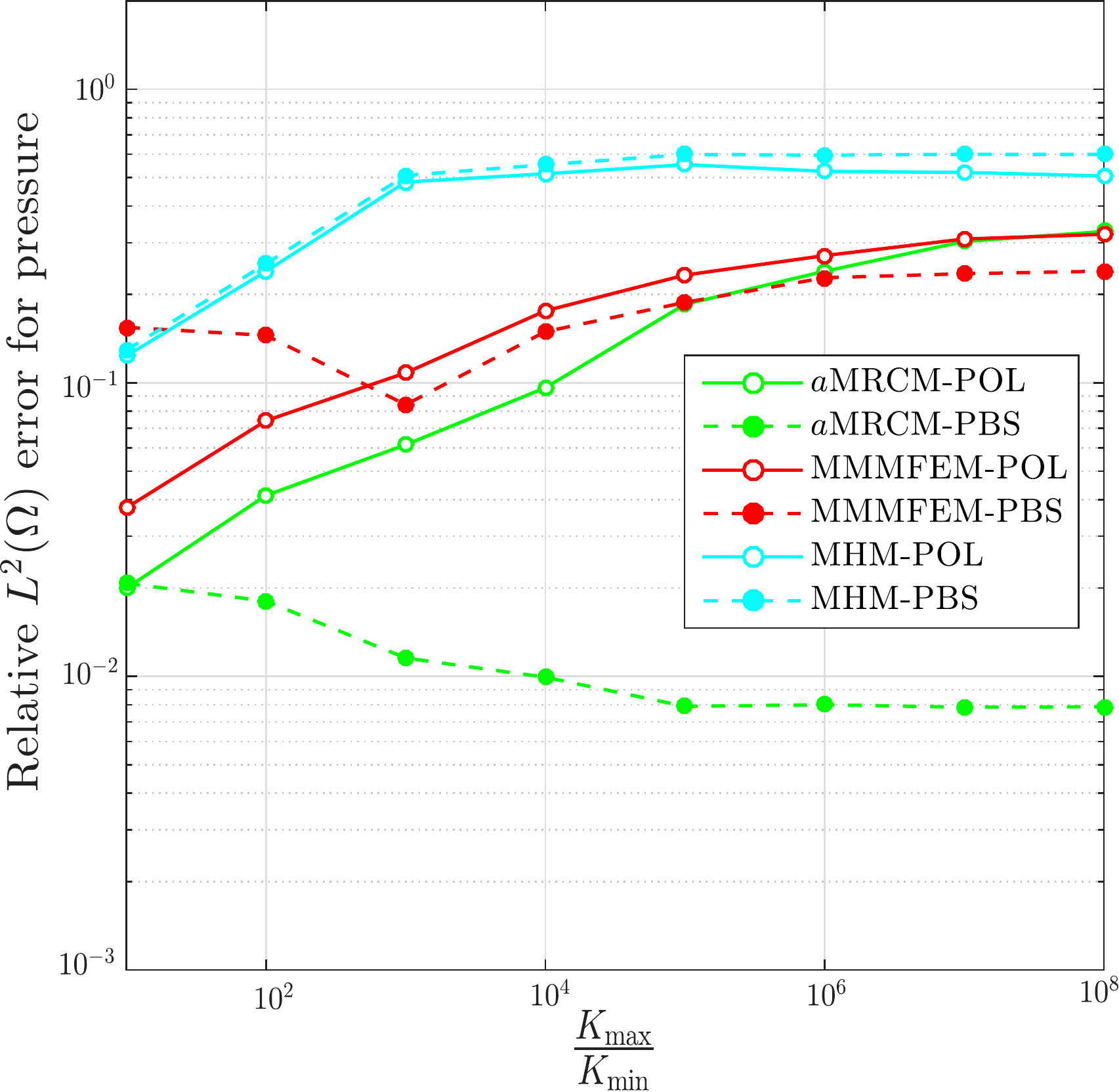}
    \includegraphics[scale=0.38]{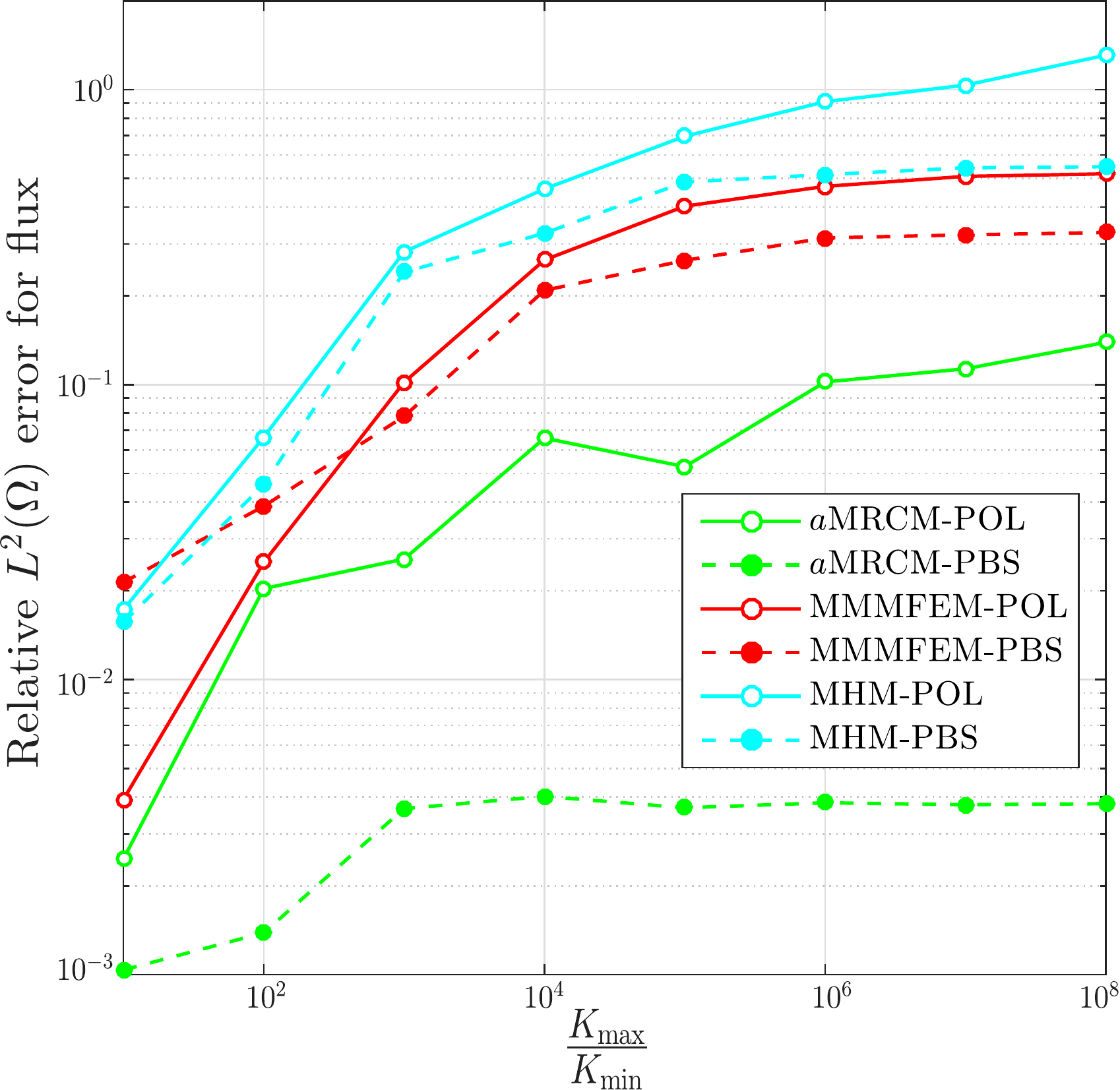}
       \caption{Relative $L^2(\Omega)$ pressure (left) and flux (right) errors as function of the contrast are shown for the $a$MRCM, MMMFEM and MHM considering both the linear and physics-based spaces. Three domain decompositions are considered: with $4\times4$, $8\times8$ and $16\times16$ coarse cells. We note that the $a$MRCM-PBS is more accurate than the $a$MRCM-POL in all the contrast and meshes considered.}
       \label{fig:comp_frac_and_barriers}
\end{figure}

\begin{figure}[htbp]
    \centering 
        \includegraphics[scale=0.155]{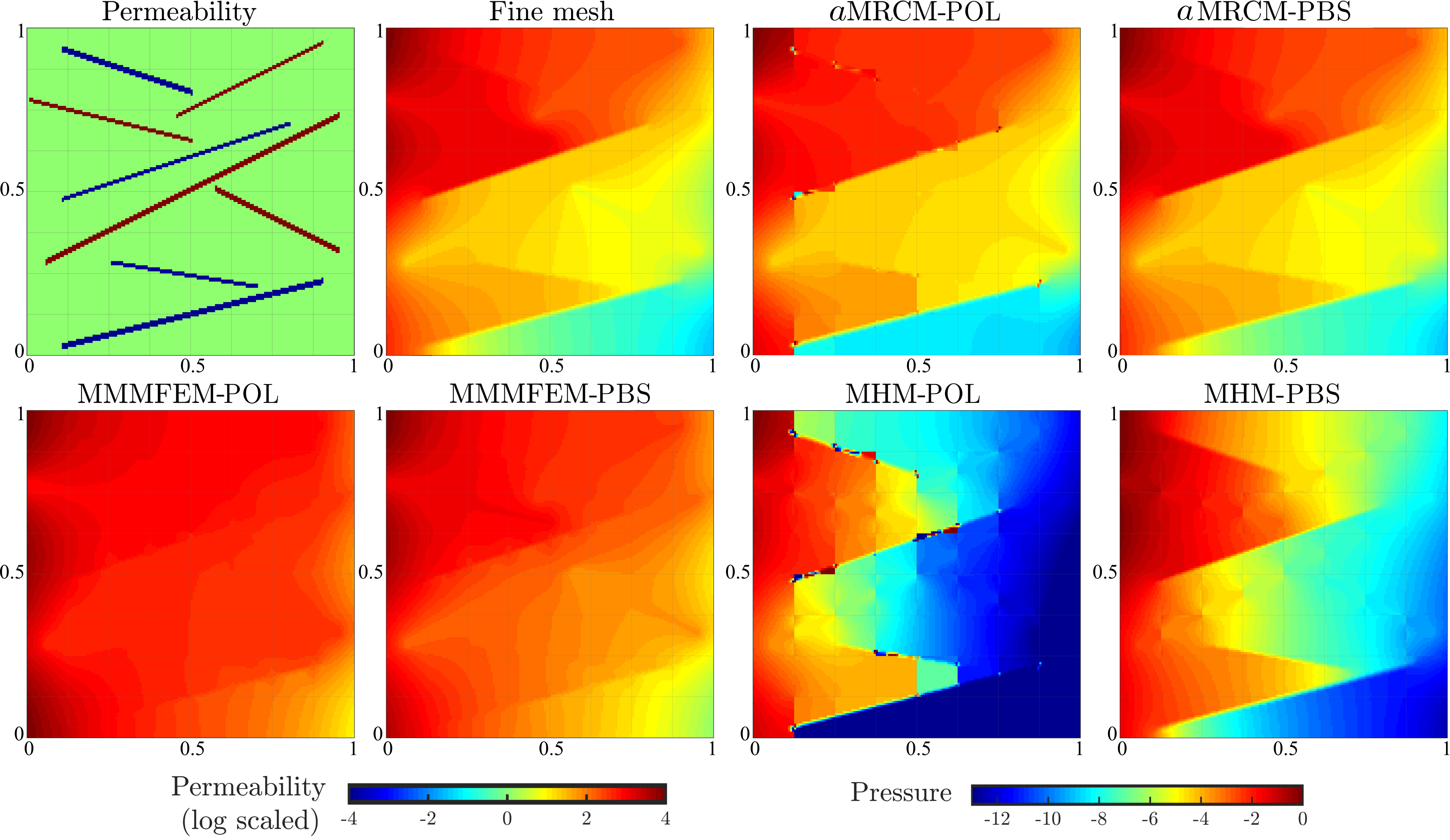}
       \caption{Pressure approximations considering the contrast of ${K_{\max}}/{K_{\min}}=10^8$. First line, left to right: permeability field, fine mesh, $a$MRCM-POL, and $a$MRCM-PBS solutions. Second line, left to right: MMMFEM-POL, MMMFEM-PBS, MHM-POL, and MHM-PBS solutions. The domain decomposition considered contains $8\times8$ coarse cells and is illustrated by the lines in the plot. We note that the $a$MRCM-PBS approximation is the most closely related to the reference solution followed by the $a$MRCM-POL.}
       \label{fig:press_frac_and_barriers}
\end{figure}

\begin{figure}[htbp]
    \centering 
    \includegraphics[scale=0.156]{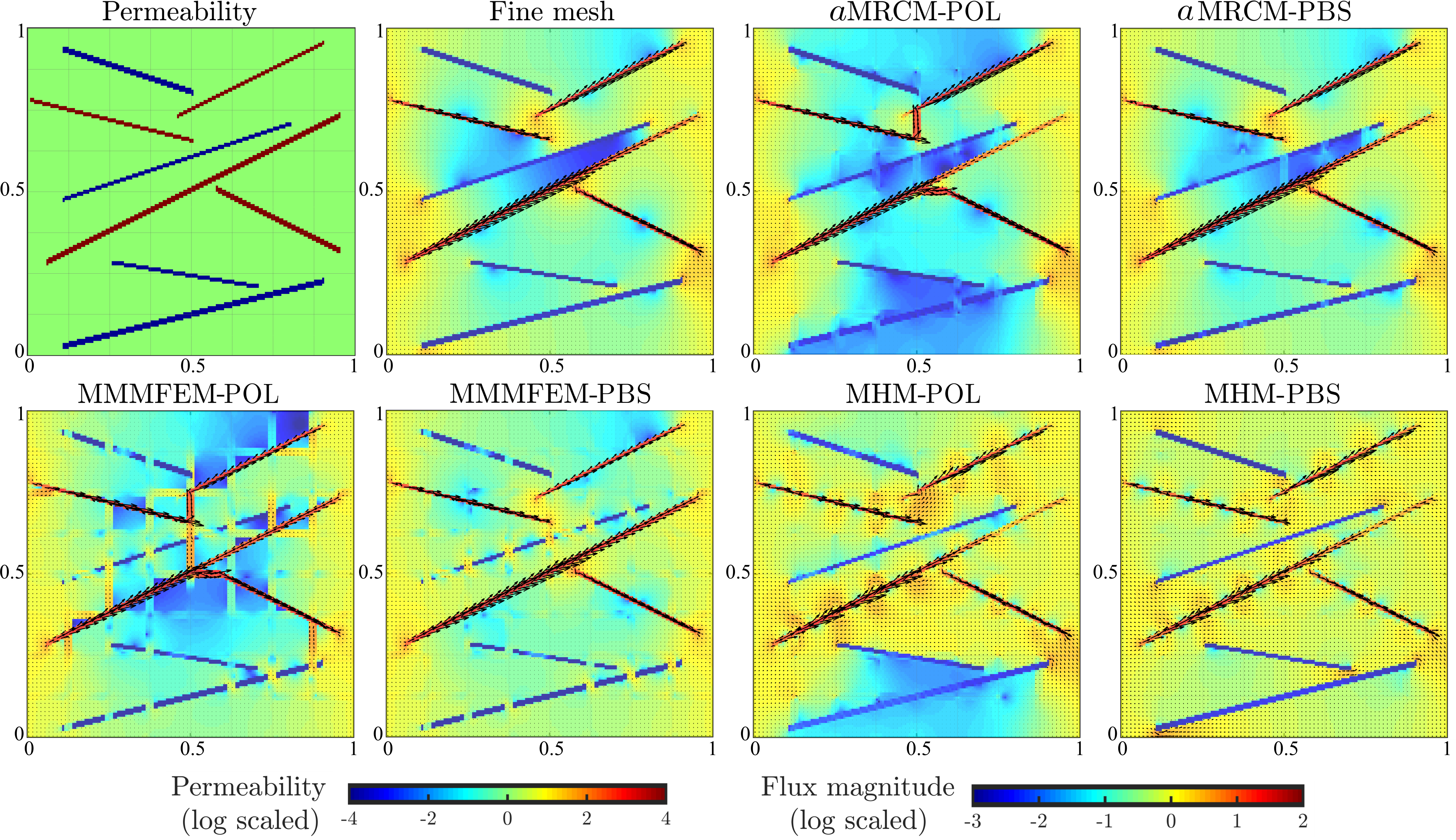}
       \caption{Flux approximations considering the contrast of ${K_{\max}}/{K_{\min}}=10^8$. First line, left to right: permeability field, fine mesh, $a$MRCM-POL, and $a$MRCM-PBS solutions. Second line, left to right: MMMFEM-POL, MMMFEM-PBS, MHM-POL, and MHM-PBS solutions. The colors refer to the log-scale flux magnitude. The domain decomposition considered contains $8\times8$ coarse cells and is illustrated by the lines in the plot. The only accurate procedure is the $a$MRCM-PBS.}
       \label{fig:flux_frac_and_barriers}
\end{figure}

The relative $L^2(\Omega)$ error norms for pressure and flux as function of the contrast are shown in Figure \ref{fig:erro_frac_and_barriers} for three domain decompositions: with $4\times4$, $8\times8$ and $16\times16$ coarse cells. We observe that the pressure and flux errors related to the linear interface spaces increase quickly with the contrast. However, the errors provided by the physics-based spaces are moderate even for very large permeability contrasts. The results are consistent for all the domain decomposition tested, where the smaller errors are attained by the decompositions with more subdomains.
In Figure \ref{fig:comp_frac_and_barriers} we compare the solutions provided by the $a$MRCM, MMMFEM and MHM considering the domain decomposition of $8\times8$ subdomains having $20\times20$ fine cells into each one. 
The first observation is that the MMMFEM and MHM approximations are not accurate for the permeability field with fractures and barriers even using the physics-based interface spaces. On the other hand, the $a$MRCM-PBS solutions are expressively more accurate than the ones obtained with the $a$MRCM-POL.  
The pressure and flux approximations for the permeability contrast of $10^{8}$ are shown in Figure \ref{fig:press_frac_and_barriers} and Figure \ref{fig:flux_frac_and_barriers}, respectively. The plots confirm that the $a$MRCM combined with the physics-based interface spaces produces the most accurate solutions. We observe that all methods fail when using the linear interface spaces. We can conclude that the physics-based interface spaces are indispensable for a reasonable approximation of pressure, flux and further applications in two-phase flows.

\begin{figure}[htbp]
    \centering 
    \includegraphics[scale=0.37]{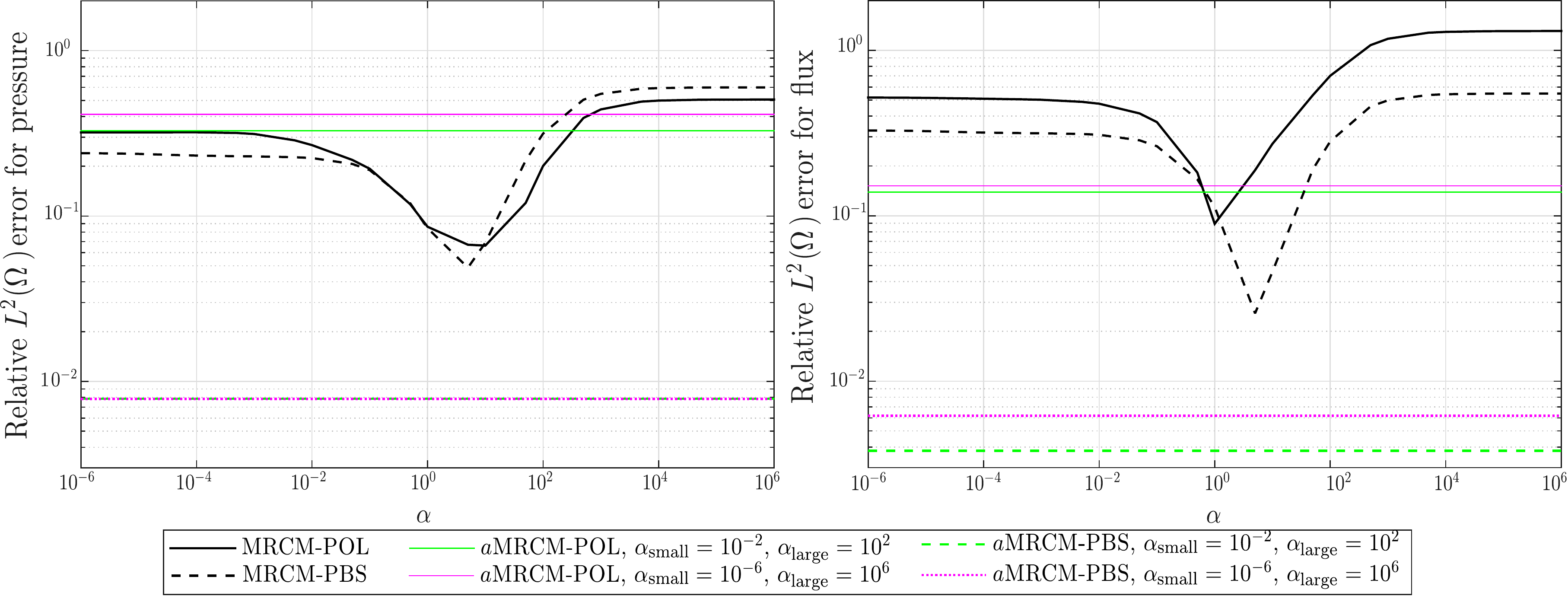}
       \caption{Relative $L^2(\Omega)$ errors as a function of $\alpha$ for pressure (left) and flux (right) considering the permeability field plotted in Figure \ref{fig:perm_frac_and_barriers} with contrast of $10^8$ and the domain decomposition that contains $8\times8$ coarse cells. The physics-based and linear spaces are compared. We include the errors for two choices in the $a$MRCM: setting $\alpha_{\text{small}}=10^{-2}$ and $\alpha_{\text{large}}=10^{2}$ or $\alpha_{\text{small}}=10^{-6}$ and $\alpha_{\text{large}}=10^{6}$ (illustrated as horizontal lines). The improvement with the combination of the $a$MRCM and the physics-based spaces is expressive.}
       \label{fig:varrido_frac_and_barriers}
\end{figure}

Finally, we present in Figure \ref{fig:varrido_frac_and_barriers} 
the behavior of the errors for pressure and flux as a function of $\alpha$ (varying from $10^{-6}$ to $10^6$) for the permeability contrast of $10^8$ and maintaining the domain decomposition of $8\times8$ coarse cells. We compare the MRCM errors with linear and physics-based interface spaces. Two choices for the $a$MRCM are presented: setting $\alpha_{\text{small}}=10^{-2}$ and $\alpha_{\text{large}}=10^{2}$ or $\alpha_{\text{small}}=10^{-6}$ and $\alpha_{\text{large}}=10^{6}$. The errors for the adaptive version of the MRCM are illustrated as horizontal lines. We note a strong dependence on the parameter $\alpha$, where the minimum errors are attained at intermediate values (for choosing $\alpha$ constant), similarly to the previous works \cite{guiraldello2018multiscale, guiraldello2019interface, bifasico}. But the approximations by choosing any constant value of $\alpha$, even for the MRCM-PBS are inaccurate. We remark that the MMMFEM and MHM are also included in this observation. One can conclude that, besides the physics-based interface spaces, the $a$MRCM is necessary to obtain more accurate solutions. We note that the choice of $\alpha_{\text{small}}=10^{-2}$ and $\alpha_{\text{large}}=10^{2}$ or $\alpha_{\text{small}}=10^{-6}$ and $\alpha_{\text{large}}=10^{6}$ does not affect significantly the error.

These results indicate that the linear interface spaces are not robust for capturing the effects of features as fractures and barriers. Even the physics-based spaces are not enough for complex fields if not combined with an appropriate multiscale method. 
In order to study how these results are reflected in the corresponding saturation fields, in the next subsection we present numerical results for two-phase flows.


\subsection{Two-phase flow and transport problems}  
 
Now we focus on a study of the MRCM performance for two-phase flows using the physics-based interface spaces for high-contrast fractured-like permeability fields. We introduce the saturation comparisons with a detailed analysis of fingering instabilities to show the impact of the design of the basis functions on the transport of the water saturation. Firstly we present numerical results for the permeability field with fractures and barriers of the previous section. Then we consider a high-contrast permeability field with channels and isolated inclusions that has frequently appeared in the literature \cite{jiang2017model}. 

In all simulations the reservoir is initially filled with oil ($s^0=0$) and water is injected at a constant rate.
The relative permeabilities are given by $k_{ro}=(1-s)^2$ and $k_{rw}=s^2$, such that the fractional flow is written as
\begin{equation}
f(s) = \dfrac{M s^2}{M s^2 + (1-s)^2},
 \label{frac_flow}
\end{equation}
where $ M = {\mu_o}/{\mu_w}$ is set as $M=10$.

In the operator splitting approach, we take the number of transport steps between the elliptic updates at most $C=20$. The time units employed are in pore volumes injected (PVI) that refers to the fraction of the total accessible pore volume injected into the domain \cite{chen2006computational}
\begin{equation}
T_{\text{PVI}}=-V_p^{-1} \int_0^t  \int_{\partial\Omega_{\text{in}}}\mathbf{u}(\mathbf{x},\tau)\cdot \mathbf{n}\ dl\  d\tau,
\end{equation}
where $V_p$ is the total pore-volume of the reservoir, $t$ is the time taken for injection and $\partial\Omega_{\text{in}}$ the inflow well boundaries with the outward unit normal $\mathbf{n}$.

\subsubsection{High-contrast permeability field with fractures and barriers} 

The objective of this study is to investigate the $a$MRCM, MMMFEM and MHM saturation solutions considering the physics-based interface spaces for the high-contrast permeability field with fractures and barriers. We consider the permeability field of the previous experiment, illustrated in Figure \ref{fig:perm_frac_and_barriers}. We fix the highest permeability contrast ${K_{\max}}/{K_{\min}}=10^8$ and the domain decomposition having $8\times8$ subdomains, each one discretized by $20\times20$ cells.

\begin{figure}[htbp]
    \centering
    \includegraphics[scale=0.52]{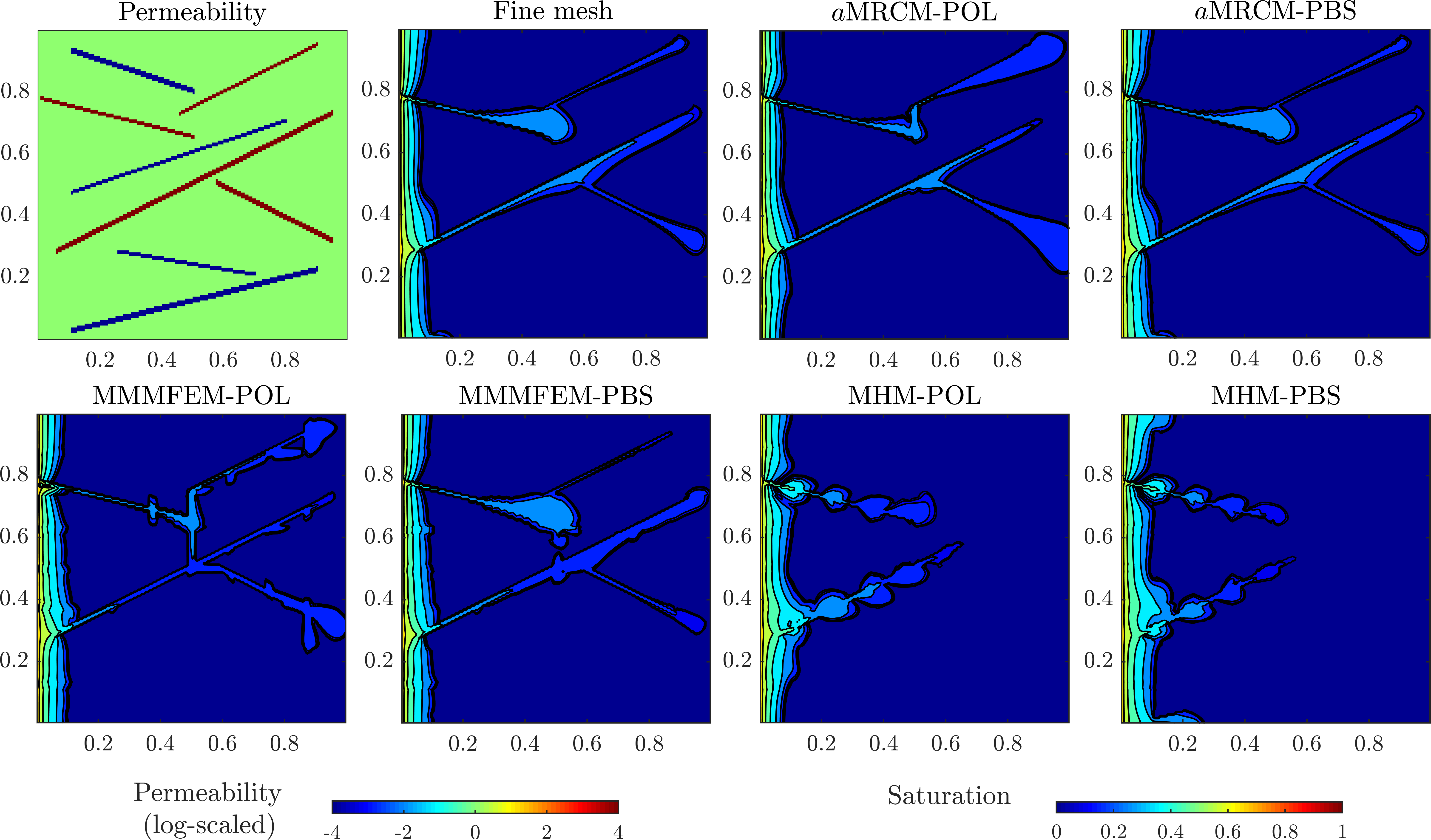}
       \caption{Comparison of multiscale methods for the slab geometry with flux boundary conditions on the left and right. Saturation profiles at $T_{PVI}=0.06$ for the permeability field with fractures and barriers are shown. First line, left to right: high-contrast permeability field (log-scaled); reference fine grid solution; $a$MRCM-POL saturation profile; $a$MRCM-PBS saturation profile. Second line, left to right: MMMFEM-POL saturation profile; MMMFEM-PBS saturation profile; MHM-POL saturation profile; MHM-PBS saturation profile. The $a$MRCM-PBS provides the most accurate approximation.}
       \label{fig:sat_maps_slab_N}
\end{figure}

\begin{figure}[htbp]
    \centering
    \includegraphics[scale=0.55]{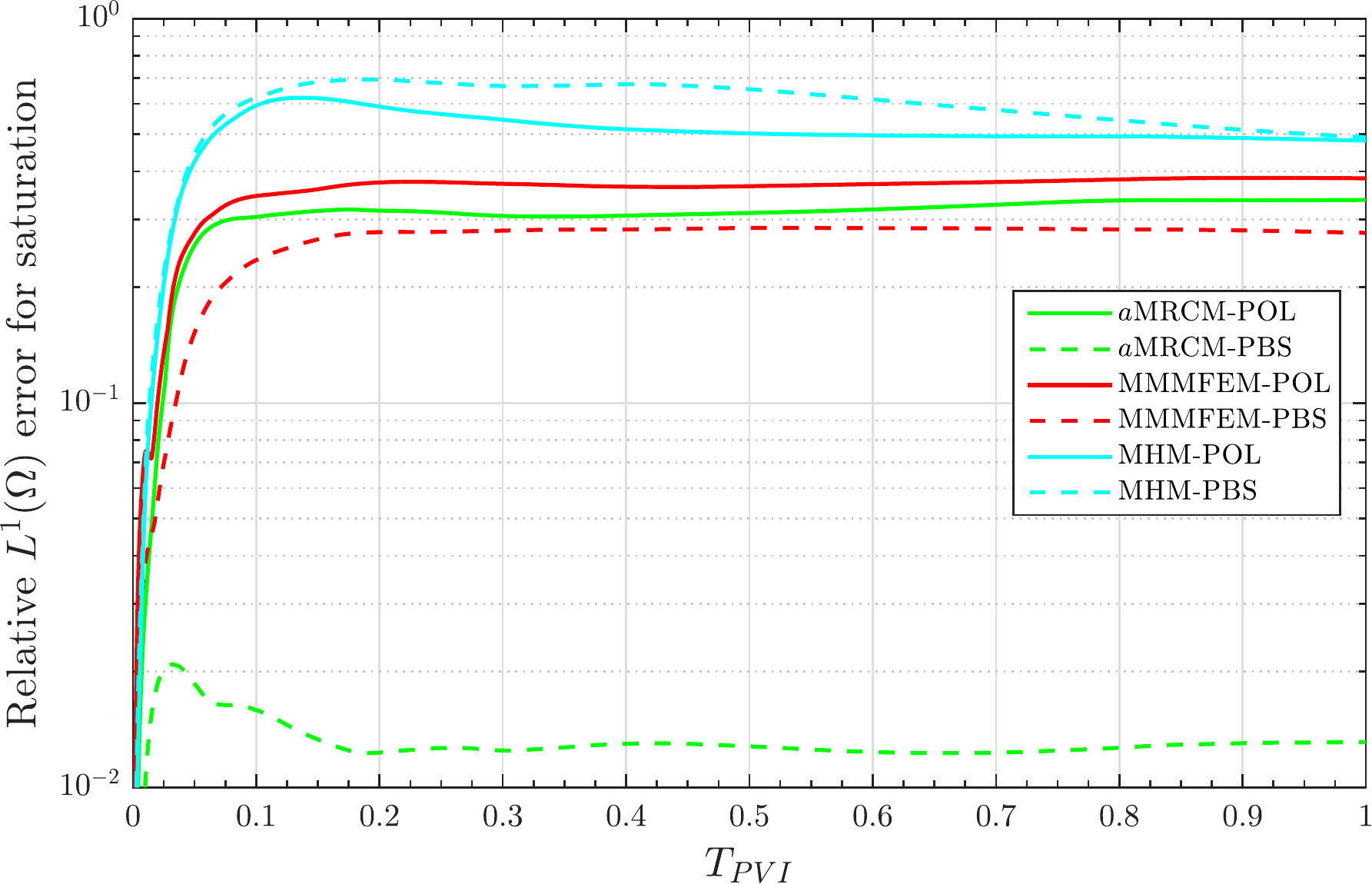} 
    \caption{Relative $L^{1}(\Omega)$ saturation errors as a function of time for the slab geometry with flux boundary conditions considering the field with fractures and barriers. We compare the $a$MRCM, MMMFEM and MHM with both physics-based and linear spaces. We note that the errors associated with the $a$MRCM-PBS are much smaller than all the others.}
       \label{fig:sat_errors_slab_N}
\end{figure}

The first example considers the same slab geometry of the previous experiments, by imposing the flux on the left and right boundaries and no source terms. Figure \ref{fig:sat_maps_slab_N} shows the permeability field (log-scaled) along with the saturation profiles at $T_{\text{PVI}}=0.06$ (before breakthrough time) approximated by the multiscale methods with the linear and the physics-based interface spaces. The procedure that produces a saturation solution closer to the reference one is the $a$MRCM-PBS. The use of the physics-based spaces enables more accurate solutions for the $a$MRCM and the MMMFEM. However, the last one is still inaccurate. The MHM solutions are unacceptable either for linear or physics-based spaces. The corresponding relative $L^1(\Omega)$ errors throughout the simulation are presented in Figure \ref{fig:sat_errors_slab_N}. We note an expressive improvement provided by the physics-based interface spaces combined with the $a$MRCM. This combination enables the error to drop by one order of magnitude.

\begin{figure}[htbp]
    \centering
    \includegraphics[scale=0.518]{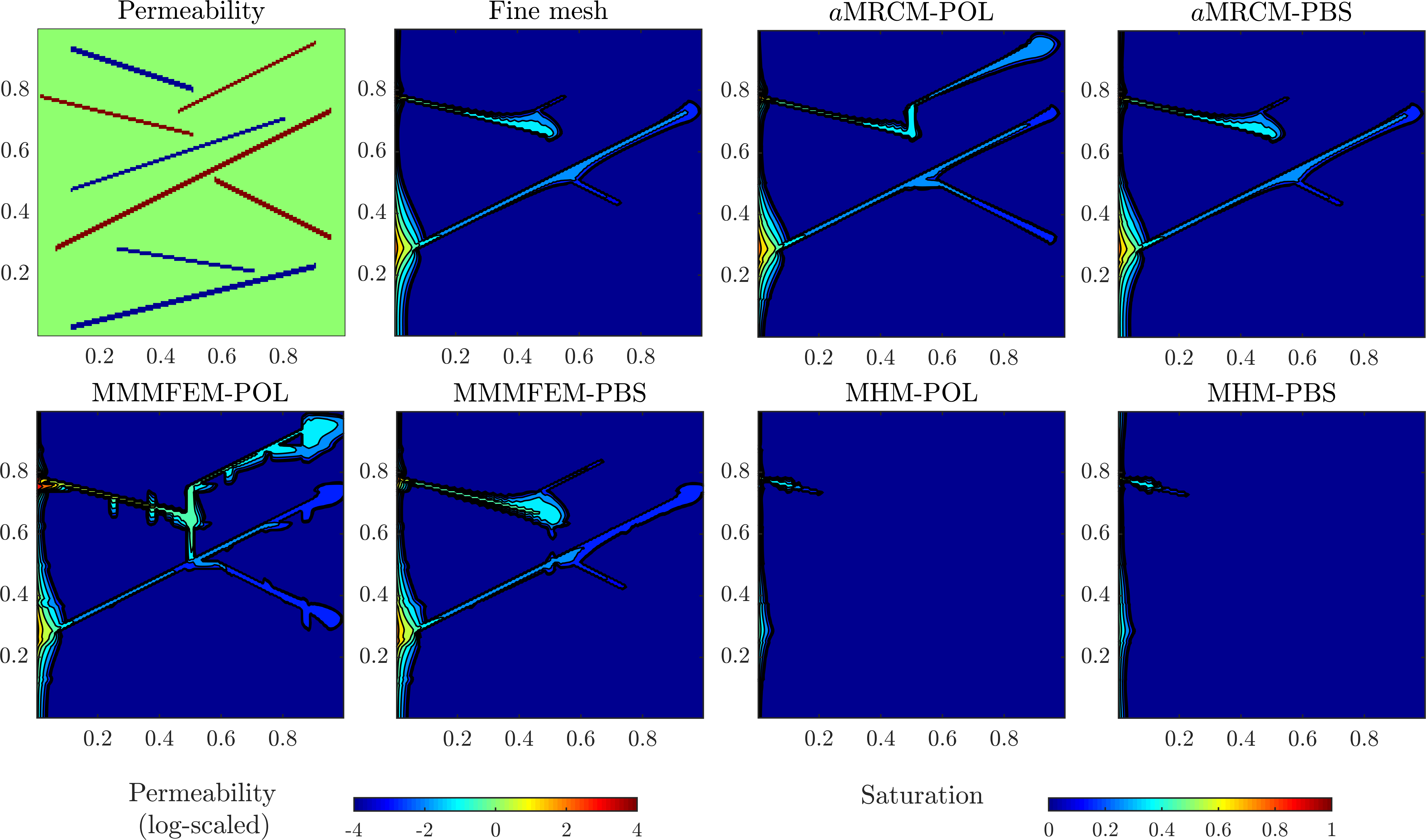} 
       \caption{Comparison of multiscale methods for the slab geometry with pressure boundary conditions on the left and right. Saturation profiles at $T_{PVI}=0.0001$ for the permeability field with fractures and barriers are shown. First line, left to right: high-contrast permeability field (log-scaled); reference fine grid solution; $a$MRCM-POL saturation profile; $a$MRCM-PBS saturation profile. Second line, left to right: MMMFEM-POL saturation profile; MMMFEM-PBS saturation profile; MHM-POL saturation profile; MHM-PBS saturation profile. The $a$MRCM-PBS is clearly the most accurate procedure.}
       \label{fig:sat_maps_slab_D}
\end{figure}

\begin{figure}[htbp]
    \centering 
     \includegraphics[scale=0.55]{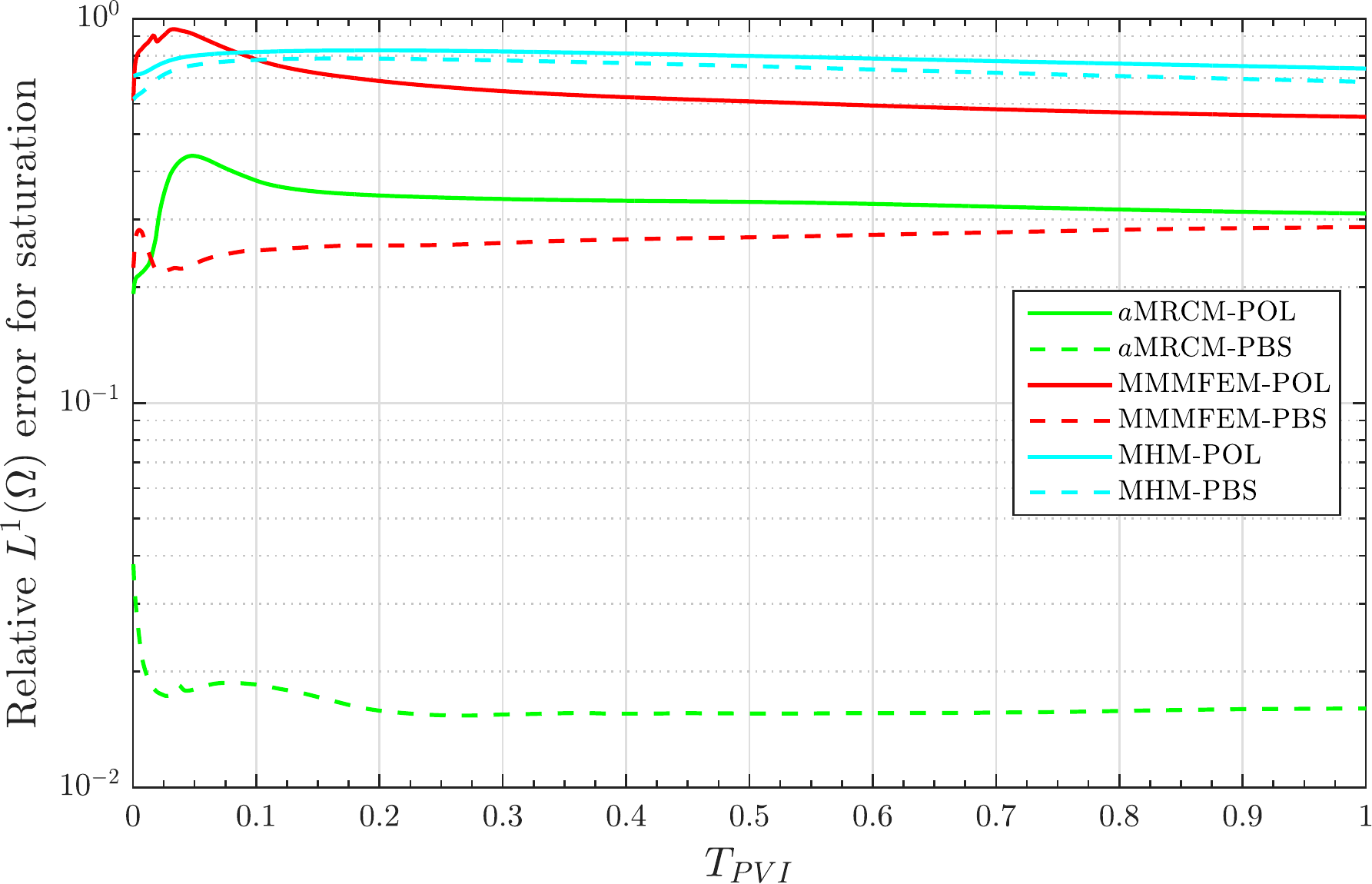}      
       \caption{Relative $L^{1}(\Omega)$ saturation errors as a function of time for the slab geometry with pressure boundary conditions considering the field with fractures and barriers. We compare the $a$MRCM, MMMFEM and MHM with both physics-based and linear spaces. Similarly to the previous example, the errors associated with the $a$MRCM-PBS are much smaller than all the others.}
       \label{fig:sat_errors_slab_D}
\end{figure}

Now we maintain the slab geometry but using global boundary conditions of no-flow at top and bottom boundaries along with imposed pressure on the left and right boundaries.
In Figure \ref{fig:sat_maps_slab_D} we present the saturation profiles at $T_{\text{PVI}}=0.0001$ (before breakthrough time) approximated by the multiscale methods with the linear and the physics-based interface spaces. Similar results to the previous example are obtained. The corresponding relative errors throughout the simulation are presented in Figure \ref{fig:sat_errors_slab_D}, where it is clear the improved accuracy provided by combining the physics-based spaces with the $a$MRCM. Again, the MHM solutions are not acceptable with both interface spaces. The MMMFEM approximations improve significantly with the physics-based spaces, however, these solutions are comparable to the $a$MRCM-POL approximations.

\begin{figure}[htbp]
    \centering
     \includegraphics[scale=0.52]{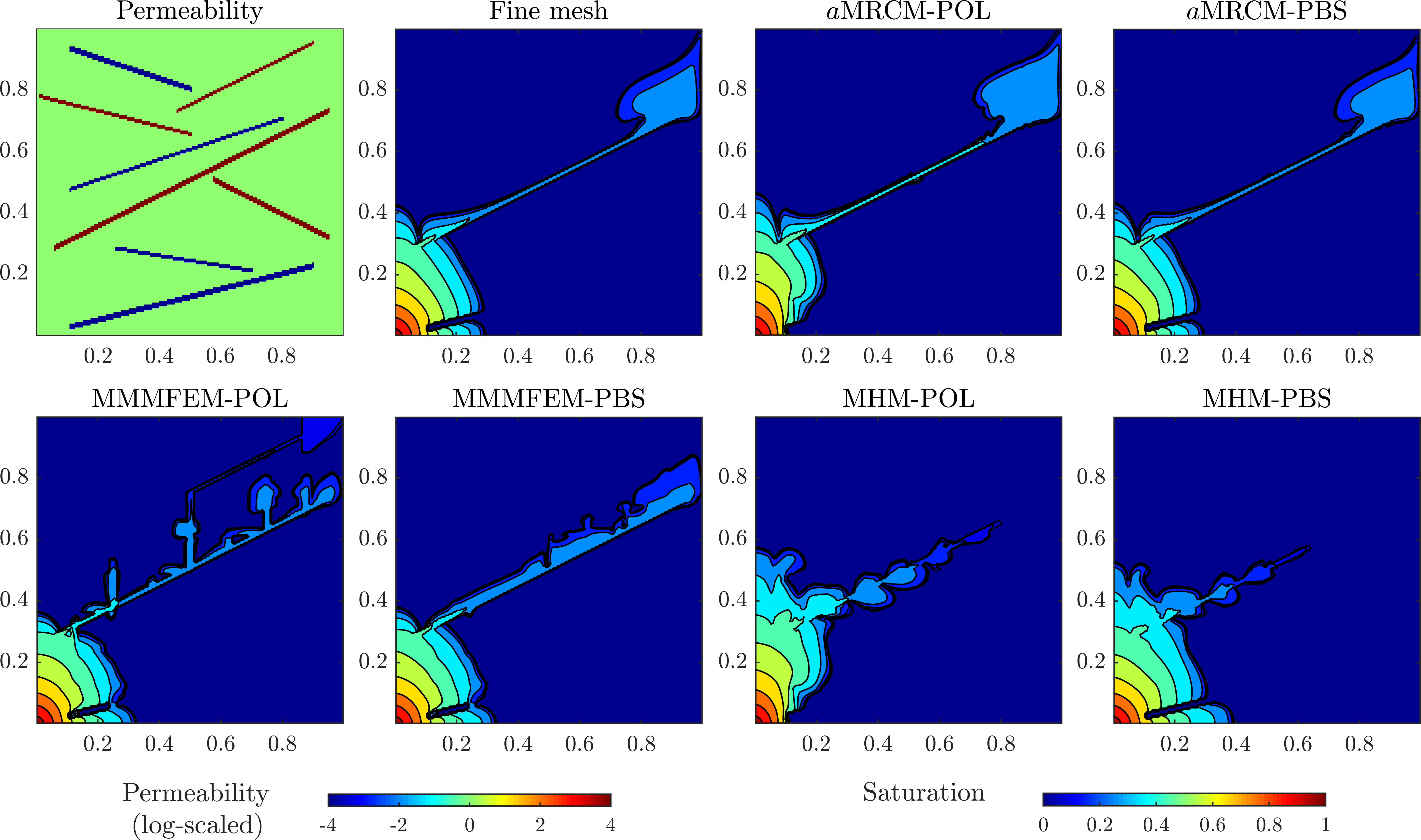} 
      \caption{Comparison of multiscale methods for the quarter of a five-spot geometry. Saturation profiles at $T_{PVI}=0.09$ for the permeability field with fractures and barriers are shown. First line, left to right: high-contrast permeability field (log-scaled); reference fine grid solution; $a$MRCM-POL saturation profile; $a$MRCM-PBS saturation profile. Second line, left to right: MMMFEM-POL saturation profile; MMMFEM-PBS saturation profile; MHM-POL saturation profile; MHM-PBS saturation profile. The $a$MRCM is the only procedure that captures the details of the fingers.}
       \label{fig:sat_maps_5spot}
\end{figure}

\begin{figure}[htbp]
    \centering 
     \includegraphics[scale=0.55]{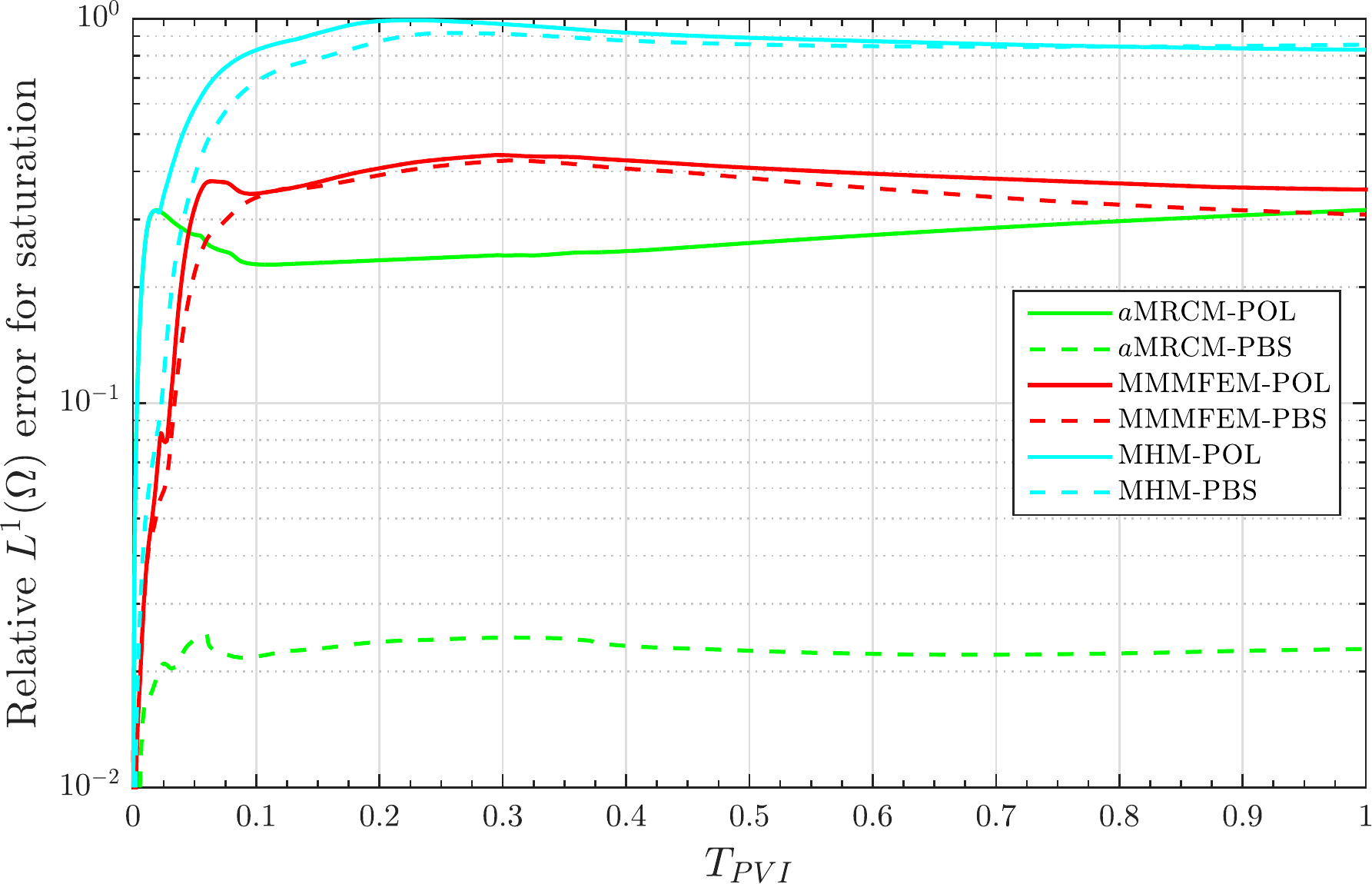}
       \caption{Relative $L^{1}(\Omega)$ saturation errors as a function of time for the quarter of a five-spot geometry on the field with fractures and barriers. We compare the $a$MRCM, MMMFEM and MHM with both physics-based and linear spaces. The errors associated with the $a$MRCM-PBS are the smallest, followed by the $a$MRCM-POL.}
       \label{fig:sat_errors_5spot}
\end{figure}

Lastly, we test the multiscale methods in a quarter of a five-spot model, where we inject the water at the bottom-left cell and the sink is located at the top-right cell. The saturation profiles at $T_{\text{PVI}}=0.09$ (before breakthrough time) are shown in Figure \ref{fig:sat_maps_5spot}. The most accurate solutions are produced by the $a$MRCM considering both types of interface spaces.  
The MMMFEM and MHM approximations present expressive fingering instabilities that are not present in the fine grid solution. In Figure \ref{fig:sat_errors_5spot}, the relative errors confirm that the $a$MRCM produces more accurate solutions. It is noticeable the poor quality solutions provided by the MMMFEM and MHM (even with the physics-based spaces), that are less accurate than the $a$MRCM-POL approximation. The solutions produced by the $a$MRCM-PBS are much more accurate than all the other approximations.

The high-contrast permeability fields are challenging for multiscale methods. We show that the methods fail in a field with both fractures and barriers. In all the models previously tested, the $a$MRCM-PBS is the only scheme that provides satisfactory approximations.

\subsubsection{High-contrast permeability field with channels and isolated inclusions}

In this subsection, we consider a high-contrast permeability field with channels and isolated inclusions that has frequently appeared in the literature \cite{efendiev2013generalized, chung2015mixed, chung2016enriched, chung2014adaptive, chung2017online, chung2018online, jiang2017model}. We consider a domain $\Omega$ containing $100\times100$ fine grid cells divided into $5\times5$ subdomains. The boundary conditions are no-flow at the top and bottom boundaries along with an imposed flux on the left and right boundaries. No source terms are considered. The permeability contrast considered is ${K_{\max}}/{K_{\min}} = 10^6$ and the cutoff values are set as $\zeta_{\max}=\zeta_{\min}=1$.

\begin{figure}[htbp]
    \centering
    \includegraphics[scale=0.52]{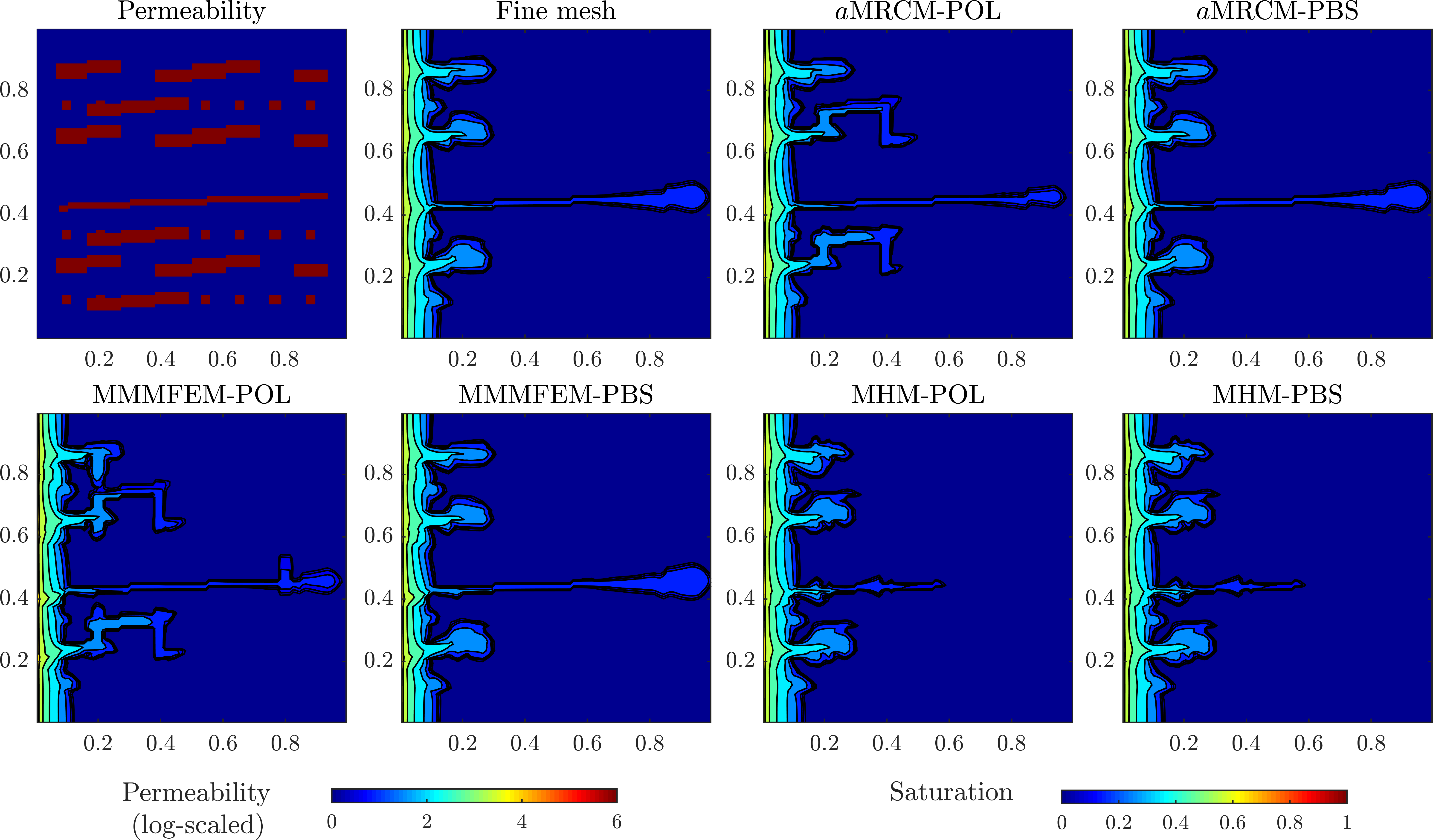} 
    \caption{Comparison of multiscale methods. Saturation profiles at $T_{PVI}=0.07$ for the permeability field with channels and isolated inclusions. First line, left to right: high-contrast permeability field (log-scaled); reference fine grid solution; $a$MRCM-POL saturation profile; $a$MRCM-PBS saturation profile. Second line, left to right: MMMFEM-POL saturation profile; MMMFEM-PBS saturation profile; MHM-POL saturation profile; MHM-PBS saturation profile. The more accurate approximations are attained by the $a$MRCM-PBS and MMMFEM-PBS.}
    \label{fig:sat_high_efendiev}
\end{figure}

\begin{figure}[htbp]
    \centering
    \includegraphics[scale=0.55]{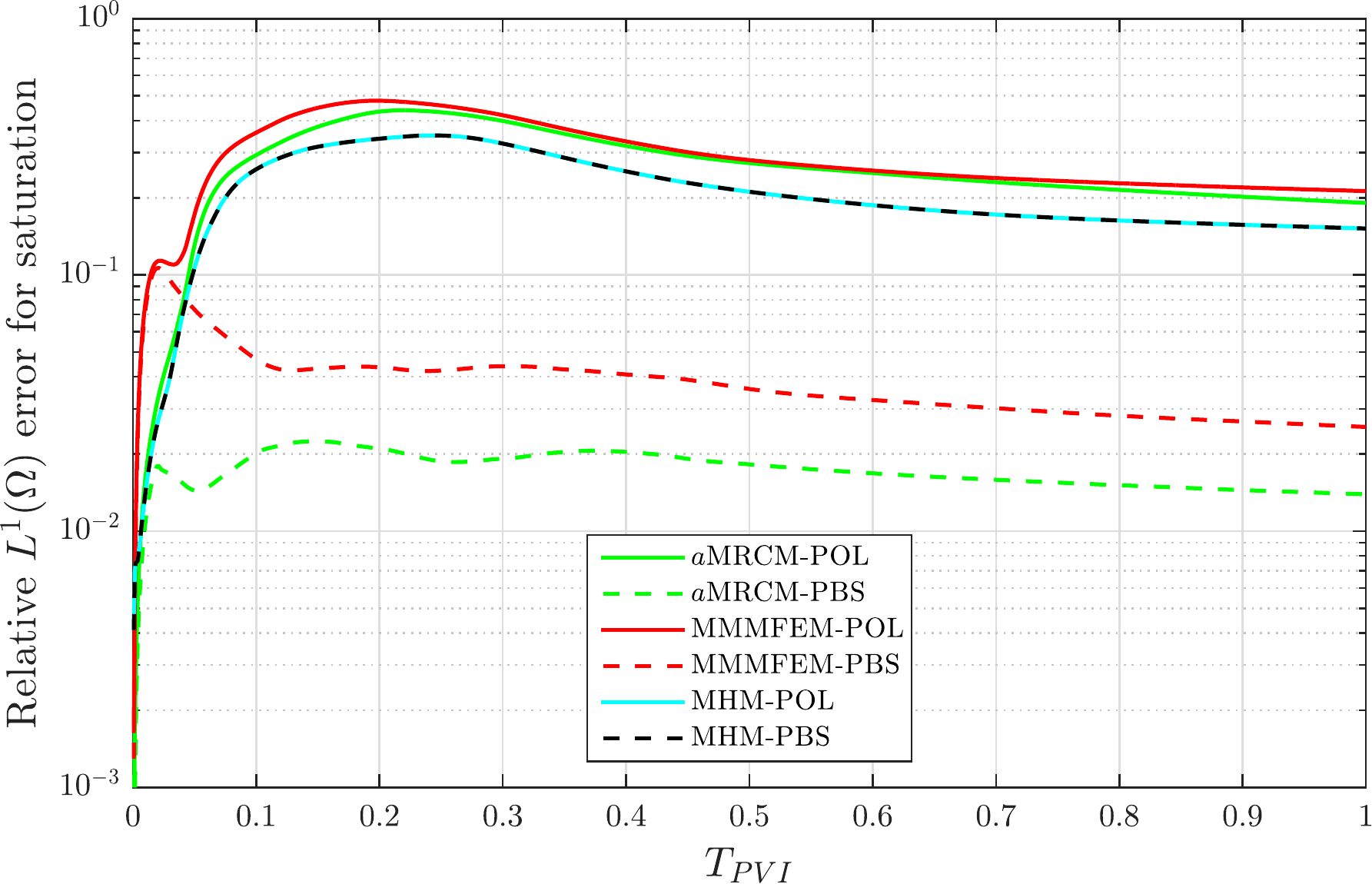} 
        \caption{Relative $L^{1}(\Omega)$ saturation errors as a function of time on the field with channels and isolated inclusions. We compare the $a$MRCM, MMMFEM and MHM with both physics-based and linear spaces. The errors associated with the $a$MRCM-PBS and MMMFEM-PBS are the smallest.}
         \label{fig:sat_erro_high_efendiev}
\end{figure}

Figure \ref{fig:sat_high_efendiev} shows the permeability field (log-scaled) containing high-permeable channels and isolated inclusions and the saturation profiles at $T_{\text{PVI}}=0.07$ (before breakthrough). In this example, only the pressure physics-based spaces are used to handle the high-permeable structures, since there are no barriers. The more accurate solutions are provided by the $a$MRCM-PBS and the MMMFEM-PBS. For these two methods, the figure shows that the imprecisions that happen by using the linear spaces have completely disappeared with the physics-based spaces.
The MHM solutions are inaccurate. We remark that the MHM provides the same solution with the linear and the physics-based interface spaces, because this method considers only the flux space, which is always maintained as linear for this permeability field.
The relative errors are shown in Figure \ref{fig:sat_erro_high_efendiev} and reflect these observations throughout the whole simulation.

\begin{figure}[htbp]
    \centering 
    \includegraphics[scale=0.52]{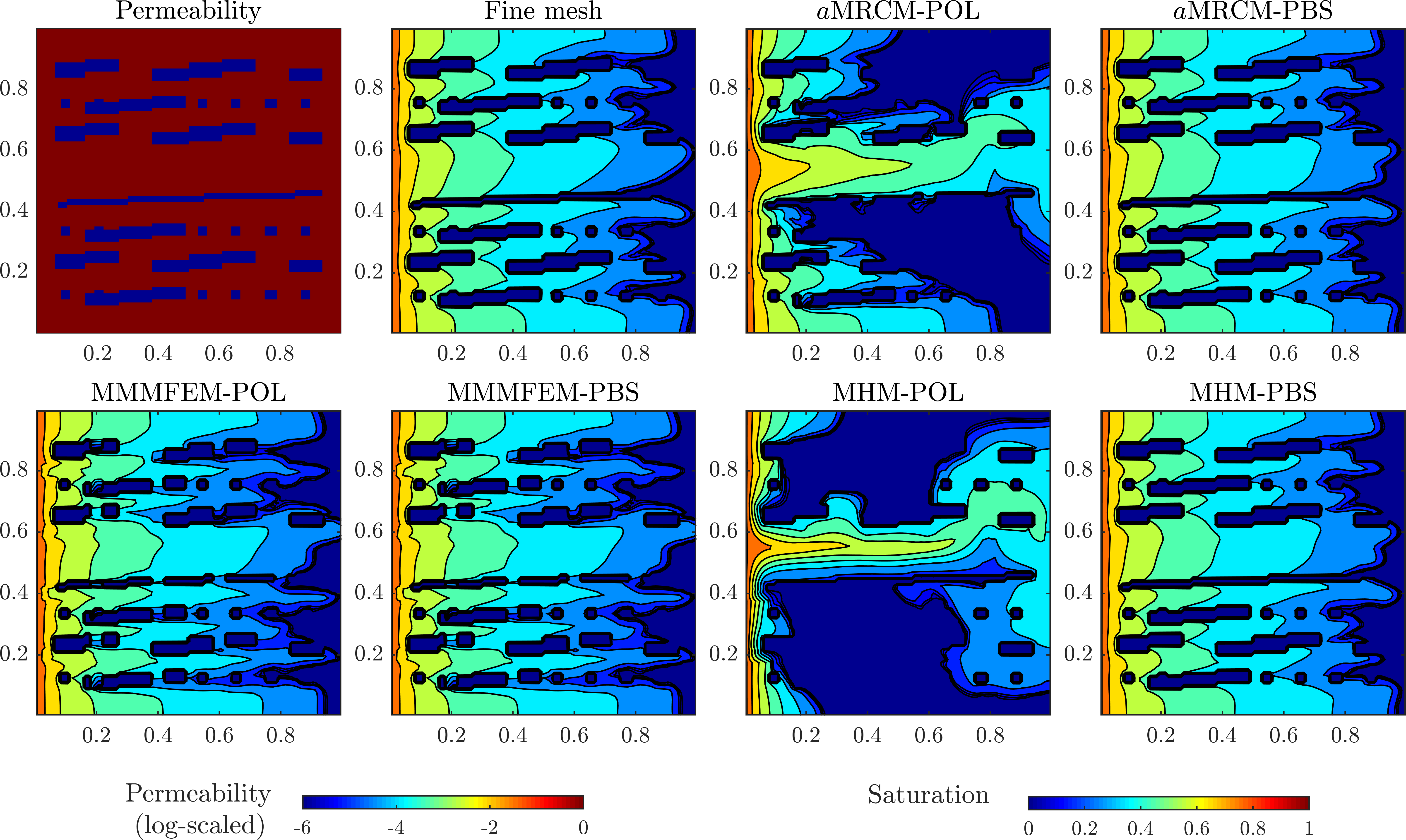} 
    \caption{Comparison of multiscale methods. Saturation profiles at $T_{PVI}=0.33$ for the permeability field with barriers are shown. First line, left to right: high-contrast permeability field (log-scaled); reference fine grid solution; $a$MRCM-POL saturation profile; $a$MRCM-PBS saturation profile. Second line, left to right: MMMFEM-POL saturation profile; MMMFEM-PBS saturation profile; MHM-POL saturation profile; MHM-PBS saturation profile. A considerable improvement is noticed in the $a$MRCM and MHM approximations replacing the linear spaces by the physics-based.}
     \label{fig:sat_low_efendiev}
\end{figure}

\begin{figure}[htbp]
    \centering
      \includegraphics[scale=0.55]{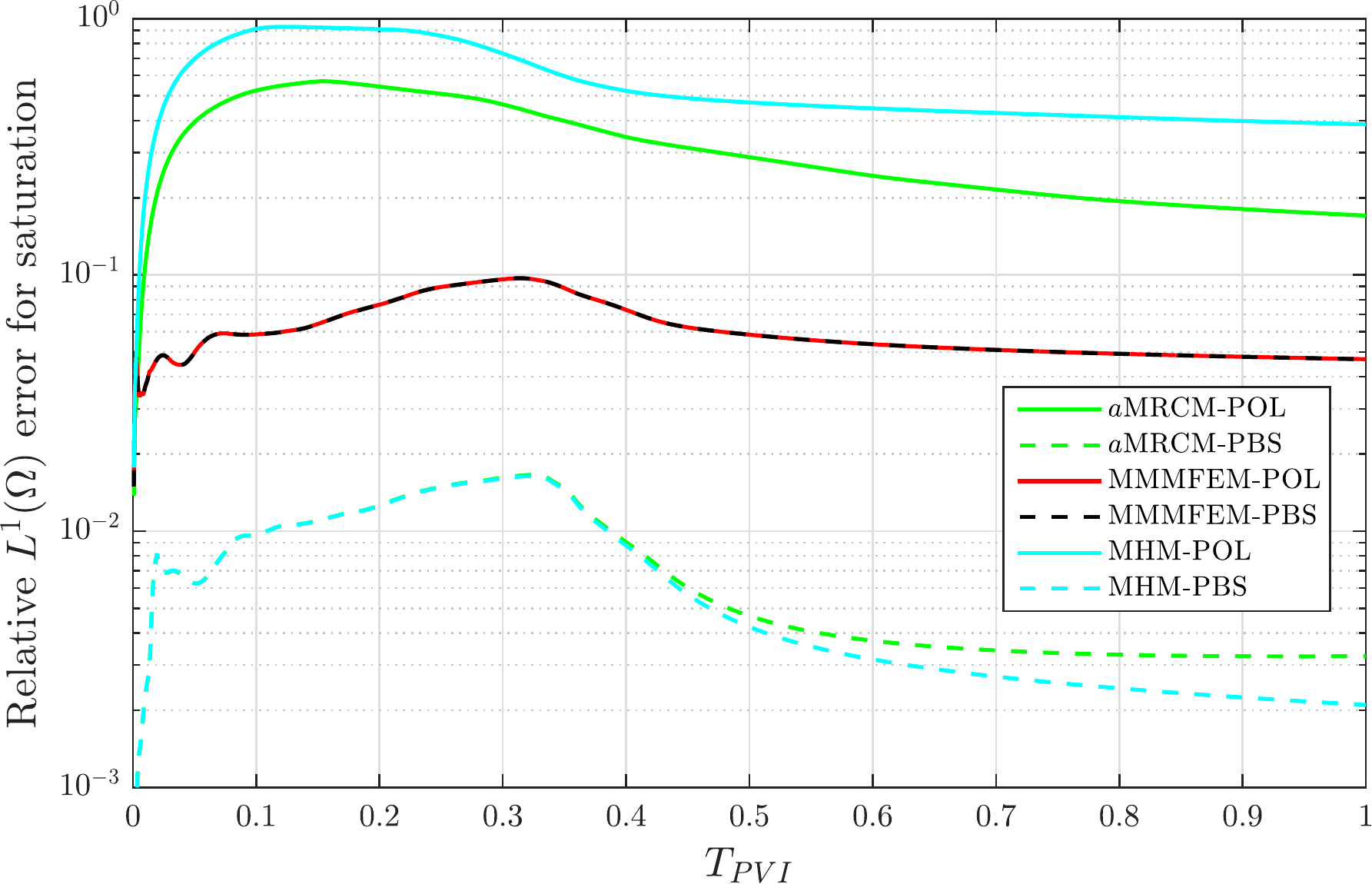} 
       \caption{Relative $L^{1}(\Omega)$ saturation errors as a function of time on the field with barriers. We compare the $a$MRCM, MMMFEM and MHM with both physics-based and linear spaces. The $a$MRCM-PBS and MMMFEM-PBS produce the smallest errors.}
        \label{fig:sat_low_erro_efendiev}
\end{figure}

As a final validation, we consider the performance of the methods in dealing with the same high-contrast permeability field replacing the type of channelized structures to barriers, as in \cite{jiang2017model}.
Figure \ref{fig:sat_low_efendiev} shows the permeability field (log-scaled) containing the low-permeable channels and the saturation profiles at $T_{\text{PVI}}=0.33$ (before breakthrough). In this case, only the flux physics-based spaces are used, since there are only low-permeable structures.
We note a considerable improvement in the $a$MRCM and MHM approximations replacing the linear interface spaces by the physics-based ones. The MMMFEM solutions present just some modest errors if compared to the $a$MRCM-POL and MHM-POL approximations. We remark that the MMMFEM solutions are the same with the linear and the physics-based interface spaces because, this method considers only the pressure space, which is always maintained as linear for this permeability field.
The relative errors throughout the whole simulation are shown in Figure \ref{fig:sat_low_erro_efendiev}. We confirm that the $a$MRCM-POL and MHM-POL are not accurate, the MMMFEM provides intermediate results and the $a$MRCM-PBS and MHM-PBS produce the most accurate approximations.

The numerical studies demonstrate the improvement obtained by using physics-based interface spaces. High-contrast fields were chosen that allowed us to assess the ability of the methods to handle problems in the presence of both high and low-permeable structures.
The results provide strong evidence that the adaptive MRCM combined with the physics-based spaces leads to improved transport approximations in high-contrast fractured-like porous media.

\section{Conclusions}\label{conclusions}

Two physics-based interface spaces (one for pressure and other for flux) have been proposed for better capturing the high-contrast effects of channelized structures. A careful investigation has been performed for the numerical solution of single and two-phase flows by combining the new spaces with multiscale mixed methods.

We show that the introduced physics-based pressure space can offer much better accuracy in comparison with the usual polynomial spaces in the presence of high-permeable structures. On the other hand, the proposed physics-based flux space can provide more accurate solutions in comparison with the polynomial spaces in the presence of low-permeable structures.
Two well known multiscale procedures have been applied to confirm these results: the Multiscale Mortar Mixed Finite Element Method (combined with the pressure space) and the Multiscale Hybrid-Mixed Finite Element Method (using the flux space). Aiming at using simultaneously both interface spaces we combine them with the MRCM, which allows for including the interface spaces independently. With this combination, we achieve the best accuracy in the approximation of challenging high-contrast flows in comparison with the other multiscale methods tested.

The MRCM formulation can take advantage of parallel computations with a computational cost comparable to existing procedures and providing a superior accuracy of the solutions in challenging high-contrast fields than the other multiscale methods that we have considered. 
The development of interface spaces based on physics for 3D reservoir flow problems
is outside the scope of the work presented here. However, this is a topic that is 
currently being considered by the authors.
Future work may include the use of the MRCM in more complex flow models, including the implicit solution of multiphase flow and transport problems.


\section*{Acknowledgments}
The authors gratefully acknowledge the financial support received from the Brazilian oil company Petrobras grant
2015/00400-4, and from the S\~ao Paulo Research Foundation FAPESP grant CEPID-CeMEAI 2013/07375-0; This study was
also financed in part by Brazilian government agencies CAPES (Finance Code 001) and CNPq grants 305599/2017-8 and 310990/2019-0; FP was
also funded in part by NSF-DMS 1514808, a Science Without Borders/CNPq-Brazil grant 400169/2014-2 and UT Dallas; FFR
acknowledges the hospitality provided by UT Dallas; Finally, the authors thank R.T. Guiraldello for providing the MRCM code and for
fruitful discussions about the method.

  \bibliographystyle{elsarticle-num} 

\bibliography{MRC_biblo,newRobertobibdata,references}

%
%
%
\end{document}